%% file: stokes_union.tex
\documentclass[1p,preprint]{elsarticle} 

\begin{document}

 \title{Stabilized isogeometric formulation of the Stokes problem on overlapping patches.
 } 

\author[1]{Xiaodong Wei
}

\author[2]{Riccardo Puppi
	} 

\author[2]{Pablo Antolin
}

\author[2]{Annalisa Buffa
}


\affiliation[1]{organization={Shanghai Jiao Tong University - University of Michigan Joint Institute}, addressline={Shanghai Jiao Tong University},
	postcode={200240}, city={Shangai}, country={China}}	
	
\affiliation[2]{organization={Institute of Mathematics}, addressline={École Polytechnique Fédérale de Lausanne},	postcode={1015}, city={Lausanne}, country={Switzerland}}

\begin{abstract}
We present a novel stabilized isogeometric formulation for the Stokes problem, where the geometry of interest is obtained via overlapping NURBS (non-uniform rational B-spline) patches, i.e., one patch on top of another in an arbitrary but predefined hierarchical order. All the visible regions constitute the computational domain, whereas independent patches are coupled through visible interfaces using Nitsche’s formulation. Such a geometric representation inevitably involves trimming, which may yield trimmed elements of extremely small measures (referred to as bad elements) and thus lead to the instability issue. Motivated by the minimal stabilization method that rigorously guarantees stability for trimmed geometries~\cite{ref:puppi20}, in this work we generalize it to the Stokes problem on overlapping patches. Central to our method is the distinct treatments for the pressure and velocity spaces: Stabilization for velocity is carried out for the flux terms on interfaces, whereas pressure is stabilized in all the bad elements. We provide a priori error estimates with a comprehensive theoretical study. Through a suite of numerical tests, we first show that optimal convergence rates are achieved, which consistently agrees with our theoretical findings. Second, we show that the accuracy of pressure is significantly improved by several orders using the proposed stabilization method, compared to the results without stabilization. Finally, we also demonstrate the flexibility and efficiency of the proposed method in capturing local features in the solution field.

This contribution is dedicated to Thomas J.R. Hughes, as a tribute to his remarkable lifetime achievements.
\end{abstract}

 \begin{keyword}
	Boolean operations \sep Stokes problem \sep minimal stabilization \sep boundary-unfitted method \sep Nitsche's method
\end{keyword}

\maketitle	

\input{introduction.tex}

\input{numerical_analysis.tex}

\input{example.tex}
\appendix
\input{appendix.tex}

\clearpage
\bibliographystyle{elsarticle-num}
\bibliography{bibliographyx}
\end{document}

%% file: introduction.tex
\section{Introduction}	
Isogeometric analysis (IGA) was introduced to bridge the gap between computer-aided design (CAD) and computer-aided engineering (CAE)~\cite{ref:hughes05}. Its central idea lies in employing the same basis functions for both design and analysis. IGA typically uses non-uniform rational B-splines (NURBS), the industrial standard of CAD, as basis functions to approximate partial differential equations (PDEs) that govern the physics of interest. In addition to the advancement in the CAD-CAE integration~\cite{ref:hsu15, ref:lai17, ref:casquero20}, IGA has also shown great potential in a wide range of applications~\cite{ref:kamensky15, ref:ali19} due to its superior numerical performance in terms of accuracy~\cite{ref:evans09} and robustness~\cite{ref:lipton10}.

However, directly using CAD models in analysis remains a challenge in IGA due to the ubiquitous trimming operation. NURBS alone is topologically restrictive and can only model box-like geometries. Trimming, essentially a visualization trick that makes part of a geometric model invisible to users, enables modeling geometries of arbitrary topologies in a practically handy way. A CAD model is generally a collection of trimmed NURBS surface patches. In the case of solid modeling, it only represents the boundary of a solid (B-rep) \rv{and thus it has no description of its interior domain. On the other hand, the description of the interior domain is often essential for physics-based simulations. There exist several possible ways to bridge this gap. For instance, volumetric parameterization seeks to find a watertight volumetric spline representation conformal to a given B-rep~\cite{ref:martin09, ref:zhang12, ref:zhang13}, but it usually involves generation of hexahedral meshes, which has been an challenge for decades. Another possibility is through the so-called volume representation (V-rep)~\cite{ref:massarwi16, ref:antolin19}, where the building blocks are true volumes and naturally have a description of the interior domain, and in the meantime its boundary is a B-rep. The development of analysis-aware and robust volumetric Boolean operations is the key to the success of this road. The third option, which is related to V-reps, is to embed a B-rep to an auxiliary background mesh.} Such a geometric representation falls into the category of immersed/unfitted discretizations (e.g., Finite Cell Method~\cite{ref:parvizian07}, cutFEM~\cite{ref:burman15}, immersogeometric analysis~\cite{ref:kamensky15}, shifted boundary method~\cite{main2018}, aggregated unfitted FEM~\cite{badia2018}), where computational meshes do not align with geometric boundaries/interfaces. While geometric modeling is extremely flexible, various issues have to be addressed on the analysis side, such as stability~\cite{ref:burman15, ref:puppi20, badia2022, dePrenter2023}, quadrature~\cite{kudela2016, ref:divi20, gunderman2021, garhuom2022, saye2022, ref:antolin22, antolin2022}, imposition of boundary conditions~\cite{ref:ruess13, ref:wei21b}, conditioning~\cite{dePrenter2023}, etc. Various methods and applications using trimmed NURBS~\cite{marussig2018} have been investigated in shells~\cite{ref:teschemacher18, ref:leidinger19a} and solid modeling~\cite{ref:schillinger12}. Interested readers may refer to~\cite{dePrenter2023} for a dedicated review on the stability and conditioning issues as well as their treatments.

In this paper, we adopt a general way for geometric construction through overlapping NURBS patches~\cite{ref:wei21a} that involves trimming, on which we study the Stokes problem and propose a novel stabilized isogeometric method. The proposed method is a generalization of the minimal stabilization for elliptic problems~\cite{ref:puppi20, ref:wei21a}, but it needs distinct treatments for the velocity space and the pressure space. Minimal stabilization is a theoretical framework to deal with the instability issue caused by the cut elements, which are inevitable in trimmed geometries. It was further extended with a conformal layer around boundaries/interfaces to greatly enhance the flexibility of handling boundary/interface conditions~\cite{ref:wei21b}. While recently minimal stabilization has been studied for the Stokes problem in a single trimmed patch~\cite{puppi_stokes}, this paper is focused on multiple overlapping patches, where handling interfaces with a stable formulation plays the key role.

Boundary-unfitted discretizations for the Stokes problem have also been extensively studied in the context of finite element methods. For instances, the inf-sup stability was rigorously shown in~\cite{guzman2018} for triangular/tetrahedral meshes, rather than being taken as an assumption. Aggregated FEM was also extended to the Stokes problem~\cite{badia2018stokes}, where additional stabilization terms are added to the mixed finite element spaces to fix the potential deficiencies of the aggregated inf-sup. Recently, high-order FEM was devised and analyzed in~\cite{burman2021} for the Stokes interface problem. 

The most relevant FEM to our work is the multi-mesh FEM~\cite{ref:johansson19, ref:johansson20}, which employs an overlapping multi-mesh to discretize the computational domain. Both methods adopt Nitsche's method to impose the transmission condition on interfaces. The difference lies in the stabilization: The multi-mesh FEM needs to stabilize over an overlapping domain that consist of both hidden and visible meshes, whereas our method only deals with visible meshes. More specifically, the multi-mesh FEM first finds a mesh intersection for the overlapping domain, which is relatively easy to compute for low-order FEM meshes but becomes much more involving for high-order IGA meshes; see~\cite{ref:cirillo23} for recent progress. In addition, an least-squares term that involves second-order derivatives is also added to the Nitsche's formulation for the stabilization purpose. On the other hand, our method only counts visible elements/interfaces as the computational domain. It retains the original form of Nitsche's method and only locally modifies certain necessary terms.

The core of the proposed work is to deal with two types of instabilities that arise from the discretizations of the velocity field and the pressure field. Stabilization for velocity is carried out for the flux terms on the interfaces, whereas pressure is stabilized in all the \emph{bad elements}, i.e., elements with a small effective area ratio (up to a given threshold). We rigorously prove the stability of the stabilized formulation and further derive the optimal \emph{a priori} error estimates. In the end, we study several numerical examples to show the convergence and accuracy of the proposed method. Expected optimal convergence is achieved, which consistently agrees with our theoretical results. The accuracy of pressure is significantly improved by several orders compared to the case without stabilization. Furthermore, we also demonstrate the flexibility and efficiency of the proposed method in capturing local features in the solution field.

The paper is structured as follows. After introducing in Section~\ref{sec:notation_modelpb} some notations and the model problem, in Section~\ref{sec:iga_discretization} we provide the basic notions that are specific to IGA, describe the construction of overlapping domains through the union operation, and introduce the approximation spaces. In Section~\ref{sec:ncmp} we discretize the Stokes problem in a collection of coupled domains with Nitsche's method, and introduce our stabilization procedure together with its key properties and the stabilized discrete formulation. Then, we develop our stability analysis and prove that optimal \emph{a priori} error estimates hold. We provide several numerical experiments in Section~\ref{sec:example} to validate the theory. Section~\ref{sec:con} concludes the paper and suggests several future directions.
\IfStandalone
{
	\bibliographystyle{../elsarticle-num}
	\bibliography{../bibliographyx}
}{}

%% file: numerical_analysis.tex
\section{Notation and model problem}\label{sec:notation_modelpb}
Let us briefly introduce the definitions and notations which will be frequently employed in this paper. With a slight abuse of notation, we will use the same symbol $\abs{\cdot}$ to denote both the $d$-dimensional Lebesgue measure and the $(d-1)$-dimensional Hausdorff measure. Given $D\subset \R^d$ and $\Sigma$ a hypersurface of $\R^d$ or a subset of it, $\abs{D}$ and $\abs{\Sigma}$ denote the $d$-dimensional Lebesgue measure of $D$ and the $(d-1)$-dimensional Hausdorff measure of $\Sigma$, respectively. The symbol $\# \cdot$ denotes the cardinality of a set. Given $E\subset \R^d$, the notations $E^\circ$ and $\operatorname{int} E$ denote its interior.

A domain is an open, bounded, subset of $\R^d$, $d\in\{2,3\}$. A domain $D$ with boundary $\partial D$ is said to be Lipschitz if for every $x\in\partial D$ there exists a neighborhood $U$ of $x$ such that $U\cap\partial D$ is the graph of a Lipschitz function. In the following $D$ denotes a Lipschitz domain with boundary $\partial D$, and $\Sigma$ a Lipschitz continuous surface contained in $\partial D$. The unit outer normal on $\partial D$ is denoted by $\n$.

We will denote as $\mathbb Q_{r,s,t}$ the vector space of trivariate polynomials of degree at most $r$, $s$, and $t$ in the first, second, and third variables, respectively (analogously for the case $d=2$), $\mathbb P_u$ the vector space of trivariate polynomials of degree at most $u$. We may write $\mathbb Q_k$ instead of $\mathbb Q_{k,k}$ or $\mathbb Q_{k,k,k}$. We will often consider the restriction of a polynomial space to a given domain $D$ and write, for instance, $\mathbb Q_k (D)$ instead of $\restr{\mathbb Q_k }{D}$.

We denote by $L^2(D)$ the space of square integrable functions on the domain $D$, equipped with the usual norm $\norm{\cdot}_{L^2(D)}$. We denote by $L^2_0(D)$ the subspace of $L^2(D)$ of functions with zero average, where the average of $v\in L^2(D)$ is $\overline v:=\abs{D}^{-1}\int_D v$.

For a given $\varphi:D\to \R$ sufficiently regular, $\bm\alpha$ a multi-index with $\abs{\bm\alpha}:=\sum_{i=1}^d\alpha_1$, and $j\in\N$, we define $D^{\bm\alpha}\varphi:={\frac{\partial^{\abs{\bm\alpha}}\varphi}{\partial x_1^{\alpha_1}\dots\partial x_d^{\alpha_d}}}$ and $\partial_n^j \varphi:=\sum_{\abs{\bm\alpha}=j}D^{\bm\alpha} \varphi  \n ^{\bm\alpha}$, where $\n^{\bm\alpha}:=n_1^{\alpha_1}\dots n_d^{\alpha_d}$. We indicate by $H^k(D)$, for $k\in\N$, the standard Sobolev space of functions in $L^2(D)$ whose $k$-th order weak derivatives belong to $L^2(D)$, equipped with the norm $\norm{\varphi}^2_{H^k(D)}:=\sum_{\abs{\bm \eta}\le k} \norm{D^{\bm\eta}\varphi}^2_{L^2(D)}$. Sobolev spaces of fractional order $H^r(D)$, $r\in\R$, can be defined by interpolation techniques, see~\cite{MR2424078}.

The space $H^1_{0,\Sigma}(D)$ consists of functions in $H^1(D)$ with vanishing trace on $\Sigma$. We write $H^1_{0}(D)$ instead of $H^1_{0,\partial D}(D)$. 
Let $H^\frac{1}{2}(\partial D)$ be the range of the trace operator of functions in $H^1(D)$ and we define its restriction to $\Sigma$ as $H^{\frac{1}{2}}(\Sigma)$. Both $H^\frac{1}{2}(\partial D)$ and $H^{\frac{1}{2}}(\Sigma)$ can be endowed with an intrinsic norm, see~\cite{MR2328004}. The dual space of $H^{\frac{1}{2}}(\Sigma)$ is denoted by $H^{-\frac{1}{2}} (\Sigma)$. The duality pairing between $H^{\frac{1}{2}}(\Sigma)$ and $H^{-\frac{1}{2}}(\Sigma)$ will be denoted with a formal integral notation.
Bold letters will be used for the  spaces of vector valued functions.

$C$ will denote generic positive constants that may change with each occurrence throughout the paper but are always independent of the local mesh size, the position of the visible interfaces with respect to the meshes, and the number of patches, unless otherwise specified. 

Let $\Omega$ be a domain with Lipschitz boundary $\Gamma$, with unit outer normal $\n$, such that $\Gamma=\overline\Gamma_D\cup\overline\Gamma_N$, where $\Gamma_D$ and $\Gamma_N$ are non-empty, open, and disjoint. The Stokes equations are a linear system that can be derived as a simplification of the Navier-Stokes equations. They describe the flow of a fluid under incompressibility and slow motion regimes. Given the body force $\f\in \bm L^2({\Omega})$, the Dirichlet datum $\u_D\in\bm H^{\frac{1}{2}}(\Gamma_D)$ and the Neumann datum $\u_N\in \bm H^{-\frac{1}{2}}(\Gamma_N)$, we look for the \emph{velocity} $\u:{\Omega}\to\R^d$ and \emph{pressure} $p:{\Omega}\to\R$ such that
\begin{equation}
	\begin{aligned}\label{stokes}
		-\mu \bm\Delta \u + \nabla p = \f,  \qquad&\text{in}\;{\Omega}, \\
		\dive \u = 0,\qquad &\text{in}\;{\Omega}, \\
		\u  = \u_D,\qquad &\text{on}\;\Gamma_D,\\
		\bm\sigma(\u,p) \n   = \u_N,\qquad &\text{on}\;\Gamma_N,
	\end{aligned}
\end{equation}
where $\mu>0$ is the \emph{viscosity coefficient}, $\bm\sigma(\u,p):= \mu D\u - p \mathbf{I}$ is the \emph{Cauchy stress tensor}, $\left(D\u\right)_{ij}:=\displaystyle{\frac{\partial u_i}{\partial x_j}}$, $i,j=1,\dots,d$. The first equation is known as the \emph{conservation of the momentum} and is nothing else than Newton's Second Law, relating the external forces acting on the fluid to the rate of change of its momentum, the second one is the \emph{conservation of mass}. For the sake of simplicity of the notation, let us set $\mu\equiv 1$.	
\section{Isogeometric analysis on overlapping multipach domains}\label{sec:iga_discretization}
Before explaining what we mean by ``overlapping multipach domains'', let us briefly introduce the basic notions in IGA.
\subsection{Univariate splines}
Given two positive integers $k$ and $n$, we say that $\Xi:=\{\xi_1,\dots,\xi_{n+k+1} \}$ is a $k$-\emph{open knot vector} if
\begin{equation*}
	\xi_1=\dots=\xi_{k+1}<\xi_{k+2}\le\dots \le \xi_n<\xi_{n+1}=\dots=\xi_{n+k+1}.
\end{equation*}
For the sake of simplicity, we assume $\xi_1=0$ and $\xi_{n+k+1}=1$. We also introduce $Z:=\{\zeta_1,\dots,\zeta_M \}$, the set of \emph{breakpoints}, or knots without repetitions, which forms a partition of the unit interval $(0,1)$. Note that
$$
\Xi=\{\underbrace{\zeta_1,\dots,\zeta_1}_{m_1\;\text{times}},\underbrace{\zeta_2,\dots,\zeta_2}_{m_2\;\text{times}},\dots,\underbrace{\zeta_M,\dots,\zeta_M }_{m_M\;\text{times}}\},
$$
where $m_j$ is the multiplicity of the breakpoint $\zeta_j$ and $\sum_{i=1}^Mm_i=n+k+1$. Moreover, we assume $m_j\le k$ for every internal knot, and we denote $I_i:=(\zeta_i,\zeta_{i+1})$ and its measure $h_i:=\zeta_{i+1}-\zeta_{i}$, $i=1,\dots, M-1$.

We denote as $\hat{B}_{i,k}:[0,1]\to\R$ the $i$th \emph{B-spline}, $1\le i\le n$, obtained using the \emph{Cox-de Boor formula}
\begin{align*}
	\hat B_{i,0}(\zeta):=&
	\begin{cases}
		1,\qquad&\text{if}\ \zeta \in [\xi_i,\xi_{i+1}),	\\
		0,\quad &\text{otherwise},
	\end{cases}	\\
	\hat B_{i,k}(\zeta) := &\frac{\zeta-\xi_{i}}{\xi_{i+k}-\xi_i}\hat B_{i,k-1}(\zeta) + \frac{\xi_{i+k+1}-\zeta}{\xi_{i+k+1}-\xi_{i+1}} \hat B_{i+1,k-1}(\zeta),\qquad k\ge 1,
\end{align*}
with the convention that $\displaystyle \frac{0}{0}=0$.
Moreover, let $S^k_{\bm\alpha}(\Xi):=\operatorname{span}\{\hat{B}_{i,k}: 1\le i\le n \}$ be the vector space of univariate splines of degree $k$. $S^k_{\bm\alpha}(\Xi)$ can also be characterized as the space of piecewise polynomials of degree $k$ with $\alpha_j:=k-m_j$ continuous derivatives at the breakpoints $\zeta_j$, $1\le j\le M$ (\emph{Curry-Schoenberg Theorem}~\cite{spline_draft}). The number of continuous derivatives at the breakpoints is collected in the \emph{regularity vector} $\bm\alpha:=\left(\alpha_j\right)_{j=1}^M$. A knot multiplicity $m_j=k+1$ corresponds to a regularity $\alpha_j=-1$, i.e., a discontinuity at the breakpoint $\zeta_j$. Since the knot vector is open, it holds $\alpha_1 = \alpha_M = -1$. For the sake of simplicity of the notation we assume that the basis functions have the same regularity at the internal knots, namely $\alpha_j = \alpha$ for $2 \le j \le  M-1$. Note that the derivatives of splines are splines too when $k \ge 1$ and $ \alpha \ge 0$ and, for $\Xi' := \{\xi_2,\dots.,\xi_{n+k}\}$, the operator
$\frac{d}{dx} : S^k_{\bm\alpha}(\Xi) \to S^{k-1}_{\bm\alpha-1}(\Xi')$ is surjective., where $\bm\alpha-1$ denotes the regularity vector $\left( \alpha_j-1\right)_{j=1}^M$.

Moreover, given an interval $I_j=\left( \zeta_j,\zeta_{j+1}\right)=(\xi_i,\xi_{i+1})$, we define its \emph{support extension} $\tilde I_j$ as
\begin{equation*}
	\tilde I_j:=\operatorname{int}\bigcup\{\operatorname{supp}(\hat{B}_{\ell,k}): \operatorname{supp}(\hat{B}_{\ell,k})\cap I_j\ne\emptyset, 1\le \ell\le n \}=\left(\xi_{i-k},\xi_{i+k+1}\right).
\end{equation*}
\subsection{Multivariate splines}
Let $d\in\{2,3\}$ denote the space dimension and $M_\ell,n_\ell,k_\ell\in\N$,\\ $\Xi_\ell=\{\xi_{\ell,1},\dots,\xi_{\ell,n_\ell+k_\ell+1} \}$, $Z_\ell=\{\zeta_{\ell,1},\dots,\zeta_{\ell,M_\ell} \}$ be given, for every $1\le \ell\le d$. We set the degree vector $\mathbf{k}:=(k_1,\dots,k_d)$, the regularity vectors $\bm\alpha_\ell$, $1\le \ell\le d$, and $\mathbf{\Xi}:=\Xi_1\times\dots\times\Xi_d$. As in the univariate case, we assume that the same regularity holds at the internal knots for every parametric direction, hence we drop the bold font once for all and write $\alpha_\ell$, $1\le \ell \le d$. Note that the breakpoints of $Z_\ell$ form a Cartesian grid in the parametric domain $\hat{\Omega}$, namely the \emph{parametric B\'ezier mesh}
\begin{equation*}
	\hat{\mathcal{M}}_{h}:=\{Q_{\mathbf{j}}=I_{1,j_1}\times\dots\times I_{d,j_d}: I_{\ell,j_\ell}=(\zeta_{\ell,j_\ell},\zeta_{\ell,j_\ell+1}): 1\le j_\ell\le M_\ell-1 \},
\end{equation*}
where each $Q_{\mathbf{j}}$ is called a \emph{parametric B\'ezier element}, with $h_{Q_{\mathbf{j}}}:=\operatorname{diam}\left( Q_{\mathbf{j}}\right)$. Let $h:=\max\{h_Q: Q\in\hat{\mathcal M}_{h}\}$.
\begin{assumption}\label{preliminaries:shape_regularity}
	The family of meshes $\{\hat{\mathcal M}_{h}\}_h$ is assumed to be \emph{shape-regular}, that is, the ratio between the smallest edge of $Q\in\hat{\mathcal M}_{h}$ and its diameter $h_Q$ is uniformly bounded with respect to $Q$ and $h$.
\end{assumption}
\begin{remark}
	Shape-regularity implies that the mesh is \emph{locally quasi-uniform}, \emph{i.e.}, the ratio of the sizes of two neighboring elements is uniformly bounded (see~\cite{MR2250029}). Also note that it allows us to assign $h_Q$ as the unique size of the element, without the necessity of dealing with the length of its edges separately.
\end{remark}

Let $\mathbf{I}:=\{\mathbf{i}=(i_1,\dots,i_d): 1\le i_\ell\le n_\ell \}$ be a set of multi-indices. For each $\mathbf{i}=(i_1,\dots,i_d)$, we define the set of \emph{multivariate B-splines} $\{\hat{B}_{\mathbf{i},\k}({ \bm\zeta})=\hat{B}_{i_1,k_1}( \zeta_1)\dots\hat{B}_{i_d,k_d}( \zeta_d): \mathbf{i}\in\mathbf{I}  \}$.
The \emph{multivariate spline space} in $\hat{\Omega}:=(0,1)^d$ is defined as
$
S^{\k}_{\alpha_1,\dots,\alpha_d}(\hat{\mathcal{M}}_h):=\operatorname{span}\{\hat{B}_{\mathbf{i},\k}:\mathbf{i}\in\mathbf{I} \},
$
which can also be seen as the space of piecewise multivariate polynomials of degree $\k$ and regularity $\alpha_1,\dots,\alpha_d$. Note that $S^{\k}_{\alpha_1,\dots,\alpha_d}(\hat{\mathcal{M}}_h)=\bigotimes_{\ell=1}^d S^{k_\ell}_{\alpha_\ell}(\Xi_\ell)$. 

Moreover, for an arbitrary B\'ezier element $Q_{\mathbf j}\in\hat {\mathcal M}_{h}$, we define its \emph{support extension} $\tilde Q_{\mathbf j}:=\tilde I_{1,j_1}\times\dots\times\tilde I_{d,j_d}$,
where $\tilde I_{\ell,j_\ell}$ is the univariate support extension of the univariate case defined above. 
\begin{remark}
The previous constructions can be easily generalized to the case of Non-Uniform Rational B-Splines (NURBS) basis functions; see, for instance,~\cite{MR3618875}. For the sake of simplicity of the exposition, in this manuscript we restrict ourselves to the case of B-splines bases.
\end{remark}	
\subsection{Construction of the union domain}
Let $\Omega_i^\ast \subset\R^d$, $0\le i\le N$, with $N\in\N$, $d\in\{2,3\}$, be spline patches, \emph{i.e.}, $\Omega_i^\ast=\F_i(\hat\Omega)$, where $\mathbf{F}_i\in \left[ S^{\mathbf{k}}_{\alpha_1,\dots,\alpha_d}(\hat{\mathcal{M}}_h) \right]^d$ is the geometric mapping defined by multivariate B-splines and corresponding control points. Each patch has an underlying \emph{premesh} $\mathcal{M}_i^{\ast}$, naturally induced by the map $\F_i$, namely
\begin{equation*}
	\mathcal{M}_i^{\ast}:=\{K\subset\Omega_i^{\ast}: K=\mathbf{F}_i(Q), Q\in\hat{\mathcal{M}}_h \}.
\end{equation*}

For the sake of simplicity of the notation and the analysis, the following simplifications are made.	
\begin{assumption}\label{assumption_isotropic_degrees}
	We assume that the degree-vector is \emph{isotropic} and that all predomains are parametrized by splines of the same degree, \emph{i.e.}, $\mathbf{k}=\left(k,\dots, k \right)$. Hence, we may write $k$ instead of $\k$. Similarly for the regularity vectors, that is, $\alpha=\alpha_1=\dots=\alpha_d$ for all predomains.
\end{assumption}	
Moreover, to prevent the existence of singularities in the parametrizations, we make the following assumption.
\begin{assumption}\label{preliminaries:1}
	The parametrizations $\mathbf{F}_i:\hat{\Omega}\to\Omega_i$, $0
\le i\le N$, are bi-Lipschitz. Moreover, $\restr{\mathbf{F}_i}{\overline{Q}}\in C^\infty(\overline{Q})$ for every $Q\in\hat{\mathcal{M}}_{h}$ and $\restr{\mathbf{F}_i^{-1}}{\overline{K}}\in C^\infty(\overline{K})$ for every $K\in\mathcal{M}_i^\ast$, $0\le i\le N$.
\end{assumption}
Assumption~\ref{preliminaries:1} together with Assumption~\ref{preliminaries:shape_regularity} implies that the premeshes $\mathcal{M}_i^{\ast}$, $0\le i\le N$, are shape-regular too.

\begin{figure}[!ht]
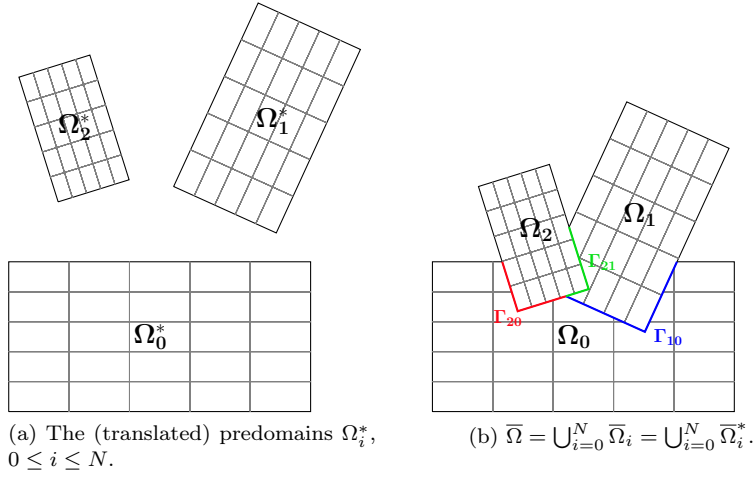

	\centering
	\subfloat[][The (translated) predomains $\Omega_i^\ast$, $0\le i \le N$.\label{figure1_a}]
	{\includestandalone[width=0.35\textwidth]{overlap1}
	}
	\hspace{0.5cm}
	\subfloat[][$\overline\Omega=\bigcup_{i=0}^N\overline\Omega_i=\bigcup_{i=0}^N\overline\Omega_i^\ast$.\label{figure1_b}]
	{
		\includestandalone[width=0.35\textwidth]{overlap2}
	}
	\caption{Definitions of predomains (a), visible parts of predomains, and local interfaces ($\Gamma_{10}$, $\Gamma_{20}$ and $\Gamma_{21}$) of predomain boundaries (b).}\label{figure1}
\end{figure}

Let us assume that our computational domain $\Omega$ is a \emph{union domain}, namely $\overline\Omega = \cup_{i=0}^N\overline\Omega_i^\ast$. We define $\Omega_i$ as the \emph{visible part} of the predomain $\Omega_i^\ast$
\begin{equation*}
	\Omega_i:= \Omega_i^\ast\setminus \bigcup_{\ell=i+1}^N\overline\Omega_\ell^\ast,\qquad i=0,\dots, N.
\end{equation*}
It holds $\Omega_N = \Omega_N^\ast$ and $\overline\Omega = \bigcup_{i=0}^N \overline\Omega_i=\bigcup_{i=0}^N \overline\Omega_i^\ast$, see Figure~\ref{figure1}. Note that this choice of definition of the $\Omega_i$'s follows~\cite{ref:wei21a, MR4154374} and implies a hierarchy of predomains. In particular, if $i>j$ then $\Omega_i^\ast$ is on top of $\Omega_j^\ast$, in the sense that $\Omega_i^\ast\cap\Omega_j^\ast$ is hidden by $\displaystyle\strut\cup_{k\ge i}\Omega_k$. We define 
\begin{equation*}
	\Gamma_i:=\partial\Omega_i^\ast\setminus \bigcup_{\ell=i+1}^N\overline \Omega_\ell^\ast,\qquad i=0,\dots, N,
\end{equation*}
\emph{i.e.}, the interface $\Gamma_i$ is the \emph{visible part of the external boundary} of $\Omega_i^\ast$ with outer unit normal $n_i$. Moreover, we define the \emph{local interfaces} $\Gamma_{ij}$ as
\begin{equation*}
	\Gamma_{ij}:=\Gamma_i\cap\overline \Omega_j,\qquad 0\le j<i\le N,
\end{equation*}
\emph{i.e.}, $\Gamma_{ij}$ is the subset of the visible boundary of $\Omega_i^\ast$ that intersects $\Omega_j$, see Figure~\ref{figure1}\subref{figure1_b}. We assume that each interface $\Gamma_{ij}$ either has non-zero $\left(d-1\right)$-measure or is the empty set. We also assume that $\Gamma_{ij}$ inherits the orientation of $\Gamma_i$, hence it has outer unit normal $\n_i$, also denoted as $\n$ when it is clear from the context {which domain is referred to}. Note that $\Gamma_{ij}$ is not connected in general. 

Finally, let us make a mild assumption on the roughness of the interfaces $\Gamma_{ij}$.
\begin{assumption}\label{mesh_assumption}
	All the interfaces $\Gamma_{ij}$, $1\le j< i \le N$, are Lipschitz-regular.
\end{assumption}	
\begin{lemma}\label{lemma1}
	There exists $C>0$ such that, for every $1\le j< i \le N$ and $K\in	\mathcal M_j $, it holds $\abs{ \Gamma_{ij}\cap \overline K}\le C \restr{\mathsf h_j^{d-1}}{K}$.
\end{lemma}
\begin{proof}
	This result holds since $\Gamma_{ij}$ is assumed to be Lipschitz-regular, hence not too oscillating; see~\cite{MR3806650}.
\end{proof}
Note that in what follows, we are going to refer to elements $K\in\mathcal M_i$, $i=0,\dots, N$, such that $\abs{\Omega_i\cap K} < \abs{K}$ as \emph{cut elements}. Moreover, integrals and norms will be defined on sets like $\Gamma_{ij}$, $\Gamma_{ij}\cap\overline K$ and they are meant to be on their interior in a suitable sense.
\begin{definition}\label{new_constant_1}
	Let $0\le j <i \le N$, and
	\begin{equation*}
		\delta_{ij}:=
		\begin{cases}
			1 \qquad\text{if}\;\;  \Gamma_{ij}\ne\emptyset, \\
			0 \qquad \text{otherwise}.
		\end{cases} 
	\end{equation*}
	For each $1\le i\le N$, the quantity $\sum_{j=0}^{i-1}\delta_{ij}$ counts the number of visible parts $\Omega_j$ whose boundaries are overlapped by $\Gamma_i$ and, for $0\le j\le N-1$, $\sum_{i=j+1}^N\delta_{ij}$ the number of visible parts whose boundaries overlap $\Gamma_j$. We further define $N^\downarrow_{\Gamma}:=\max_{1\le i \le N}\sum_{j=0}^{i-1}\delta_{ij}$, the maximum number of visible parts whose boundaries are overlapped by any visible boundary, and $N^\uparrow_{\Gamma}:=\max_{0\le j \le N-1}\sum_{i=j+1}^{N}\delta_{ij}$, the maximum number of visible parts whose boundaries cover any visible boundary. Hence, we let $N_\Gamma:=\max\{N^\downarrow_{\Gamma},N^\uparrow_{\Gamma}\}$ be the \emph{maximum number of boundary overlaps} in the current configuration.
\end{definition}
\begin{definition}\label{new_constant_2}
	We let $O_{ij}:=\Omega_i\cap \Omega_j^\ast$, $0\le j < i \le N$, be the \emph{overlap} between the $j$-th predomain and the $i$-th visible part. For every $0\le j<i\le N$, we define
	\begin{equation*}
		\eta_{ij}:=
		\begin{cases}
			1 \qquad\text{if}\;\;  O_{ij}\ne\emptyset, \\
			0 \qquad \text{otherwise}.
		\end{cases} 
	\end{equation*}
	For each $1\le i\le N$, the quantity $\sum_{j=0}^{i-1}\eta_{ij}$ counts the number of predomains covered by the visible part $\Omega_i$. We further define $N_{\mathcal O}:=\max_{1\le i \le N}\sum_{j=0}^{i-1}\eta_{ij}$, the maximum number of predomains covered by any visible part.
\end{definition}

\begin{figure}[!ht]
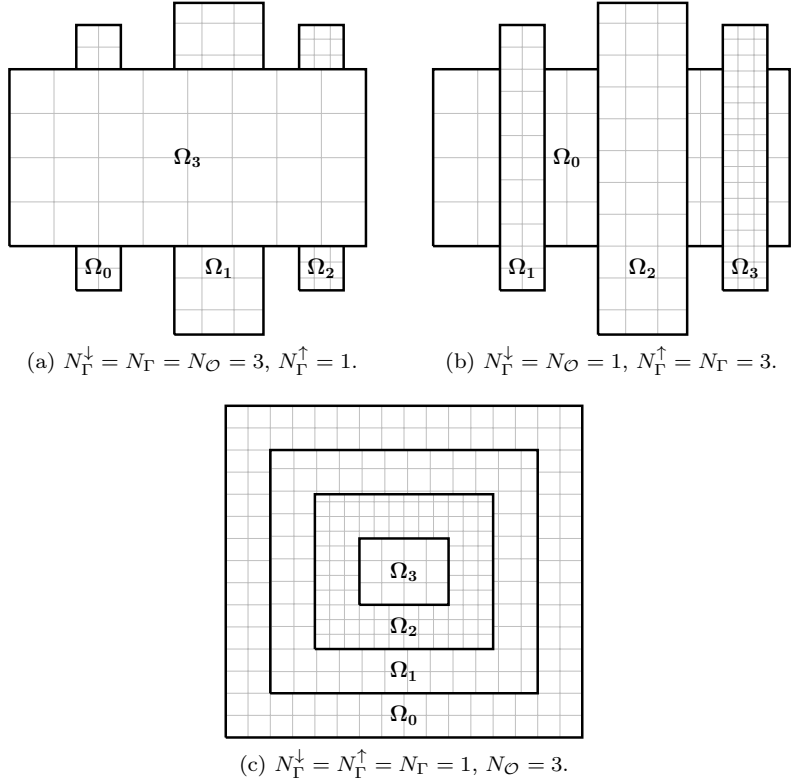

	\centering
	\subfloat[][$N_\Gamma^\downarrow=N_\Gamma=N_{\mathcal O}=3$, $N_\Gamma^\uparrow=1$.\label{figure2_a}]
	{\includestandalone[width=0.35\textwidth]{N_gamma_down}
	}
	\hspace{0.5cm}
	\subfloat[][$N_\Gamma^\downarrow=N_{\mathcal O}=1$, $N_\Gamma^\uparrow=N_\Gamma=3$.\label{figure2_b}]
	{
		\includestandalone[width=0.35\textwidth]{N_gamma_up}
	}\\
	\subfloat[][$N_\Gamma^\downarrow=N_\Gamma^\uparrow=N_\Gamma=1$, $N_{\mathcal O}=3$.\label{figure2_c}]
	{
		\includestandalone[width=0.35\textwidth]{N_O}
	}
	\caption{Illustration of $N_\Gamma^\downarrow$, $N_\Gamma^\uparrow$, $N_\Gamma$, $N_{\mathcal O}$.}\label{figure2}
\end{figure}

We observe that $N_\Gamma, N_{\mathcal O}\le N$; see Figure~\ref{figure2}. Moreover, in applications we expect $N_{\Gamma}, N_{\mathcal O}\ll N$.

We refer to $\mathcal M_i :=\{ K\in \mathcal M^\ast_i: K\cap\Omega_i\ne\emptyset  \}$, $ i=0,\dots,N,$
as the $i$-th \emph{active mesh}, consisting of all visible elements (not necessarily fully visible) of the $i$-th premesh $\mathcal M^\ast_i$. We define $\mathsf h_i:\Omega_i\to\mathbb{R}^+$ to be the piecewise constant mesh size function of $\mathcal M_i$ assigning to each visible element its whole diameter (rather than the diameter of its visible part), namely $\restr{\mathsf h_i}{K\cap\Omega}:=h_{i,K}$, where $h_{i,K}:=\operatorname{diam}\left(K\right)$ for every $K\in\mathcal M_i$, $0\le i\le N$. Moreover let us denote $h_i:=\max_{K\in\mathcal M_i}h_{i,K}$ and $h:=\max_{0\le i \le N}h_i$ (we use the same notation as for the mesh size of the parametric domain; this is not a dangerous ambiguity since Assumption~\ref{preliminaries:1} implies that the local mesh sizes are comparable to the one of the parametric domain). Finally, we denote by $\h:\Omega\to\R^{+}$ the piecewise constant function defined as $\restr{\h}{\Omega_i}:=\h_i$. Throughout the paper, we are going to rely {on the} assumption that adjacent sub-domains are discretized with similar mesh sizes.
\begin{assumption}\label{assumption_uniformity}
	The meshes locally have compatible sizes in the following sense. There exist $c,C>0$ such that for every $\Gamma_{ij}\ne\emptyset$, with $1\le j< i \le N$, $K_i\in	\mathcal M_i$ such that  $\overline K_i\cap\Gamma_{ij}\ne\emptyset$ and $ K_j\in	\mathcal M_j$ such that $\overline K_j\cap\Gamma_{ij}\ne\emptyset$, it holds
	\begin{equation*}
		c  \restr{\mathsf h_j}{K_j}\le \restr{\mathsf h_i}{K_i} \le C \restr{\mathsf h_j}{K_j}.
	\end{equation*}		
\end{assumption}	

\subsection{Construction of the local approximation spaces}
In this section, we construct the isogeometric spaces in the parametric domain.
\begin{align*}
	\hat V_{h}^\text{RT}:=&
	\begin{cases}
		S^{k+1,k}_{ \alpha+1,\alpha}(\hat{\mathcal{M}}_h)\times S^{k,k+1}_{\alpha,\alpha+1}(\hat{\mathcal{M}}_h),&\text{if}\;d=2,\\
		S^{k+1,k,k}_{\alpha+1,\alpha, \alpha}(\hat{\mathcal{M}}_h)\times S^{k,k+1,k}_{\alpha,\alpha+1,\alpha}(\hat{\mathcal{M}}_h)\times S^{k,k,k+1}_{\alpha,\alpha,\alpha+1}(\hat{\mathcal{M}}_h),\;\;\;\;\;\qquad\qquad\qquad&\text{if}\;d=3,
	\end{cases}\\
	\hat V_{h}^\text{N}:=&
	\begin{cases}
		S^{k+1,k+1}_{\alpha+1,\alpha}(\hat{\mathcal{M}}_h)\times S^{k+1,k+1}_{\alpha,\alpha+1}(\hat{\mathcal{M}}_h),&\text{if}\;d=2,\\
		S^{k+1,k+1,k+1}_{\alpha+1,\alpha,\alpha}(\hat{\mathcal{M}}_h)\times S^{k+1,k+1,k+1}_{\alpha,\alpha+1,\alpha}(\hat{\mathcal{M}}_h)\times S^{k+1,k+1,k+1}_{\alpha,\alpha,\alpha+1}(\hat{\mathcal{M}}_h),&\text{if}\;d=3,
	\end{cases}\\
	\hat V_{h}^\text{TH}:=&
	\begin{cases}
		S^{k+1,k+1}_{\alpha,\alpha}(\hat{\mathcal{M}}_h)\times S^{k+1,k+1}_{\alpha,\alpha}(\hat{\mathcal{M}}_h),&\text{if}\;d=2,\\
		S^{k+1,k+1,k+1}_{\alpha,\alpha,\alpha}(\hat{\mathcal{M}}_h)\times S^{k+1,k+1,k+1}_{\alpha,\alpha,\alpha}(\hat{\mathcal{M}}_h)\times S^{k+1,k+1,k+1}_{\alpha,\alpha,\alpha}(\hat{\mathcal{M}}_h),&\text{if}\;d=3,
	\end{cases}\\
	\hat Q_{h}:=&
	\begin{cases}
		S^{k,k}_{\alpha,\alpha}(\hat{\mathcal{M}}_h), \qquad&\text{if}\;d=2,\\
		S^{k,k,k}_{\alpha,\alpha,\alpha}(\hat{\mathcal{M}}_h), \qquad&\text{if}\;d=3.
	\end{cases}
\end{align*}
Note that $\hat{\mathcal{M}}_h$ depends on $d$. Moreover, spaces of different degrees and regularities on the same parametric mesh $\hat{\mathcal{M}}_h$ will correspond to different knot vectors. But we do not need to include these details here and we refer the reader to~\cite{MR3202239} for an in-depth discussion.
It holds that $\hat V_{h}^\text{RT}\subset \hat V_{h}^\text{N}\subset \hat V_{h}^\text{TH}$; see~\cite{MR2808112}. We note that for $\alpha=-1$, $\hat V^\text{RT}_{h}$ and $\hat V^\text{N}_{h}$ recover the classical Raviart-Thomas finite element and N\'ed\'elec finite element of the second kind, respectively. For $\alpha=0$, $\hat V_{h}^\text{TH}$ represents the classical Taylor-Hood finite element space. Henceforth we assume $\alpha\ge 0$, and adopt the notation $\hat V_{h}^\square$, $\square\in\{\text{RT},\text{N},\text{TH}\}$.

 For every $0\le i\le N$, we define the local approximation spaces in the patches 
\begin{align*}
	V_{h,i}^{\ast,\square}:=&  \{ \vv_h : \iota^i_v (\vv_h) \in \hat V_{h}^\square   \}, \qquad \square\in\{ \text{RT},\text{N}\},\qquad 
	 V_{h,i}^{\ast,\mathrm{TH}} :=  \{ \vv_h: \vv_h\circ \F_i \in \hat V_{h}^\mathrm{TH}   \},\\
	Q^{\ast,\square}_{h,i}:=&	\{q_h:  \iota_p^i(q_h) \in \hat Q_{h}\},\quad 	Q^{\ast,\mathrm{TH}}_{h,i}:=  \{q_h: q_h\circ \F_i \in \hat Q_{h}\},
\end{align*}
where 
\begin{alignat*}{3}
	\iota_v^i:&\bm H(\dive;\Omega_i^\ast)\to \bm H(\dive;\hat\Omega),\qquad  &&\iota^i_v \left(\vv\right):=\det \left( D\F_i\right) D\F_i^{-1} \left( \vv\circ\F_i\right),\\
	\iota_p^i: & L^2(\Omega_i^\ast)\to L^2(\hat \Omega),\qquad &&\iota^i_p(q):=\det \left( D\F_i\right)\left( q\circ\F_i\right).
\end{alignat*}
\section{Isogeometric discretization on overlapping multipatch domains}\label{sec:ncmp}
\subsection{Nitsche's method for the Stokes problem}
We rewrite problem~\eqref{stokes} in the following multipatch form. Find $\u:\Omega\to\R^d$ and $p:\Omega\to\R$ such that
\begin{subequations}
	\begin{align}	
		-\dive\bm\sigma(\u_i,p_i) = \f, &\qquad\text{in}\; \Omega_i,\;  i=0,\dots,N,\label{problem2.1}\\
		-\dive \u_i = 0,&\qquad\text{in}\; \Omega_i,\;  i=0,\dots,N,\label{problem2.2}\\
		\u_i - \u_j = \bm 0,  &\qquad\text{on}\; \Gamma_{ij},\;  0\le j<i\le N,\label{problem2.3}\\
		\bm\sigma(\u_i,p_i)\n_i +\bm\sigma(\u_j,p_j)\n_j= \bm 0, &\qquad \text{on}\; \Gamma_{ij},\; 0\le j<i\le N,\label{problem2.4}\\
		\u_i=\u_D, &\qquad\text{on}\; \Gamma_D\cap \Gamma_i,\; i=0,\dots,N,\label{problem2.5}\\
		\bm\sigma(\u_i,p_i)\n_i=\u_N,&\qquad \text{on}\; \Gamma_N\cap\Gamma_i,\; i=0,\dots,N,\label{problem2.6}
	\end{align}
\end{subequations}
where $\u_i:=\restr{\u}{\Omega_i}$, $p_i:=\restr{p}{\Omega_i}$, $i=0,\dots, N$. Equations~\eqref{problem2.3} and~\eqref{problem2.4} are commonly known as~\emph{transmission conditions} at the local interfaces.
\begin{proposition}\label{prop1}
	Problems~\eqref{stokes} and~\eqref{problem2.1}--\eqref{problem2.6} are equivalent. 
\end{proposition}
\begin{proof}
	To demonstrate this result, the variational formulation of the two problems must be used. We refer the interested reader to Chapter~5 of~\cite{MR1857663}.
\end{proof}	

Let us introduce, for each visible part $\Omega_i$, the \emph{local isogeometric spaces}
\begin{align*}
	V_{h,i}^\square:=&  \{ \restr{\vv_h}{\Omega_i} : \vv_h \in  V_{h,i}^{\ast,\square}   \}, \quad
	Q_{h,i}^\square:=  \{\restr{q_h}{\Omega_i}: q_{h,i} \in Q^{\ast,\square}_{h,i}\},
\end{align*}
and glue them to form the \emph{union isogeometric spaces}
\begin{align*}
	V_h^\square:=\bigoplus_{i=0}^N  V^{\square}_{h,i},\qquad Q_h^\square:=\bigoplus_{i=0}^N  Q^\square_{h,i},
\end{align*}
where $\square\in\{\mathrm{RT},\mathrm{N},\mathrm{TH}\}$. To further alleviate the notation, we adopt the convention to denote as $V_h$ the space of the velocities and omit
the superscript $\square\in\{\mathrm{RT},\mathrm{N},\mathrm{TH}\}$ when what said does not depend from the particular finite element choice. Elements of $V_h$ and $Q_h$ are $(N+1)$-tuples $\vv_h=\left( \vv_0,\dots,\vv_N \right)$ and $q_h=\left( q_0,\dots,q_N \right)$, respectively. In practice, we can treat them as ordinary functions thanks to the identification,
\begin{align*}
	\vv_h(x) = \vv_i(x),\qquad  \forall x\in\Omega_i, \, i=0,\dots,N,\\
	q_h(x) = q_i(x),\qquad\  \forall x\in\Omega_i, \, i=0,\dots,N.
\end{align*} 

We assume for simplicity that $\partial\Omega$ contains the image of one or more full faces and that Dirichlet boundary conditions are set on full faces only.
\begin{assumption}\label{assumption_dirichlet_bc}
There exists $0\le i\le N$ such that $\hat \Gamma_i$, $\hat \Gamma_i:=\F_i^{-1}(\Gamma_i)$, contains a full face of the parametric domain $\hat\Omega$. Moreover, for every $0\le i\le N$, $\hat \Gamma_D := \F^{-1}(\Gamma_D)$, $\hat\Gamma_D\cap\hat\Gamma_i$ is either empty or the union of full faces of $\hat\Omega$.
\end{assumption}
In case Assumption~\ref{assumption_dirichlet_bc} does not hold, we can combine the technique in this paper with that detailed in~\cite{puppi_stokes} to deal with the imposition of Dirichlet boundary conditions in a weak sense.

Let $\varphi:\Omega\to\R$ be smooth enough and, for every $0\le i\le N$, we denote its restriction to $\Omega_i$ as $\varphi_i:=\restr{\varphi}{\Omega_i}$. Then, for every interface $\Gamma_{ij}$, $0\le j < i \le N$, and a.e. $x\in\Gamma_{ij}$, we define, respectively, the \emph{average} and the \emph{jump} of $\varphi$ as
\begin{align}
	\langle \varphi\rangle_{t,\Gamma_{ij}}(x) :=&t\restr{\varphi_i}{\Gamma_{ij}} (x)+(1-t) \restr{\varphi_j}{\Gamma_{ij}} (x),\qquad t\in\{\frac{1}{2},1\},\\	
	[\varphi]_{\Gamma_{ij}}(x) :=&\restr{\varphi_i}{\Gamma_{ij}} (x)- \restr{\varphi_j}{\Gamma_{ij}} (x).
\label{eq:flux}
\end{align}		
We may remove the subscript $\Gamma_{ij}$ when it is clear from the context {to which interface we refer}. The average is said to be \emph{symmetric} if $t=\frac{1}{2}$ and \emph{one-sided} when $t=1$. We define the jump and average of a vector valued function $\bm\tau:\Omega\to\R^d$ componentwise by letting $\langle \bm\tau \rangle_{t,k}:=\langle \tau_k\rangle_t$ and $[\bm\tau]_{k}:=[\tau_k]$, for $0\le k\le d$.

Let us endow $V_h$ and $Q_h$ with the mesh dependent norms:
\begin{align*}
	\norm{\vv_h}^2_{1,h}:=&\sum_{i=0}^N\norm{D \vv_i}^2_{L^2(\Omega_i)}+\sum^N_{i=1}\sum_{j=0}^{i-1}\norm{\h^{-\frac{1}{2}}\left[ \vv_h \right]}^2_{L^2(\Gamma_{ij})},\qquad\ \vv_h\in V_h,\\
	\norm{q_h}^2_{0,h}:=&\sum_{i=0}^N\norm{ q_i}^2_{L^2(\Omega_i)}+\sum^N_{i=1}\sum_{j=0}^{i-1}\norm{\h^{\frac{1}{2}}\left[ q_h \right]}^2_{L^2(\Gamma_{ij})},\qquad\qquad q_h\in Q_h.
\end{align*}
We further define
\begin{equation*}
V_h^{\u_D} = \{\vv_h \in V_h: \restr{\vv_h}{\Gamma_D} = \u_D \}, \qquad V_h^{\bm 0} = \{\vv_h \in V_h: \restr{\vv_h}{\Gamma_D} = \bm 0 \}.
\end{equation*}
By enforcing the transmission conditions in a weak sense using Nitsche's method, we obtain the following discrete problem.

Find $\left(\u_h,p_h\right)\in V_h^{\u_D}\times Q_h$ such that
\begin{equation}
	\begin{aligned}\label{non_stabilized_formulation}
		a_h(\u_h,\vv_h) + b(\vv_h,p_h) &= F(\vv_h),\qquad &\forall\ \vv_h\in V_h^{\bm 0},\\
		b(\u_h,q_h) &= 0,\qquad &\forall\ q_h\in Q_h,
	\end{aligned}
\end{equation}
where
\begin{equation*}
	\begin{aligned}
		a_h(\w_h,\vv_h):=&\sum_{i=0}^N\int_{\Omega_i}D \w_i: D \vv_i
		-\sum_{i=1}^N\sum_{j=0}^{i-1}\int_{\Gamma_{ij}}\left(\langle D\w_h\n \rangle_t\left[\vv_h\right]+\left[\w_h\right]\langle D\vv_h \n \rangle_t\right)\\
		&+\gamma\sum_{i=1}^N\sum_{j=0}^{i-1}\int_{\Gamma_{ij}}\h^{-1}\left[\w_h\right]\left[\vv_h\right], \qquad\w_h,\vv_h\in V_h,  \\
		b(\vv_h,q_h) : = & - \sum_{i=0}^N \int_{\Omega_i} q_i \dive \vv_i + \sum_{i=1}^N \sum_{j=0}^{i-1} \int_{\Gamma_{ij}} \langle q_h \rangle_t [\vv_h\cdot\n ],\qquad\ \vv_h\in V_h, q_h\in Q_h,
	\end{aligned}
\end{equation*}
with $t\in\{\frac{1}{2},1\}$, and
\begin{equation*}
	\begin{aligned}
		F(\vv_h):=&\sum_{i=0}^N\int_{\Omega_i}\f\cdot\vv_i + \int_{\Gamma_N} \u_N\cdot \vv_h,\qquad \vv_h\in V_h.
	\end{aligned}
\end{equation*}

\begin{proposition}\label{proposition:consistency}
	The discrete variational formulation in Equation~\eqref{non_stabilized_formulation} is consistent, \emph{i.e.}, the solution $(\u,p)$ of problem~\eqref{stokes} satisfies problem~\eqref{non_stabilized_formulation} as well. 
\end{proposition}
\begin{proof}
	The proof is quite standard. See, for instance, Chapter~5 of~\cite{MR1857663}.	
\end{proof}	

Unfortunately, the choice of the parameter $\gamma$ that guarantees well-posedness of~\eqref{non_stabilized_formulation} cannot be done independently of how meshes are cut. In the next section, we discuss this issue and propose a suitable stabilization approach.

\subsection{Stabilization procedure} \label{sec:ministab}
As it is well documented in the literature (see e.g., \cite{ref:johansson19, ref:puppi20, ref:wei21a}), when there are cut elements, the lack of stability for problem~\eqref{non_stabilized_formulation} depends on two main issues:
\begin{enumerate}[i)]
\item The Nitsche's method requires stabilization. We will adopt the minimal stabilization approach proposed in~\cite{ref:puppi20} and successfully used in~\cite{ref:wei21a}.
\item The formulation~\eqref{non_stabilized_formulation} is not inf-sup stable in presence of small and slim elements. This was observed in e.g., \cite{puppi_stokes}; see also~\cite{ref:johansson20}.
\end{enumerate}

In what follows, our approach to both i) and ii) will be based on the local extrapolation operators introduced in~\cite{ref:puppi20}. We start by introducing the concept of bad element.

\begin{figure}[!ht]
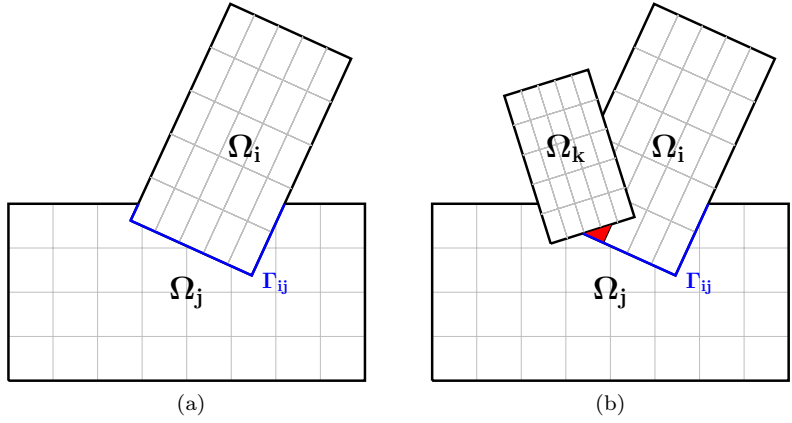

	\centering
	\subfloat[][\label{interface_overlap_a}]
	{\includestandalone[width=0.35\textwidth]{small_intersection_j_i_bis}
	}
	\hspace{0.5cm}
	\subfloat[][\label{interface_overlap_b}]
	{
		\includestandalone[width=0.35\textwidth]{small_intersection_j_i}
	}
	\caption{Two-patch and three-patch overlapping along the interface $\Gamma_{ij}$.}
	\label{interface_overlap_madre}
\end{figure}

 \begin{definition}\label{def_goodandbad}
	Fix $\theta\in (0,1]$, the \emph{area-ratio threshold}. For every $K\in{\mathcal{M}}_i$, $i=0,\dots, N$, we say that $K$ is a \emph{good element} if
	\begin{equation*}
		\frac{\abs{\Omega_i\cap K}}{\abs{K}}\ge\theta.
	\end{equation*}
	Otherwise, $K$ is a \emph{bad element}. 
	We denote as $\mathcal M^g_i$ and $\mathcal M^b_i$ the collection of good and bad physical B\'ezier elements, respectively. Note that all the uncut elements in $\mathcal M_i$ are good elements, and all its bad elements are cut elements. Moreover, it holds $\mathcal M^g_N =\mathcal M_N$ (or $\mathcal M_N^b=\emptyset$).
	
\end{definition}
\begin{definition}\label{def_neigh}
	Given $K\in\mathcal M_i$, $i=0,\dots,N-1$, the set of its \emph{neighbors} is
	\begin{equation} \label{constante_scuola}
		\mathcal N (K):=\{K'\in\mathcal M_k : \operatorname{dist}\left(K,K'\right)\le C\restr{\mathsf h}{K},\; k=0,1,\ldots,N \}\setminus\{K\},
	\end{equation}
	where $C>0$ does not depend on the mesh sizes.
\end{definition}

Next, for each bad cut element $K\in\mathcal M_i^b$, $0\le i < N$,  we want to associate a \emph{good neighbor} $K'$ (a neighbor that is a good element). Note that in principle we allow $K'\in \mathcal M_k$ with $i\neq k$, i.e., a good neighbor can belong to the mesh of another domain. For every $K\in\mathcal M_i^b$, $0\le i < N$, its associated good neighbor $K'$ is chosen according to the procedure in Algorithm~\ref{algorithm1}.

\begin{algorithm}[H]
	\SetAlgoLined
	Given $K\in\mathcal M_i^b$, $0\le i < N$\;
	\For{$k=i,\dots, N$}{
		\If{$\mathcal N(K)\cap\mathcal M_k^g\ne\emptyset$}{
			$K'\gets$ any element of $\mathcal N(K)\cap \mathcal M_k^g$\;  
			\textbf{break} 
		}
	}
	\Return{$K'$}\;
	\caption{Find good neighbor}\label{algorithm1}
\end{algorithm}

\begin{figure}
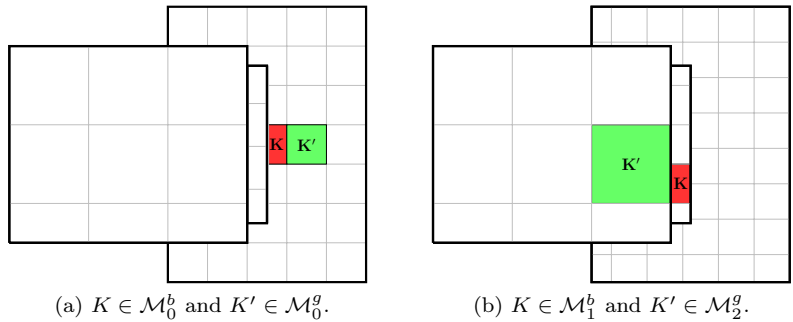

	\centering
	\subfloat[][$K\in\mathcal M_0^b$ and $K'\in\mathcal M_0^g$.]
	{\includestandalone[width=0.35\textwidth]{good_bad_choice_multi1}
	}
	\hspace{0.5cm}
	\subfloat[][$K\in\mathcal M_1^b$ and $K'\in\mathcal M_2^g$.]
	{
		\includestandalone[width=0.35\textwidth]{good_bad_choice_multi}\label{fig:goodnb:bad}
	}
	\caption{The good neighbor $K'$ located in the same domain as $K$ (a), and in a different (top) domain (b).}
	\label{fig:goodnb}
\end{figure}

If Algorithm~\ref{algorithm1} does not produce any output, then it suffices to relax the definition of the good neighbor by taking a larger constant $C$ in Definition~\ref{def_neigh}. Figure~\ref{fig:goodnb} shows two choices of good neighbor. In Figure~\ref{fig:goodnb}\subref{fig:goodnb:bad} there is clearly no neighbor of $K$ in $\mathcal M_1^g$, and the algorithm chooses it in $\mathcal M_2^g$.

We remark that, thanks to Assumption~\ref{assumption_uniformity}, the number of times an active element can be chosen as a good neighbor of any other element stays bounded. For this reason, we will not track the dependence of our constant from this number.

%

The following assumption is not restrictive and is satisfied whenever the meshes of the patches are sufficiently refined, and we take the constant $C$ in~\eqref{constante_scuola} large enough. 	

\begin{assumption}
For every $K\in\mathcal M_i^b$, $0\le i < N$, there exist $i\le k\le N$ and $K'\in\mathcal N(K)\cap\mathcal M_k^g$. From now on we will refer to such $K'$ as a \emph{good neighbor} of $K$.	
\end{assumption}	

For every $K\in\mathcal M_i^b$, $0\le i < N$, its associated good neighbor $K'$ is chosen according to the procedure described in Algorithm~\ref{algorithm1} in Section~\ref{sec:ministab}. As already observed, if Algorithm~\ref{algorithm1} does not produce any output, then it suffices to relax the definition of the good neighbor by taking a larger constant $C$ in~\eqref{constante_scuola}.

We observe that formulation~\eqref{non_stabilized_formulation} is well-posed if there are no cut elements, i.e., {$\mathcal M_{k}^b=\emptyset$}, for every $0\le k\le N$. In the general case {$\mathcal M_{k}^b\ne\emptyset$}, for some $0\le k< N$, the goal of the stabilization is, informally speaking, to extend the stability of the discrete formulation from the internal elements of the visible parts of the patches up to their cut elements.

Let us start by addressing i), i.e., stabilizing the pressure. For $0\le j< i \le N$, $\ell\in\{i,j\}$, let us define $R_{\ell}^p:Q_{h}\to L^2(\Omega_\ell )$ locally. For every $q_h\in Q_h$ and $K\in\mathcal M_\ell$, we distinguish two cases:
\begin{itemize} 
	\item if $K\in\mathcal M_\ell^g$, then
	\begin{equation*}
		\restr{R_\ell^p(q_h)}{K}:=\restr{q_\ell}{K},
	\end{equation*}	
	\item if $K\in\mathcal M_\ell^b$, then
	\begin{equation*}
		\restr{R_\ell^p(q_h)}{K}:=\restr{{\mathcal E}_{K',K}\left( \Pi_{K'}\left(\restr{q_k}{K'} \right)\right)}{K},
	\end{equation*}	
	where {$K'\in\mathcal M_k^g$ ($\ell\leq k\leq N$)}, $\Pi_{K'}: L^2(K')\to \mathbb Q_k (K')$ is the $L^2$-orthogonal projection onto the space of tensor product polynomials defined on the good neighbor $K'$, and $\mathcal E_{K',K}: \mathbb Q_k(K')\to\mathbb Q_k (K\cup K')$ is the canonical polynomial extension. {Note that thanks to the identification $\restr{q_h}{K'}=\restr{q_k}{K'}$. The projection essentially transforms composite spline functions (i.e., splines composite with geometric mapping) to a polynomial representation, which makes the subsequent polynomial extension readily applicable directly in the physical domain. Moreover, functions supported exclusively on badly cut element are removed from the pressure space, and thus projection is necessary.}
\end{itemize}
\begin{lemma}[Stability properties of $R^p_\ell$]\label{lemma:stability_pressure}
	Given $\theta\in (0,1]$, there exist $C_{S,1},C_{S,2}>0$ such that, for every $0\le j<i\le N$, $K\in\mathcal M_\ell$, $\ell\in\{i,j\}$, and $q_h\in Q_h$, we have
	\begin{align*}
		\norm{\mathsf h_{\ell}^{\frac{1}{2}}R^p_{\ell}(q_h)}_{L^2(\Gamma_{ij}\cap\overline K)}\le  C_{S,1} \norm{q_k}_{L^2( K'\cap\Omega)},\quad	\norm{ R^p_{\ell}(q_h)}_{L^2(K\cap\Omega_\ell)}\le  C_{S,2} \norm{q_k}_{L^2( K'\cap\Omega)},
	\end{align*}	
	where $K'\in\mathcal M_k^g$ is a good neighbor of $K$ if $K\in\mathcal M_\ell^b$, or $K'=K$ if $K\in\mathcal M_\ell^g$.
\end{lemma}	
\begin{proof}
	The proof is analogous to Theorem 6.4 in~\cite{ref:puppi20}, so we skip it.
\end{proof}

The following corollary is an easy consequence of Lemma~\ref{lemma:stability_pressure}.
\begin{corollary}\label{prop_press_stability}
	There exists $C>0$, depending on $N_\Gamma$, such that for every $q_h\in Q_h$, we have
	\begin{align*}
		\sum_{i=1}^N \sum_{j=0}^{i-1} \norm{\h^{\frac{1}{2}} [R_{ij}^p(q_h)]}^2_{L^2(\Gamma_{ij})}\le C \sum_{i=0}^N \norm{q_h}^2_{L^2(\Omega_i)},
	\end{align*}	
	where $[R^p_{ij}(q_h)]:= \restr{R^p_i(q_h)}{\Gamma_{ij}}-\restr{R_j^p(q_h)}{\Gamma_{ij}}$ for $0\le j< i\le N$.
\end{corollary}	

Now, we propose a minimal stabilization of the bilinear form $a_h(\cdot,\cdot)$. For $0\le j< i \le N$, $\ell\in\{i,j\}$, let us define $R_{\ell}^v:V_{h}\to \bm L^2(\Omega_\ell )$ locally. For every $\vv_h\in V_h$ and {$K\in\mathcal M_\ell$}, we distinguish two cases:
\begin{itemize} 
	\item if $K\in\mathcal M_\ell^g$, then
	\begin{equation*}
		\restr{R_\ell^v(\vv_h)}{K}:=\restr{\vv_\ell}{K},
	\end{equation*}	
	\item if $K\in\mathcal M_\ell^b$, then
	\begin{equation*}
		\restr{R_\ell^v(\vv_h)}{K}:=\restr{\bm{\mathcal E}_{K',K}\left(\bm \Pi_{K'}\left(\restr{\vv_k}{K'} \right)\right)}{K},
	\end{equation*}	
	where $\bm\Pi_{K'}:\bm L^2(K')\to \mathbb V_k (K')$ is the $L^2$-orthogonal projection,
	\begin{align*}
		\mathbb V_k (K):=	&
		\begin{cases}	
			\mathbb S_k(K)\qquad\qquad\qquad\qquad\qquad\qquad\qquad\qquad\;\;\;\,\;&\text{if}\ \square= \mathrm{RT},\\
			\left(\mathbb Q_{k+1}(K)\right)^d\qquad&\text{if}\ \square\in \{\mathrm{N},\mathrm{TH}\},\\
		\end{cases}\\
		\mathbb S_k(K) :=&
		\begin{cases}
			\mathbb Q_{k+1,k}(K)\times \mathbb Q_{k,k+1}(K)\qquad\qquad\qquad\qquad\qquad\, &\text{if}\ d=2, \\
			\mathbb Q_{k+1,k,k}(K)\times \mathbb Q_{k,k+1,k}(K) \times \mathbb Q_{k,k,k+1} (K)\qquad &\text{if}\ d=3, \\	
		\end{cases}	
	\end{align*}
	and $\bm{\mathcal E}_{K',K}: \mathbb V_k(K')\to \mathbb V_k(K\cup K')$ is the canonical polynomial extension. 
\end{itemize}
We denote, for $0\le j< i\le N$ such that $\Gamma_{ij}\ne\emptyset$ and $t\in\{\frac{1}{2},1\}$,
\begin{equation}\label{eq:stabilized_vel}
	\langle D R^v_{ij} (\vv_h)\n  \rangle_t : = t \restr{D R^v_i\left(\vv_h\right)\n_i}{\Gamma_{ij}}+\restr{\left(1-t\right) D R^v_j\left(\vv_h\right) \n_j}{\Gamma_{ij}}.
\end{equation}
In the previous definition, we used the spaces $V_h$ and $V_{h,\ell}$. Let us point out that the stabilization operators $R^v_\ell$ are equally defined for the spaces with zero boundary conditions $V_h^{\bm 0}$ and $ V^{\bm 0}_{h,\ell}$, $0\le \ell\le N$, as well. 

{
\begin{remark}
Stabilization of the flux term $\langle D R^v_{ij} (\vv_h)\n  \rangle_t$ is generally needed for both the symmetric case ($t=\frac{1}{2}$) and the one-sided case ($t=1$). The symmetric case is straightforward as it involves the flux from a (trimmed) low-level patch. In the one-sided case, however, whether stabilization is needed depends on how many patches are overlapped. If only two patches are overlapped, the flux comes from the (non-trimmed) high-level patch, so there is no need to stabilize. However, if multiple overlapping patches are involved, a certain middle-level patch serves as a high-level patch to the patches below it, so it provides the flux. But at the same time it may be trimmed by the patches on top of it. As a result, the one-sided flux may come from a trimmed patch, so stabilization is needed.
\end{remark}
}

\begin{lemma}[Stability property of $R^v_\ell$]\label{lemma:stability_velocity}
	Given $\theta\in (0,1]$, there exists $C>0$ such that for every $0\le j<i\le N$, $K\in\mathcal M_\ell$, $\ell\in\{i,j\}$, and $\vv_h\in V_h$, we have
	\begin{equation*}
		\norm{\mathsf h_{\ell}^{\frac{1}{2}}D R^v_{\ell}(\vv_h)\n_\ell}_{L^2(\Gamma_{ij}\cap\overline K)}\le C \norm{D \vv_k}_{L^2( K'\cap\Omega)},
	\end{equation*}	
	where $K'\in\mathcal M_k^g$ is a good neighbor of $K$ if $K\in\mathcal M_\ell^b$, or $K'=K$ if $K\in\mathcal M_\ell^g$.
\end{lemma}	

\begin{proof}
	We refer the interested reader to its scalar counterpart, i.e., the proof of Theorem 6.4 in~\cite{ref:puppi20}, or Lemma 3.6 of~\cite{ref:wei21a}. 
\end{proof}

\begin{corollary}\label{proposition_stability_velocity}
	There exists $C>0$, depending on {$N_\Gamma$} such that for every $\vv_h\in V_h$,
	\begin{equation*}
		\sum_{i=1}^N \sum_{j=0}^{i-1}\norm{\h^{\frac{1}{2}}\langle D {R^v_{ij}}(\vv_h)\n\rangle_t}_{L^2(\Gamma_{ij})}\le C \sum_{i=0}^N \norm{D \vv_i}_{L^2(\Omega_i)}
	\end{equation*}	
\end{corollary}	
\begin{proof}
	This result is a direct consequence of Lemma~\ref{lemma:stability_velocity}.
\end{proof}

We introduce the following stabilized space for the pressures.
\begin{align*}
	\overline Q_{h}=\{ \varphi_h\in L^2(\Omega):&\ \exists\ q_h\in Q_{h}\ \text{such that}\ \restr{\varphi_h}{K}=\restr{q_h}{K}\ \forall\ K\in\mathcal M_i^g\\
	 &\text{and}\ \restr{\varphi_h}{K} =\restr{R^p_h(q_h)} {K}\ \forall\ K\in\mathcal M_i^b,
	 \forall\ 0\le i \le N  \}.	
\end{align*}	
We can finally propose our stabilized weak formulation.

Find $\left(\u_h,p_h\right)\in V_h^{\u_D}\times \overline Q_h$ such that
\begin{equation}
	\begin{aligned}\label{stabilized_formulation}
		\overline	a_h(\u_h,\vv_h) + b(\vv_h,p_h) = F(\vv_h),\qquad &\forall\ \vv_h\in V_h^{\bm{0}},\\
		b(\u_h,q_h) = 0,\qquad &\forall\ q_h\in \overline Q_h,
	\end{aligned}
\end{equation}
where
\begin{equation}
	\begin{aligned}
		\overline a_h(\w_h,\vv_h):=&\sum_{i=0}^N\int_{\Omega_i}D \w_i: D \vv_i
		-\sum_{i=1}^N\sum_{j=0}^{i-1}\int_{\Gamma_{ij}}\left(\langle DR^v_{ij}(\w_h)\n \rangle_t\left[\vv_h\right]+\left[\w_h\right]\langle DR^v_{ij}(\vv_h) \n \rangle_t\right)\\
		&+\gamma\sum_{i=1}^N\sum_{j=0}^{i-1}\int_{\Gamma_{ij}}\h^{-1}\left[\w_h\right]\left[\vv_h\right], \qquad\ \w_h,\vv_h\in V_h,  \\
		b(\vv_h,q_h) : = & - \sum_{i=0}^N \int_{\Omega_i} q_i \dive \vv_i + \sum_{i=1}^N \sum_{j=0}^{i-1} \int_{\Gamma_{ij}} \langle q_h \rangle_t [\vv_h\cdot\n ],\qquad\ \vv_h\in V_h, q_h\in \overline Q_h,
	\end{aligned}
\end{equation}
with $t\in\{\frac{1}{2},1\}$, and
\begin{equation*}
	\begin{aligned}
		F(\vv_h):=&\sum_{i=0}^N\int_{\Omega_i}\f\cdot\vv_i + \int_{\Gamma_N} \u_N\cdot \vv_h,\qquad\ \vv_h\in V_h.
	\end{aligned}
\end{equation*}
As before $\gamma>0$ is a penalty parameter scaling as the spline degree of the velocity.

{
\begin{remark}
We note that after the discrete pressure space $Q_h$ is altered to the stabilized counterpart $\overline{Q}_h$, the structure-preserving property of isogeometric Raviart-Thomas elements~\cite{MR3048532} no longer holds. On the other hand, the difference between $Q_h$  and  $\overline{Q}_h$ is localised at cut elements, and so will be the effect of lack of structure preserving properties. 
\end{remark}
}

\subsection{Interpolation and approximation properties of the discrete spaces}
Before analyzing problem~\eqref{stabilized_formulation}, we need some technical results. {Given a Sobolev} function living in the whole physical domain $\Omega$, we consider its restrictions to the predomains $\Omega_i^\ast$ in order to be able to interpolate on each premesh $\mathcal M^\ast_i$, restrict them in their turn to the visible parts $\Omega_i$, and finally glue together the interpolated functions. 

Let us start with the velocity. We construct a quasi-interpolant operator for each local space $V_{h,i}$. Given $\square\in\{\mathrm{RT},\mathrm{N},\mathrm{TH}\}$ and $m\ge 1$, for every $i\in\{0,\dots,N\}$, we denote
\begin{equation*}
	\Pi^{i,\square}_{V_h}:\bm H^{m}(\Omega_i)\to V^\square_{h,i},\qquad
	\vv\mapsto \restr{\Pi_{V_h}^{i,\ast,\square}\left(\restr{\vv}{\Omega_i^\ast}\right)}{\Omega_i} ,
\end{equation*}
where $\Pi_{V_h}^{i,\ast,\square}: \bm H^{m}(\Omega_i^{\ast})\to V^{\ast,\square}_{h,i}$ is a standard quasi-interpolant operator~\cite{MR3202239,MR3585787}. Then, we glue together the local operators as
\begin{equation*}
	\Pi^\square_{V_h}:\bm H^{m}(\Omega)\to V_h^\square,\qquad \vv\mapsto \bigoplus_{i=0}^N  \Pi_{V_h}^{i,\square}\left(\vv_i\right),
\end{equation*}
where $\vv_i(x):=\vv(x)$ for every $x\in\Omega_i$, $i=0,\dots,N$. For the pressures, given $r\ge 1$, we introduce the local quasi-interpolants
\begin{align*}
	\Pi^{i,\square}_{Q_h}: H^{r}(\Omega_i)\to Q_{h,i}^\square,\qquad
	q \mapsto  \restr{\Pi_{Q_h}^{i,\ast,\square}\left(\restr{q}{\Omega_i^\ast}\right)}{\Omega_i} ,
\end{align*}
where $\Pi_{Q_h}^{i,\ast,\square}:  H^{r}(\Omega_i^{\ast})\to Q^{\ast,\square}_{h,i}$ is a standard quasi-interpolant operator. We glue together the local operators for the pressures as
\begin{equation*}
	\Pi^\square_{Q_h}: H^{r}(\Omega)\to \overline Q_h^\square,\qquad q\mapsto R_h^p \left( \bigoplus_{i=0}^N  \Pi_{Q_h}^{i,\square}\left(q_i\right)\right),
\end{equation*}
where $q_i(x):=q(x)$ for every $x\in\Omega_i$, $i=0,\dots,N$.

\begin{proposition}[Interpolation error estimate]\label{inter_err_est}
	There exist $C_v,C_p>0$ such that, for every $(\vv,q)\in \bm H^{m}(\Omega)\times H^{r}(\Omega)$, $m\ge 1$ and $r\ge 1$, it holds
	\begin{equation*}
		\norm{\vv-\Pi_{V_h}^\square \vv}_{1,h}\le C_v h^s\norm{\vv}_{H^m(\Omega)},\qquad \norm{q-\Pi_{Q_h} q}_{0,h} \le C_p{h}^\ell \norm{q}_{H^{r}(\Omega)},
	\end{equation*}	
	where $s:=\min \{k,m-1\}$ if $\square=\mathrm{RT}$, $s:=\min\{k+1,m-1\}$ if $\square \in\{\mathrm{N},\mathrm{TH} \}$, and $\ell:=\min\{k+1,r\}$.
\end{proposition}	
\begin{proof}
	The proof follows from the best approximation properties of the local quasi-interpolant operators $\Pi_{V_h}$; see e.g.,~\cite{MR3202239}. On the other hand, for $\Pi_{Q_h}$, the statement is a direct consequence of Lemma~\ref{lemma:stability_velocity}, and the fact that, by construction, $R_h^p(\varphi)=\varphi$ for all $\varphi\in \mathbb{Q}_k$.

\end{proof}

\begin{lemma}[Approximation property of $R^v_\ell$]\label{lemma_approx}
	There exists $C>0$, depending on $N_\Gamma$ and $N_{\mathcal O}$, such that, for every $\vv\in \bm H^m(\Omega)$, $m\ge 2$, it holds
	\begin{equation*}
		\sum_{i=1}^N \sum_{j=0}^{i-1}	\norm{\h^{\frac{1}{2}}\langle D R^v_{ij} \left( \Pi_{V_h}^{\square}\left(\vv \right)\right)\n -D\vv\n \rangle_t }^2_{L^2(\Gamma_{ij})}\le C h^{2s}\norm{ \vv}^2_{H^{m}\left(\Omega\right)},
	\end{equation*}
	where $s:=\min\{k,m-1\}$ if $\square=\mathrm{RT}$ and $s:=\min\{k+1,m-1\}$ if $\square \in\{\mathrm{N},\mathrm{TH} \}$.
\end{lemma}
\begin{proof}
	The proof of such a statement in the scalar case can be found in Proposition 6.9 in~\cite{ref:puppi20} and Proposition 4.11 in~\cite{puppi_stokes}.
\end{proof}

\subsection{Well-posedness of the stabilized formulation}
{In order to avoid a direct proof of the inf-sup condition for  $b(\cdot,\cdot)$, following ~\cite{ref:stenberg95} and also \cite{puppi_stokes}, we 
  reframe the Nitsche formulation as a stabilized Lagrange multiplier method}. In doing so, we find a formulation equivalent to~\eqref{stabilized_formulation}, but whose well-posedness is easier to prove. 
First of all, we observe that the transmission conditions at the interfaces~\eqref{problem2.2} and~\eqref{problem2.3} can be enforced by introducing a Lagrange multiplier living in $\Lambda:=\bigoplus_{0\le j<i\le N} H_{00}^{-\frac{1}{2}}(\Gamma_{ij})$, where, for $0\le j<i\le N$, $H_{00}^{-\frac{1}{2}}(\Gamma_{ij}):= \left( H_{00}^{\frac{1}{2}} (\Gamma_{ij}) \right)'$ and $H^{\frac{1}{2}}_{00}(\Gamma_{ij}):=\{\varphi\in L^2(\Gamma_{ij}): \tilde \varphi \in H^{\frac{1}{2}}(\Gamma_j) \}$, $\tilde \varphi$ denoting the extension by zero of $\varphi\in L^2(\Gamma_{ij})$ on $\Gamma_j\setminus\Gamma_{ij}$ (see Chapter~11 of~\cite{MR3930592}). Let $\Lambda_h$ be a discrete subspace of $\Lambda$ large enough such that, for every $\left(\vv_h,q_h \right) \in V_h \times \overline Q_h$, it holds
\begin{align}\label{condition:multiplier_space}
	\h^{-1} \restr{\left[\vv_h\right]}{\Gamma_{ij}}, \restr{\langle D R^v_{ij}(\vv_h)\n \rangle_t }{\Gamma_{ij}}, \restr{\langle q_h\n \rangle_t }{\Gamma_{ij}}   \in\Lambda_h,\qquad \forall\ 0\le j<i\le N, t\in\{\frac{1}{2},1 \}.
\end{align}
Our choice is:
\begin{align}\label{choice_multiplier}
	\Lambda_h: = \bigoplus_{0\le j<i\le N} \Lambda_{ij} , \qquad \Lambda_{ij}:= W_h (\Gamma_{ij}) + L_h (\Gamma_{ij}) \n + N_h (\Gamma_{ij}) \n,
\end{align}

\begin{align*}
	W_h (\Gamma_{ij}):=&\{\h^{-1}\restr{\w_h}{\Gamma_{ij}}: \w_h \in V_h \},\\	
	L_h (\Gamma_{ij}):=& \{\restr{D\w_h}{\Gamma_{ij}}: \w_h \in V_h \},\\	
	N_h(\Gamma_{ij}):=& \{\restr{\varphi_h}{\Gamma_{ij}}: \varphi_h \in \overline{Q}_h \}.
\end{align*}
It is easy to check that~\eqref{choice_multiplier} satisfies conditions~\eqref{condition:multiplier_space}. Moreover, $\Lambda_h\neq \emptyset$. We endow $\Lambda_h$ with the mesh-dependent norm
\begin{align*}
	\norm{\bmu_h}^2_{-\frac{1}{2},h}:= \sum_{i=1}^N\sum_{j=0}^{i-1} \norm{\h^{\frac{1}{2}} \bmu_h}^2_{L^2(\Gamma_{ij})}, \qquad \ \bmu_h\in \Lambda_h.
\end{align*}
Let us introduce the following stabilized augmented Lagrangian formulation.

Find $\left( \u_h,p_h,\blambda_h\right) \in V_h^{\u_D}\times \overline Q_h\times \Lambda_h$ such that
\begin{align}\label{stabilized_augmented_formulation}
	\overline{\mathcal A}_h \left( \left( \u_h,p_h,\blambda_h \right) ; \left(\vv_h,q_h,\bmu_h \right) \right) = \mathcal F \left( \vv_h,q_h,\bmu_h\right), \quad \forall\ (\vv_h,q_h,\bmu_h)\in V_h^{\bm 0}\times \overline Q_h\times \Lambda_h,
\end{align}	
where, for $\left( \w_h,r_h,\bm\eta_h \right),\left(\vv_h,q_h,\bmu_h \right)\in V_h\times \overline Q_h\times \Lambda_h$,
\begin{align*}
	\overline {\mathcal A}_h ( \left( \w_h,r_h,\bm\eta_h \right) ;& \left(\vv_h,q_h,\bmu_h \right) ) := 	a(\w_h,\vv_h)+ b_0(\vv_h,r_h) + b_\Gamma (\vv_h,\bm\eta_h)\\
	&+ b_0(\w_h,q_h) + b_\Gamma (\w_h,\bmu_h)\\
	&-\sum_{i=1}^{N} \sum_{j=0}^{i-1}\gamma^{-1}\int_{\Gamma_{ij}}\h  \left(\bm\eta_h +\langle D R_{ij}^v(\w_h)\n \rangle_t - \langle r_h\n\rangle_t\right) \left(\bm\bmu_h +\langle D R_{ij}^v(\vv_h)\n \rangle_t - \langle q_h\n\rangle_t\right),\\
	\mathcal F \left( \vv_h,q_h,\bmu_h\right) :=& F(\vv_h),\qquad 	a(\w_h,\vv_h) := \sum_{i=0}^N \int_{\Omega_i} D\w_i: D\vv_i,\\
	b_0(\vv_h,q_h):= &-\sum_{i=0}^N \int_{\Omega_i} q_i \dive \vv_i,\qquad 	b_\Gamma (\vv_h,\bmu_h):= \sum_{i=1}^{N} \sum_{j=0}^{i-1} \int_{\Gamma_{ij}} \bmu_h [\vv_h].
\end{align*}
\begin{proposition}\label{stability:prop1_stab}
	Let $\Lambda_h$ be  defined as in~\eqref{choice_multiplier}. Then, problem~\eqref{stabilized_augmented_formulation} is equivalent to the stabilized Nitsche's formulation~\eqref{stabilized_formulation}.
\end{proposition}
\begin{proof}
	We can reformulate problem~\eqref{stabilized_augmented_formulation} as follows.
	
	Find $(\u_h,p_h,\blambda_h)\in V_h^{\u_D} \times \overline Q_h \times \Lambda_h$ such that
	\begin{alignat*}{3}
		a(\u_h,\vv_h)+ b_0(\vv_h,p_h) + b_\Gamma (\vv_h,\blambda_h)\\
		-\sum_{i=1}^{N} \sum_{j=0}^{i-1}\gamma^{-1}\int_{\Gamma_{ij}}\h \left(\blambda_h +\langle D R_{ij}^v(\u_h)\n \rangle_t - \langle p_h\n\rangle_t\right) \langle D R_{ij}^v(\vv_h)\n \rangle_t &= F(\vv_h), \qquad &&\forall\ \vv_h\in V_h^{\bm 0},\\
		b_0(\u_h,q_h) + \sum_{i=1}^{N} \sum_{j=0}^{i-1} \gamma^{-1} \int_{\Gamma_{ij}}\h \left(\blambda_h + \langle D R_{ij}^v(\u_h)\n \rangle_t - \langle p_h\n\rangle_t\right) \langle q_h\rangle_t &=0,\qquad &&\forall\ q_h\in \overline Q_h,\\
		b_\Gamma (\u_h,\bmu_h) -\sum_{i=1}^{N} \sum_{j=0}^{i-1} \gamma^{-1} \int_{\Gamma_{ij}}\h \left(\blambda_h + \langle D R_{ij}^v(\u_h)\n \rangle_t - \langle p_h\n\rangle_t\right) \bmu_h&=0, \qquad &&\forall\ \mu_h\in \Lambda_h.
	\end{alignat*}
	From the third equation, we can \emph{statically condensate} the Lagrange multiplier as $\displaystyle\bm\lambda_h=\sum_{i=1}^N \sum_{j=0}^{i-1}\bm\lambda_{ij}$, where
	\begin{align*}
		\blambda_{ij} = \sum_{i=1}^N \sum_{j=0}^{i-1} \gamma  P_{ij}\left(\h^{-1}[\u_h]\right)-  P_{ij}\left(\langle DR_{ij}^v\left(\u_h\right)\n\rangle_t \right) +P_{ij}\left(\langle p_h\n\rangle_t\right),
	\end{align*}
	where $P_{ij}:\bm L^2(\Gamma_{ij})\to \Lambda_{ij}$ is the $L^2$-orthogonal projection. By substituting it back to the first and second equations and using~\eqref{choice_multiplier}, we recover the desired formulations. 
\end{proof}

\begin{remark}
	{In what follows, we prove the wellposedness of ~\eqref{stabilized_augmented_formulation} by showing the invertibility of the arising global systems through the so-called \emph{Banach-Neças-Babuška Theorem}~\cite{MR2050138}. On the other hand, we are no longer bound to satisfy an inf-sup condition on the bilinear form $b(\cdot,\cdot)$, and for this reason, we are free to choose the Lagrange multiplier space large enough so that the hypotheses of Proposition~\ref{stability:prop1_stab} are satisfied. }
\end{remark}

A fundamental ingredient for our numerical analysis is that the following \emph{local inf-sup conditions} are satisfied on each visible part of the domain. 
\begin{assumption}\label{assumption:local_is}
	For every $0\le i\le N$, let us equip $V_{h,i}$ with 
	\begin{align}\label{caprone}
		\norm{\vv_i}^2_{1,h,i}:=\norm{D\vv_i}^2_{L^2(\Omega_i)}+\sum_{k=i+1}^N \norm{\h^{-\frac{1}{2}} \vv_i}^2_{L^2(\Gamma_{ki})} +\sum_{j=0}^{i-1}\norm{\h^{-\frac{1}{2}} \vv_i}^2_{L^2(\Gamma_{ij})},\vv_i\in V_{h,i},
	\end{align}	
	and $\restr{\overline Q_{h}}{\Omega_i}$ with $\norm{\cdot}_{L^2(\Omega_i)}$. We assume that, for every $0\le i\le N$, there exist $\beta_i$ such that $\forall q_i \in \restr{\overline Q_{h}}{\Omega_i}$,
\begin{equation}
\underset{\vv_i \in V_{h,i}}{\sup} \frac{\int_{\Omega_i} q_i \dive \vv_i}{\norm{\vv_i}_{1,h,i}} \ge {\beta_i} \norm{q_i}_{L^2(\Omega_i)}.
\end{equation}
\end{assumption}	

Although we are not able to show a proof of Assumption~\ref{assumption:local_is}, it is a local result that can be tested. Extensive inf-sup tests are provided in~\cite{puppi_stokes}, which provide a solid numerical basis for this paper.

The following Lemma provides a norm equivalence result for the space of stabilized pressures, useful for the proof of the inf-sup stability of $b_0(\cdot,\cdot)$. 
\begin{lemma}\label{lemma:norm_equivalence}
	There exists $C>0$ such that 
	\begin{align}\label{eq:norm_equivalence}
		\sum_{i=0}^N \norm{q_i}^2_{L^2(\Omega_i)}\le \norm{q_h}^2_{0,h}\le C\sum_{i=0}^N \norm{q_i}^2_{L^2(\Omega_i)}, \qquad\forall\ q_h\in \overline Q_h.
	\end{align}		
\end{lemma}	
\begin{proof}
	Let us take $q_h\in \overline Q_h$. By definition,
	\begin{align*}
		\norm{q_h}^2_{0,h}= &\sum_{i=0}^N \norm{q_i}^2_{L^2(\Omega_i)} + \sum_{i=1}^N \sum_{j=0}^{i-1} \norm{\h^{\frac{1}{2}} [q_h]}^2_{L^2(\Gamma_{ij})}.
	\end{align*}	
	Hence, we use Proposition~\ref{prop_press_stability} so that
	\begin{align*}
		\norm{q_h}^2_{0,h}\le & C \sum_{i=0}^N \norm{q_i}^2_{L^2(\Omega_i)},
	\end{align*}
	where, in particular, $C$ depends on $N_\Gamma$ and $\theta$. The other inequality trivially holds.
\end{proof}

The following inf-sup condition for $b_0(\cdot,\cdot)$ is a key ingredient for the proof of the well-posedness of formulation~\eqref{stabilized_augmented_formulation}.
\begin{lemma}\label{lemma:is}
	Under Assumption~\ref{assumption:local_is}, given $\theta\in(0,1]$, there exists $\beta_0>0$ such that
	\begin{equation}\label{eq:is}
		\begin{aligned}
			\inf_{q_h\in \overline{Q}_h} \sup_{\vv_h \in V_h} \frac{b_0(\vv_h,q_h)}{\norm{\vv_h}_{1,h} \norm{q_h}_{0,h}} \ge \beta_0.	
		\end{aligned}	
	\end{equation}
\end{lemma}	
\begin{proof}
	We prove that there exist $C_1,C_2>0$ such that, for every $q_h\in \overline Q_h$, there exists $\vv_h\in V_h$ such that
	\begin{equation*}
		\begin{aligned}
			b_0(\vv_h,q_h)\ge C_{1} \norm{q_h}^2_{0,h},\qquad
			\norm{\vv_h}_{1,h} \le C_{2} \norm{q_h}_{0,h}.
		\end{aligned}
	\end{equation*}
	Let us fix $q_h:=\left(q_0,\dots,q_N\right)\in \overline Q_h$. From~Assumption~\ref{assumption:local_is}, there exist $C_{1,i},C_{2,i}>0$ and $\vv_i\in V_{h,i}$, $0\le i\le N$, such that
\begin{equation*}
\int_{\Omega_i} q_i \dive \vv_i \ge C_{1,i} \norm{q_i}_{L^2(\Omega_i)}, \qquad
\norm{\vv_i}_{1,h,i} \le C_{2,i} \norm{q_i}_{L^2(\Omega_i)}.
\end{equation*}
Let $\vv_h:=\left(\vv_0,\dots,\vv_N\right)$. By letting $\overline C_1:=\min_{0\le i\le N} C_{1,i}$, we have
	\begin{align*}
		b_0(\vv_h,q_h)\ge \overline C_1\sum_{i=0}^N\norm{q_i}^2_{L^2(\Omega_i)}.
	\end{align*}
	On the other hand, from the Young inequality,
	\begin{align*}
		\norm{\vv_h}^2_{1,h} \le \sum_{i=0}^N \norm{D \vv_i}^2_{L^2(\Omega_i)} + 2\left( \sum_{i=1}^N \sum_{j=0}^{i-1} \norm{\h^{-\frac{1}{2}}\vv_i}^2_{L^2(\Gamma_{ij})} +\sum_{i=1}^N \sum_{j=0}^{i-1} \norm{\h^{-\frac{1}{2}}\vv_j}^2_{L^2(\Gamma_{ij})} \right).
	\end{align*}
	Notice that $\sum_{i=1}^N \sum_{j=0}^{i-1}\norm{\h^{-\frac{1}{2}} \vv_j}^2_{L^2(\Gamma_{ij})} =\sum_{i=0}^{N-1}\sum_{k=i+1}^N \norm{\h^{-\frac{1}{2}} \vv_i}^2_{L^2(\Gamma_{ki})}$. Hence, from the definition of $\norm{\cdot}_{1,h,i}$~\eqref{caprone} and Assumption~\ref{assumption:local_is}, it holds
	\begin{align}\label{eq:proof_is1}
		\norm{\vv_h}^2_{1,h} \le 	C	\sum_{i=0}^N\norm{\vv_i}^2_{1,h,i}\le \overline C_2^2\sum_{i=0}^N
		\norm{q_i}^2_{L^2(\Omega_i)},
	\end{align}	
	where $\overline C_2:=(C\tilde C)^{\frac{1}{2}}$ and $\tilde C:=\max_{0\le i\le N} C_{2,i}$. We conclude by taking $C_1:=\overline C_1 C^{-1}$ and $C_2:= \overline C_2$, where $C>0$ come from Lemma~\ref{lemma:norm_equivalence}. In particular, note that $C$ depends on $\theta$.
\end{proof}	

Let us indirectly study the well-posedness of the problem~\eqref{stabilized_formulation} by showing that~\eqref{stabilized_augmented_formulation} verifies the hypotheses of the so-called Banach-Neças-Babuška Theorem~\cite{MR2050138}. The proof of the next result is given in~\ref{appendix:A1}.

\begin{theorem}\label{thm:wellposedness}
	Under Assumption~\ref{assumption:local_is}, there exists $\overline \gamma>0$ and $C>0$ such that, for every $\gamma\ge\overline\gamma$,
	\begin{align}\label{wellposed:eq1}
		\inf_{\left( \vv_h,q_h,\bmu_h\right) \in V_h\times \overline Q_h\times \Lambda_h} \sup_{\left( \w_h, r_h,\bm\eta_h \right) \in V_h\times \overline Q_h\times \Lambda_h} \frac{\overline {\mathcal A}_h \left( \left( \w_h,r_h,\bm\eta_h \right) ; \left(\vv_h,q_h,\bmu_h \right) \right) }{ \vertiii{ (\w_h,r_h,\bm\eta_h)}  \vertiii{ (\vv_h,q_h,\bmu_h)}  }\ge C,
	\end{align}
	where $\vertiii{\cdot}$ is the Euclidean product norm on $V_h\times \overline Q_h\times \Lambda_h$. 
\end{theorem}

As it is well known (see, e.g., Theorem 3.4.1 in~\cite{MR3097958}), Theorem~\ref{thm:wellposedness} implies the well-posedness of problem~\eqref{stabilized_augmented_formulation}.

{
We proceed to characterize the existence, uniqueness, and stability for the solution of problem~\eqref{stabilized_formulation} as follows.
\begin{proposition}\label{prop:characterization_stability}
	\begin{enumerate}[(i)]
		\item\label{1} There exists $\overline\gamma>0$ such that, for every fixed $\gamma\ge\overline\gamma$, there exists $M_a>0$ such that
		\begin{align*}
			\abs{\overline a_h(\w_h,\vv_h)}\le M_a \norm{\w_h}_{1,h}\norm{\vv_h}_{1,h},\qquad\forall\ \w_h,\vv_h\in V_h.
		\end{align*}	
		\item\label{2} There exists $\alpha>0$ such that, for every $\gamma\ge\overline\gamma$, 
		\begin{align*}
			\alpha \norm{\vv_h}^2_{1,h}\le \overline a_h(\vv_h,\vv_h),\qquad\forall\ \vv_h\in \ker B,
		\end{align*}
		where $\ker B:=\{\vv_h\in V_h: b(\vv_h,q_h)=0,\quad\forall\ q_h\in \overline Q_h \}$.
		\item\label{3} There exists $M_b>0$ such that
		\begin{align*}
			\abs{b(\vv_h,q_h)} \le M_b \norm{\vv_h}_{1,h}	\norm{q_h}_{0,h},\qquad\forall\ \vv_h\in V_h, \forall \ q_h\in \overline Q_h.
		\end{align*}	
		\item\label{4} There exists $\beta>0$ such that
		\begin{align*}
			\inf_{q_h\in \overline Q_h} \sup_{\vv_h\in V_h} \frac{b(\vv_h,q_h)}{\norm{\vv_h}_{1,h} \norm{q_h}_{0,h} }\ge \beta.
		\end{align*}
		\item\label{5} There exists $M_F>0$ such that
		\begin{align*}
			\abs{F(\vv_h)} \le M_F \norm{\vv_h}_{1,h},\qquad\forall\ \vv_h\in V_h.	
		\end{align*}	 
	\end{enumerate}
	Conditions~\eqref{1}-\eqref{5} hold if and only there exists a unique solution $\left(\u_h,p_h\right)$ of~\eqref{stabilized_formulation}. Moreover,
	\begin{align}\label{eq:stability_estimates}
		\norm{\u_h}_{1,h}\le  \frac{1}{\alpha}\norm{F}_{-1,h},\qquad \norm{p_h}_{0,h}\le  \frac{1}{\beta}\left( 1+\frac{M_a}{\alpha}\right) \norm{F}_{-1,h},
	\end{align}
	where $\norm{\cdot}_{-1,h}$ denotes the dual norm with respect to $\norm{\cdot}_{1,h}$.
\end{proposition}
\begin{proof}
	This is a standard result; see, for instance, Theorem~3.4.1 in~\cite{MR3097958}.
\end{proof}
}

\begin{theorem}
	Under Assumption~\ref{assumption:local_is}, given $\theta\in (0,1]$, there exists a unique solution $\left(\u_h,p_h\right)\in V_h^{\u_D}\times\overline Q_h$ of the stabilized problem~\eqref{stabilized_formulation} satisfying {\eqref{eq:stability_estimates}}.
\end{theorem}
	
\begin{proof}
	From Theorem~\ref{thm:wellposedness} there exists a unique solution $\left(\u_h,p_h,\blambda_h\right)\in V_h^{\u_D}\times \overline Q_h\times \Lambda_h$ of~\eqref{stabilized_augmented_formulation}. Thanks to Proposition~\ref{stability:prop1_stab}, $\left(\u_h,p_h\right)$ is the unique solution of problem~\eqref{stabilized_formulation}. 
\end{proof}

\subsection{A priori error estimates}
The goal of this section is to demonstrate that the errors, for both the velocity and pressure fields, achieve optimal \emph{a priori} convergence rates in the topologies induced by the norms $\norm{\cdot}_{1,h}$ and $\norm{\cdot}_{0,h}$, respectively.
\begin{lemma}\label{lemma:apriori}
	Given $\theta\in (0,1]$, let $\left(\u,p\right)\in \bm H^{\frac{3}{2}+\eps}(\Omega)\times L^2(\Omega)$, $\eps> 0$, and $\left(\u_h,p_h\right)\in  V_h \times \overline Q_h$  be the solutions to~\eqref{stokes} and~\eqref{stabilized_formulation}, respectively. Then, for every $\left(\u^I,p^I \right)\in  V_h\times  \overline Q_h$, the following estimates hold.
	\begin{align*}
		\norm{\u_h - \u^I}_{1,h} \le& \frac{1}{\alpha}\bigg(M_a\norm{\u-\u^I}_{1,h} +\sum_{i=1}^N\sum_{j=0}^{i-1} \norm{\h^{\frac{1}{2}}\langle D\left(\u-R_{ij}^v(\u^I)\right)\n\rangle_t}_{L^2(\Gamma_{ij})}\\
		&+M_{b}\norm{p-p^I}_{0,h}  \bigg)+ \frac{1}{\beta} \left(1+\frac{M_a}{\alpha}\right) M_{b}\norm{\u-\u^I}_{1,h}, \\
		\norm{p_h-p^I}_{0,h} \le & \frac{1}{\beta}\left(1+\frac{M_a}{\alpha} \right) \bigg(M_a\norm{\u-\u^I}_{1,h} +\sum_{i=1}^N\sum_{j=0}^{i-1} \norm{\h^{\frac{1}{2}}\langle D\left(\u-R_{ij}^v(\u^I)\right)\n\rangle_t}_{L^2(\Gamma_{ij})}\\
		&+M_{b}\norm{p-p^I}_{0,h} \bigg) + \frac{M_a}{\beta^2}\left(1+\frac{M_a}{\alpha}\right) M_{b}\norm{\u-\u^I}_{1,h}.
	\end{align*}
\end{lemma}
\begin{proof}
	Let us fix $\left(\u^I,p^I\right)\in V_h\times \overline Q_h$. By linearity $\left(\u_h-\u^I,p_h-p^I \right)\in  V_h\times \overline  Q_h$ satisfies the saddle point problem
	\begin{equation}
		\begin{aligned}\label{eq:weak_cons}
			\overline a_h(\u_h - \u^I,\vv_h) + b(\vv_h,p_h-p^I) =  F^I(\vv_h), \qquad&\forall\ \vv_h\in V_h^{\bm 0},\\
			b(\u_h-\u^I,q_h) = G^I(q_h),\qquad&\forall\ q_h\in  \overline Q_h,
		\end{aligned}
	\end{equation}
	where, for every $\left(\vv_h,q_h\right)\in V_h^{\bm{0}}\times\overline Q_h$,
	\begin{align*}
		F^I(\vv_h):= & \sum_{i=0}^N \int_{\Omega_i} D(\u-\u^I):D\vv_h -\sum_{i=1}^N\sum_{j=0}^{i-1}\int_{\Gamma_{ij}}\big( \langle D \left( \u- R^v_{ij}(\u^I)\right)\n \rangle_t [\vv_h]\\
		& + \langle DR_{ij}^v(\vv_h)\n\rangle_t [\u-\u^I]  \big)+ \gamma\sum_{i=1}^N\sum_{j=0}^{i-1} \int_{\Gamma_{ij}} \h^{-1} [\u-\u^I] [\vv_h] + b(\vv_h,p-p^I),\\
		G^I(q_h):=& b(\u-\u^I,q_h).
	\end{align*}
It is sufficient to apply the stability estimates~\eqref{eq:stability_estimates} to the solution $(\u_h-\u^I,p_h-p_I)$ of~\eqref{eq:weak_cons}.
\end{proof}
\begin{theorem}\label{thm:apriori}
	Under Assumption~\ref{assumption:local_is}, given $\theta\in (0,1]$, let $\left(\u,p\right)\in \bm H^{m}\left(\Omega\right)\times H^{r}\left(\Omega\right)$, $m\ge 2$ and $r\ge 1$, be the solution to problem~\eqref{stokes}. Then, the discrete solution $\left(\u_h,p_h\right)\in  V_h^{\square,\u_D}\times  \overline Q_h^\square$ of the stabilized problem~\eqref{stabilized_formulation} satisfies
	\begin{equation*}
		\norm{\u-\u_h}_{1,h}+\norm{p-p_h}_{0,h} \le C h^{\min\{s,\ell\}} \left(\norm{\u}_{H^{m}(\Omega)} + \norm{p}_{H^{r}(\Omega)}  \right),
	\end{equation*}
	where $s:=\min \{k,m-1\}$ if $\square=\mathrm{RT}$, $s:=\min \{k+1,m-1\}$ if $\square\in\{\mathrm{N},\mathrm{RT}\}$, $\ell:=\min\{k+1,r\}$, and $C>0$ depends on the constants appearing in Lemma~\ref{lemma:apriori}.
\end{theorem}	
\begin{proof}
It is sufficient to proceed by triangular inequality and apply Proposition~\ref{inter_err_est}, and Lemma~\ref{lemma_approx}.
\end{proof}	
\IfStandalone
{
	\bibliographystyle{../elsarticle-num}
	\bibliography{../bibliographyx}
}{}

%% file: overlap1.tex
	\begin{tikzpicture}
	
\definecolor{color1}{rgb}{0,0,1} 
\definecolor{color2}{rgb}{1,0.019607843137255,0.0819607843137255} 
\definecolor{color3}{rgb}{0.019607843137255,0.872549019607843,0.0819607843137255} 

\coordinate (A0) at (0,-1);
\coordinate (B0) at (8,-1);
\coordinate (C0) at (8,3);
\coordinate (D0) at (0,3);

\coordinate (A1) at (4.37,5.01);
\coordinate (B1) at (7.09,3.77);
\coordinate (C1) at (9.32,8.67);
\coordinate (D1) at (6.61,9.91);

\coordinate (A2) at (1.31,4.6);
\coordinate (B2) at  (3.19,5.19);
\coordinate (C2) at  (2.15,8.52);
\coordinate (D2) at (0.27,7.93);

\draw[thick, fill = white] (A0) -- (B0) -- (C0) -- (D0) -- (A0); 

\draw[gray]
\foreach \i in {1, ..., 4} {
	\foreach \j in {1, ..., 4} {
		($(A0)!\i/5!(B0)$) -- ($(D0)!\i/5!(C0)$)
		($(B0)!\i/5!(C0)$) -- ($(A0)!\i/5!(D0)$)
	}
}
; 

\draw[thick,fill=white] (A1) -- (B1) -- (C1) -- (D1) -- (A1); 

\draw[gray]
\foreach \i in {1, ..., 4} {
	\foreach \j in {1, ..., 4} {
		($(A1)!\i/5!(B1)$) -- ($(D1)!\i/5!(C1)$)
		($(B1)!\i/5!(C1)$) -- ($(A1)!\i/5!(D1)$)
	}
}
; 

\draw[thick,fill=white] (A2) -- (B2) -- (C2) -- (D2) -- (A2); 

\draw[gray]
\foreach \i in {1, ..., 4} {
	\foreach \j in {1, ..., 4} {
		($(A2)!\i/5!(B2)$) -- ($(D2)!\i/5!(C2)$)
		($(B2)!\i/5!(C2)$) -- ($(A2)!\i/5!(D2)$)
	}
}
; 
	
	\coordinate (a) at (3.75,1);
	\coordinate (b) at (7,6.75);
	\coordinate (c) at (1.75,6.6);
	
\node at (a)   {\huge$\mathbf{\Omega_0^\ast}$};
\node at (b)   {\huge$\mathbf{\Omega_1^\ast}$};
\node at (c)   {\huge$\mathbf{\Omega_2^\ast}$};

	\coordinate (a) at (2,2);
	\begin{scope}
	\clip (0,0) -- (9.5,0) -- (9.5,10) -- (0,10) -- (0,0);
	\end{scope}

	\end{tikzpicture}
	

%% file: overlap2.tex
\begin{tikzpicture}

\definecolor{color1}{rgb}{0,0,1} 
\definecolor{color2}{rgb}{1,0.019607843137255,0.0819607843137255} 
\definecolor{color3}{rgb}{0.019607843137255,0.872549019607843,0.0819607843137255} 

\coordinate (A0) at (0,-1);
\coordinate (B0) at (8,-1);
\coordinate (C0) at (8,3);
\coordinate (D0) at (0,3);

\coordinate (A1) at (2.91,2.37);
\coordinate (B1) at (5.63,1.13);
\coordinate (C1) at (7.86,6.03);
\coordinate (D1) at (5.15,7.27);

\coordinate (A2) at (2.26,1.68);
\coordinate (B2) at  (4.14,2.27);
\coordinate (C2) at (3.1,5.6);
\coordinate (D2) at (1.22,5.01);

\draw[thick, fill = white] (A0) -- (B0) -- (C0) -- (D0) -- (A0); 

\draw[gray]
\foreach \i in {1, ..., 4} {
	\foreach \j in {1, ..., 4} {
		($(A0)!\i/5!(B0)$) -- ($(D0)!\i/5!(C0)$)
		($(B0)!\i/5!(C0)$) -- ($(A0)!\i/5!(D0)$)
	}
}
; 

\draw[thick,fill=white] (A1) -- (B1) -- (C1) -- (D1) -- (A1); 

\draw[gray]
\foreach \i in {1, ..., 4} {
	\foreach \j in {1, ..., 4} {
		($(A1)!\i/5!(B1)$) -- ($(D1)!\i/5!(C1)$)
		($(B1)!\i/5!(C1)$) -- ($(A1)!\i/5!(D1)$)
	}
}
; 

\draw[thick,fill=white] (A2) -- (B2) -- (C2) -- (D2) -- (A2);  

\draw[gray]
\foreach \i in {1, ..., 4} {
	\foreach \j in {1, ..., 4} {
		($(A2)!\i/5!(B2)$) -- ($(D2)!\i/5!(C2)$)
		($(B2)!\i/5!(C2)$) -- ($(A2)!\i/5!(D2)$)
	}
}
; 

\draw[ultra thick,color2] (1.85,3) -- (2.26,1.68) -- (3.54,2.08); 
\draw[ultra thick,color1] (3.54,2.08) -- (5.63,1.13) -- (6.48,3); 
\draw[ultra thick,color3] (3.54,2.08) --  (4.14,2.27) -- (3.62,3.93); 

	\coordinate (a) at (3.75,1);
\coordinate (b) at (5.5,4.25);
\coordinate (c) at (2.75,3.85);
	\coordinate (d) at (2,1.5);
\coordinate (e) at (6.25,1);
\coordinate (f) at (4.5,3);

\node at (a)   {\huge$\mathbf{\Omega_0}$};
\node at (b)   {\huge$\mathbf{\Omega_1}$};
\node at (c)   {\huge$\mathbf{\Omega_2}$};

\node at (d)   {\textcolor{color2}{\Large$\mathbf{\Gamma_{20}}$}};
\node at (e)   {\textcolor{color1}{\Large$\mathbf{\Gamma_{10}}$}};
\node at (f)   {\textcolor{color3}{\Large$\mathbf{\Gamma_{21}}$}};

\coordinate (a) at (2,2);
\begin{scope}
\clip (0,0) -- (9.5,0) -- (9.5,10) -- (0,10) -- (0,0);
\end{scope}

\end{tikzpicture}

%% file: N_gamma_down.tex
\begin{tikzpicture}

\definecolor{color1}{rgb}{0,0,1} 
\definecolor{color2}{rgb}{1,0.019607843137255,0.0819607843137255} 
\definecolor{color3}{rgb}{0.019607843137255,0.872549019607843,0.0819607843137255} 

\draw[step=0.5cm,lightgray,thin] (1.5,-1) grid (2.5,5) ; 
\draw[ultra thick] (1.5,-1) -- (2.5,-1) -- (2.5,5) -- (1.5,5) -- (1.5,-1); 

\draw[step=0.72cm,lightgray,thin] (3.7,-2)  grid (5.7,5.5) ; 
\draw[ultra thick] (3.7,-2) -- (5.7,-2) -- (5.7,5.5) -- (3.7,5.5) -- (3.7,-2); 

\draw[step=0.36cm,lightgray,thin] (6.5,-1)  grid (7.5,5) ; 
\draw[ultra thick] (6.5,-1) -- (7.5,-1) -- (7.5,5) -- (6.5,5) -- (6.5,-1); 

\draw[fill=white] (0,0) -- (8,0) -- (8,4) -- (0,4) -- (0,0); 
\draw[step=1cm,lightgray,thin] (0,0) grid (8,4); 
\draw[ultra thick] (0,0) -- (8,0) -- (8,4) -- (0,4) -- (0,0); 

\coordinate (a) at (4,2);
\coordinate (b) at (2,-0.5);
\coordinate (c) at (4.7,-0.5);
\coordinate (d) at (7,-0.5);

\node at (a)   {\Large$\mathbf{\Omega_3}$};
\node at (b)   {\Large$\mathbf{\Omega_0}$};
\node at (c)   {\Large$\mathbf{\Omega_1}$};
\node at (d)   {\Large$\mathbf{\Omega_2}$};

\begin{scope}
\clip (0,0) -- (8,0) -- (8,4) -- (0,4) -- (0,0);
\end{scope}

\end{tikzpicture}

%% file: N_gamma_up.tex
\begin{tikzpicture}

\definecolor{color1}{rgb}{0,0,1} 
\definecolor{color2}{rgb}{1,0.019607843137255,0.0819607843137255} 
\definecolor{color3}{rgb}{0.019607843137255,0.872549019607843,0.0819607843137255} 

\draw[step=1cm,lightgray,thin] (0,0) grid (8,4); 
\draw[ultra thick] (0,0) -- (8,0) -- (8,4) -- (0,4) -- (0,0); 

\draw[fill = white] (1.5,-1) -- (2.5,-1) -- (2.5,5) -- (1.5,5) -- (1.5,-1); 
\draw[step=0.5cm,lightgray,thin] (1.5,-1) grid (2.5,5) ; 
\draw[ultra thick] (1.5,-1) -- (2.5,-1) -- (2.5,5) -- (1.5,5) -- (1.5,-1); 

\draw[fill = white] (3.7,-2) -- (5.7,-2) -- (5.7,5.5) -- (3.7,5.5) -- (3.7,-2); 
\draw[step=0.72cm,lightgray,thin] (3.7,-2)  grid (5.7,5.5) ; 
\draw[ultra thick] (3.7,-2) -- (5.7,-2) -- (5.7,5.5) -- (3.7,5.5) -- (3.7,-2); 

\draw[fill = white] (6.5,-1) -- (7.5,-1) -- (7.5,5) -- (6.5,5) -- (6.5,-1); 
\draw[step=0.36cm,lightgray,thin] (6.5,-1)  grid (7.5,5) ; 
\draw[ultra thick] (6.5,-1) -- (7.5,-1) -- (7.5,5) -- (6.5,5) -- (6.5,-1); 

\coordinate (a) at (3,2);
\coordinate (b) at (2,-0.5);
\coordinate (c) at (4.7,-0.5);
\coordinate (d) at (7,-0.5);

\node at (a)   {\Large$\mathbf{\Omega_0}$};
\node at (b)   {\Large$\mathbf{\Omega_1}$};
\node at (c)   {\Large$\mathbf{\Omega_2}$};
\node at (d)   {\Large$\mathbf{\Omega_3}$};

\begin{scope}
\clip (0,0) -- (8,0) -- (8,4) -- (0,4) -- (0,0);
\end{scope}

\end{tikzpicture}

%% file: N_O.tex
\begin{tikzpicture}

\definecolor{color1}{rgb}{0,0,1} 
\definecolor{color2}{rgb}{1,0.019607843137255,0.0819607843137255} 
\definecolor{color3}{rgb}{0.019607843137255,0.872549019607843,0.0819607843137255} 

\draw[step=0.5cm,lightgray,thin] (0,-2) grid (8,5.5) ; 
\draw[ultra thick] (0,-2) -- (8,-2) -- (8,5.5) -- (0,5.5) -- (0,-2); 

\draw[fill=white] (1,-1) -- (7,-1) -- (7,4.5) -- (1,4.5) -- (1,-1); 
\draw[step=0.51cm,lightgray,thin] (1,-1)  grid (7,4.5) ; 
\draw[ultra thick] (1,-1) -- (7,-1) -- (7,4.5) -- (1,4.5) -- (1,-1); 

\draw[fill=white] (2,0) -- (6,0) -- (6,3.5) -- (2,3.5) -- (2,0); 
\draw[step=0.333cm,lightgray,thin] (2,0)  grid (6,3.5) ; 
\draw[ultra thick] (2,0) -- (6,0) -- (6,3.5) -- (2,3.5) -- (2,0); 

\draw[fill=white] (3,1) -- (5,1) -- (5,2.5) -- (3,2.5) -- (3,1); 
\draw[step=0.5cm,lightgray,thin] (3,1) grid (5,2.5); 
\draw[ultra thick] (3,1) -- (5,1) -- (5,2.5) -- (3,2.5) -- (3,1); 

\coordinate (a) at (4,1.75);
\coordinate (b) at (4,0.5);
\coordinate (c) at (4,-0.5);
\coordinate (d) at (4,-1.5);

\node at (a)   {\Large$\mathbf{\Omega_3}$};
\node at (b)   {\Large$\mathbf{\Omega_2}$};
\node at (c)   {\Large$\mathbf{\Omega_1}$};
\node at (d)   {\Large$\mathbf{\Omega_0}$};

\begin{scope}
\clip (0,0) -- (8,0) -- (8,4) -- (0,4) -- (0,0);
\end{scope}

\end{tikzpicture}

%% file: small_intersection_j_i_bis.tex
\begin{tikzpicture}

\definecolor{color1}{rgb}{0,0,1} 
\definecolor{color2}{rgb}{1,0.019607843137255,0.0819607843137255} 
\definecolor{color3}{rgb}{0.019607843137255,0.872549019607843,0.0819607843137255} 
\draw[step=1cm,lightgray,thin] (0,-1) grid (8,3); 
\draw[ultra thick] (0,-1) -- (8,-1) -- (8,3) -- (0,3) -- (0,-1); 

\coordinate (A) at (2.74,2.62);
\coordinate (B) at (5.46,1.38);
\coordinate (C) at (7.69,6.28);
\coordinate (D) at (4.98,7.52);

\coordinate (E) at (5.74,2);
\coordinate (F) at (6,2);
\coordinate (G) at (5.99,2.55);

\coordinate (H) at (2.66,2.1);
\coordinate (I) at (4.54,2.7);
\coordinate (L) at (3.5,6.03);
\coordinate (M) at (1.62,5.44);


\draw[ultra thick, fill = white] (A)-- (B) -- (C) -- (D) -- (A); 

\draw[lightgray]
\foreach \i in {1, ..., 4} {
	\foreach \j in {1, ..., 4} {
		($(A)!\i/5!(B)$) -- ($(D)!\i/5!(C)$)
		($(B)!\i/5!(C)$) -- ($(A)!\i/5!(D)$)
	}
}
; 



\coordinate (N) at (3.38,2.33);
\coordinate (O) at (3.85,2.12);
\coordinate (P) at (4.04,2.54);


\draw[ultra thick, blue] (2.92,3)--(A) --(N)-- (B)--(6.2,3) ; 

\coordinate (a) at (4,1);
\coordinate (b) at (5.3,4.25);
\coordinate (c) at (3,4.25);
\coordinate (d) at (5.99,1.2);

\node at (a)   {\huge$\mathbf{\Omega_j}$};
\node at (b)   {\huge$\mathbf{\Omega_i}$};
\node at (d)   {\textcolor{blue}{\Large$\mathbf{\Gamma_{ij}}$}};

\begin{scope}
\clip (0,0) -- (8,0) -- (8,7.5) -- (0,7.5) -- (0,0);
\end{scope}

\end{tikzpicture}

%% file: small_intersection_j_i.tex
\begin{tikzpicture}

\definecolor{color1}{rgb}{0,0,1} 
\definecolor{color2}{rgb}{1,0.019607843137255,0.0819607843137255} 
\definecolor{color3}{rgb}{0.019607843137255,0.872549019607843,0.0819607843137255} 
\draw[step=1cm,lightgray,thin] (0,-1) grid (8,3); 
\draw[ultra thick] (0,-1) -- (8,-1) -- (8,3) -- (0,3) -- (0,-1); 

\coordinate (A) at (2.74,2.62);
\coordinate (B) at (5.46,1.38);
\coordinate (C) at (7.69,6.28);
\coordinate (D) at (4.98,7.52);

\coordinate (E) at (5.74,2);
\coordinate (F) at (6,2);
\coordinate (G) at (5.99,2.55);

\coordinate (H) at (2.66,2.1);
\coordinate (I) at (4.54,2.7);
\coordinate (L) at (3.5,6.03);
\coordinate (M) at (1.62,5.44);


\draw[ultra thick, fill = white] (A)-- (B) -- (C) -- (D) -- (A); 

\draw[lightgray]
\foreach \i in {1, ..., 4} {
	\foreach \j in {1, ..., 4} {
		($(A)!\i/5!(B)$) -- ($(D)!\i/5!(C)$)
		($(B)!\i/5!(C)$) -- ($(A)!\i/5!(D)$)
	}
}
; 

\draw[ultra thick, fill = white] (H)-- (I) -- (L) -- (M) -- (H); 

\draw[lightgray]
\foreach \m in {1, ..., 4} {
	\foreach \n in {1, ..., 4} {
		($(H)!\m/5!(I)$) -- ($(M)!\m/5!(L)$)
		($(I)!\m/5!(L)$) -- ($(H)!\m/5!(M)$)
	}
}
;

\coordinate (N) at (3.38,2.33);
\coordinate (O) at (3.85,2.12);
\coordinate (P) at (4.04,2.54);

\draw[ fill =red] (N)-- (O) -- (P) -- (N); 

\draw[ultra thick, blue] (N)-- (B) -- (6.2,3); 

\coordinate (a) at (4,1);
\coordinate (b) at (5.3,4.25);
\coordinate (c) at (3,4.25);
\coordinate (d) at (5.99,1.2);

\node at (a)   {\huge$\mathbf{\Omega_j}$};
\node at (b)   {\huge$\mathbf{\Omega_i}$};
\node at (c)   {\huge$\mathbf{\Omega_k}$};
\node at (d)   {\textcolor{blue}{\Large$\mathbf{\Gamma_{ij}}$}};

\begin{scope}
\clip (0,0) -- (8,0) -- (8,7.5) -- (0,7.5) -- (0,0);
\end{scope}

\end{tikzpicture}

%% file: good_bad_choice_multi1.tex
	\begin{tikzpicture}
	
	\definecolor{color2}{rgb}{0.872549019607843,0.019607843137255,0.0819607843137255}
	
	\draw[step=1cm,lightgray,very thin] (4,2) grid (9,9); 
	\draw[ultra thick] (4,8) -- (4,9) -- (9,9) -- (9,2) -- (4,2) -- (4,3); 
	\filldraw[fill=red!80!white, draw=black] (6.5,5) rectangle (7,6); 
    \filldraw[fill=green!60!white, draw=black] (7,5) rectangle (8,6); 
    	\draw[ultra thick,fill=white] (6,3.5) -- (6.5,3.5) -- (6.5,7.5) -- (6,7.5) -- (6,3.5); 
    \draw[step=1.09cm,lightgray,very thin] (6,3.5) grid (6.5,7.5); 
	\draw[ultra thick] (6,3.5) -- (6.5,3.5) -- (6.5,7.5) -- (6,7.5) -- (6,3.5); 
		\draw[ultra thick,fill=white] (0,3) -- (6,3) -- (6,8) -- (0,8) -- (0,3); 
		\draw[step=2cm,lightgray,very thin] (0,3) grid (6,8); 
	\draw[ultra thick] (0,3) -- (6,3) -- (6,8) -- (0,8) -- (0,3); 
	
	\coordinate (a) at (6.75,5.5);
	\coordinate (b) at (7.5,5.5);
	\coordinate (d) at (2.5,7.5);
	\coordinate (c) at (2,4);
	
	\node at (a)   {\large$\mathbf{K}$};
	\node at (b)   {\large$\mathbf{K'}$};
	
	\coordinate (a) at (2,2);
\begin{scope}
\clip (1,2) -- (9,2) -- (9,9) -- (0,9) -- (1,2);
\end{scope}

	\end{tikzpicture}
	

%% file: good_bad_choice_multi.tex
\begin{tikzpicture}
	
\definecolor{color2}{rgb}{0.872549019607843,0.019607843137255,0.0819607843137255}

\draw[ultra thick,fill=white] (4,8) -- (4,9) -- (9,9) -- (9,2) -- (4,2) -- (4,3); 
\draw[step=0.9cm,lightgray,very thin] (4,2) grid (9,9); 
\draw[ultra thick] (4,8) -- (4,9) -- (9,9) -- (9,2) -- (4,2) -- (4,3); 
\draw[ultra thick,fill=white] (6,3.5) -- (6.5,3.5) -- (6.5,7.5) -- (6,7.5) -- (6,3.5); 
\filldraw[fill=red!80!white, draw=black] (6,4) rectangle (6.48,5); 
\draw[step=1cm,lightgray,very thin] (6,3.5) grid (6.5,7.5); 
\draw[ultra thick] (6,3.5) -- (6.5,3.5) -- (6.5,7.5) -- (6,7.5) -- (6,3.5); 
\draw[ultra thick,fill=white] (0,3) -- (6,3) -- (6,8) -- (0,8) -- (0,3); 
\filldraw[fill=green!60!white, draw=black] (4,4) rectangle (5.98,6); 
\draw[step=2cm,lightgray,very thin] (0,3) grid (6,8); 
\draw[ultra thick] (0,3) -- (6,3) -- (6,8) -- (0,8) -- (0,3); 

\coordinate (a) at (5,5);
\coordinate (b) at (6.25,4.5);

\node at (a)   {\large$\mathbf{K'}$};
\node at (b)   {\large$\mathbf{K}$};

\begin{scope}
  \clip (1,2) -- (9,2) -- (9,9) -- (0,9) -- (1,2);
\end{scope}

\end{tikzpicture}

%% file: example.tex

\section{Numerical examples}
\label{sec:example}

In this section, we present several numerical examples to demonstrate the convergence, accuracy, conditioning, and efficiency of the proposed method. The viscosity coefficient is $\mu=1$. Taylor-Hood elements are used for the discretization. The threshold of $\theta=0.1$ is used to identify bad elements. The creation of an interface quadrature mesh follows exactly the same as~\cite{ref:wei21a}, where the key is to find a mesh intersection on the boundary of the top patch such that each quadrature cell only involves polynomials. The reparameterization of cut elements for numerical integration follows that of~\cite{ref:wei21b}, where each cut element is decomposed into a collection of primitive cells (i.e., high-order triangles and/or quadrilaterals).

\subsection{Two-patch union}
\label{sec:2patch}

\begin{figure}[!ht]
\centering
\includegraphics[width=0.45\textwidth]{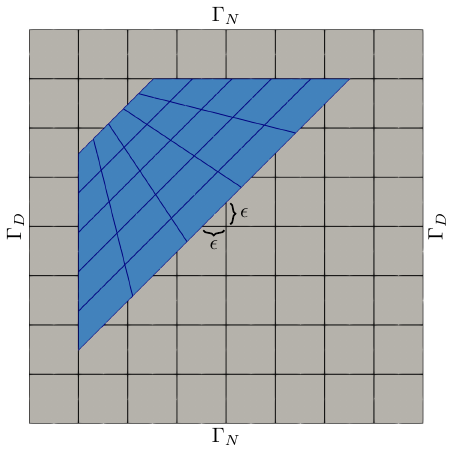}
\caption{The initial mesh and boundary conditions of the two-patch union.}
\label{fig:two_in}
\end{figure}

We start with a two-patch union that can control the trimming condition via a single parameter $\epsilon$. As the input, the bottom (B-spline) patch is simply a unit square with a grid of $8\times 8$ elements, whereas the top patch is a trapezoid with $5\times 5$ elements; see Figure~\ref{fig:two_in}. Note that the horizontal and vertical boundaries of the top patch align with certain mesh lines of the bottom patch. As $\epsilon$ decreases, the effective area of cut elements shrinks, i.e., the trimming condition deteriorates. 

We adopt the following manufactured solution for the convergence study,
$$
\mathbf{u}=
\begin{bmatrix}
2 e^x (-1+x)^2 x^2 (y^2-y) (-1+2y) \\
-e^x (-1+x) x (-2+3x+x^2) (-1+y)^2 y^2
\end{bmatrix},
$$
and
$$
\begin{aligned}
p=& -424+156e+(y^2-y) (-456 + e^x (456+x^2 (228 - 5(y^2-y)) \\
  & +2x(-228+y^2-y) +2x^3(-36+y^2-y) +x^4(12+y^2-y)  ) ).
\end{aligned}
$$
Boundary conditions are imposed on the bottom patch. The Dirichlet boundary condition ($\mathbf{u}=\bm{0}$) is imposed on the left and right boundaries, whereas the Neumann boundary condition is imposed on the top and bottom boundaries. 

\begin{figure}[!ht]
\centering
\begin{tabular}{cc}
\includegraphics[width=0.45\textwidth]{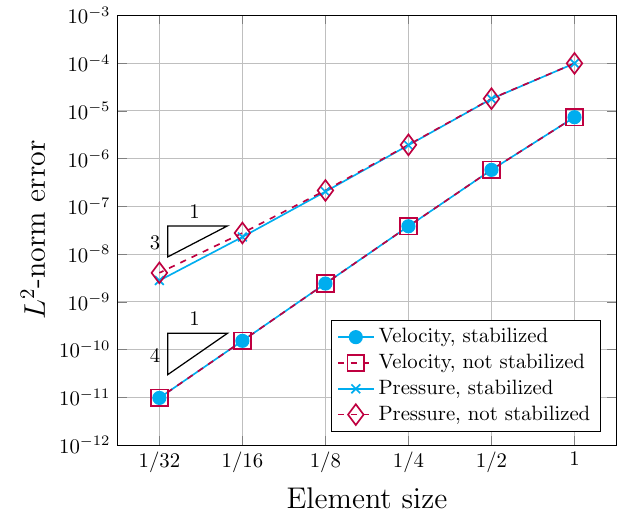}&
\includegraphics[width=0.45\textwidth]{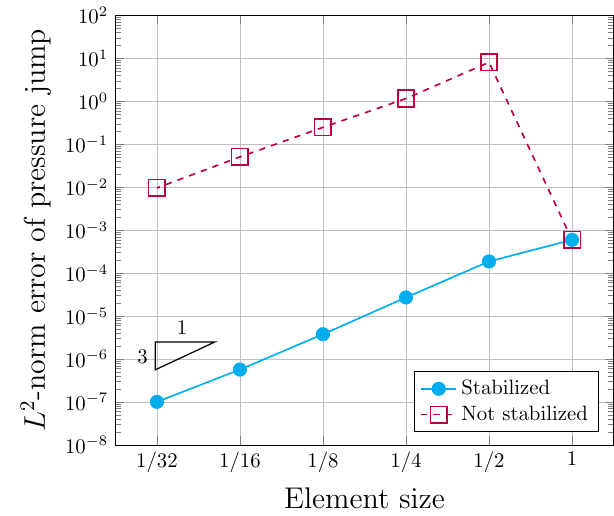}\\
(a) & (b)\\
\includegraphics[width=0.45\textwidth]{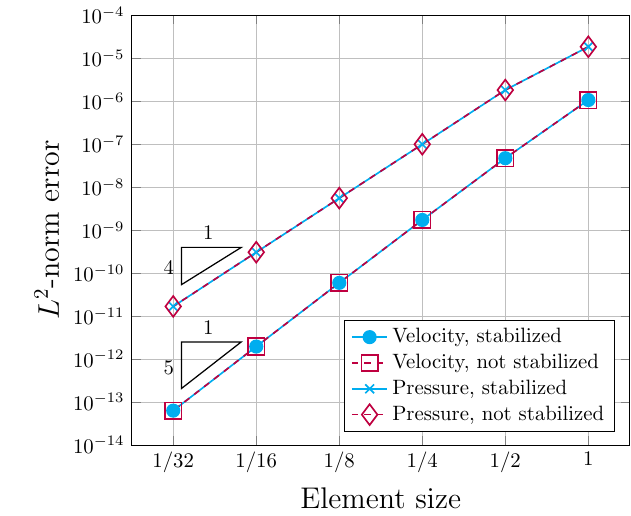}&
\includegraphics[width=0.45\textwidth]{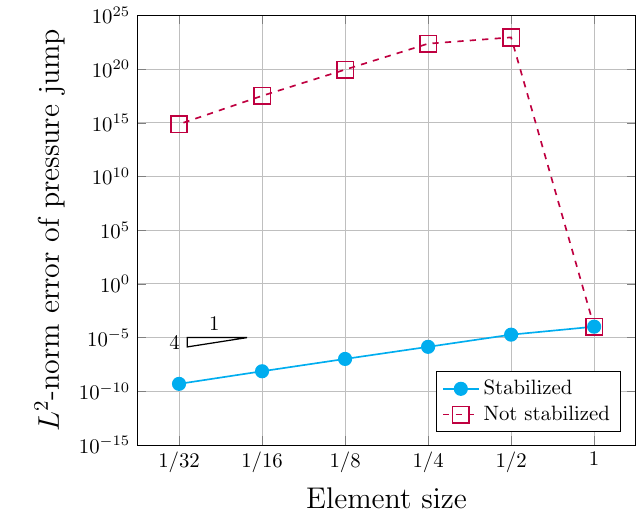}\\
(c) & (d)\\
\end{tabular}
\caption{Convergence plots of the two-patch union with $\epsilon=10^{-12}$. (a, b) Quadratic spline discretization (for pressure), and (c, d) cubic spline discretization (for pressure).}
\label{fig:two_conv}
\end{figure}


We first check the convergence in an extremely trimmed case with $\epsilon=10^{-12}$. Quadratic and cubic spline discretizations (for pressure) are studied. A series of consecutively refined meshes is obtained for the convergence study by globally refining both top and bottom patches. The symmetric flux is adopted on the interface. The results are summarized in Figure~\ref{fig:two_conv}. {We compute two kinds of error: Figure~\ref{fig:two_conv}(a, c) corresponds to the error in the patch interior without considering the interface, whereas Figure~\ref{fig:two_conv}(b, d) is the error on the interface}. In Figure~\ref{fig:two_conv}(a, c), we observe that expected optimal convergence rates are achieved in all cases, regardless of whether bad elements are stabilized or not. In fact, the accuracy (in terms of the $L^2$-norm error) is almost indistinguishable. However, when looking at the accuracy of the pressure jump across the interface in Figure~\ref{fig:two_conv}(b, d), we find a significant difference of 5 orders (the quadratic case) or 26 orders (the cubic case). This means that without stabilization, the approximation error of pressure jump blows up at the interface. This behavior was not observed in the elliptic problems~\cite{ref:wei21a}. It signifies the necessity of stabilizing bad elements in general problems also for the sake of accuracy.

\begin{figure}[!ht]
\centering
\includegraphics[width=0.45\textwidth]{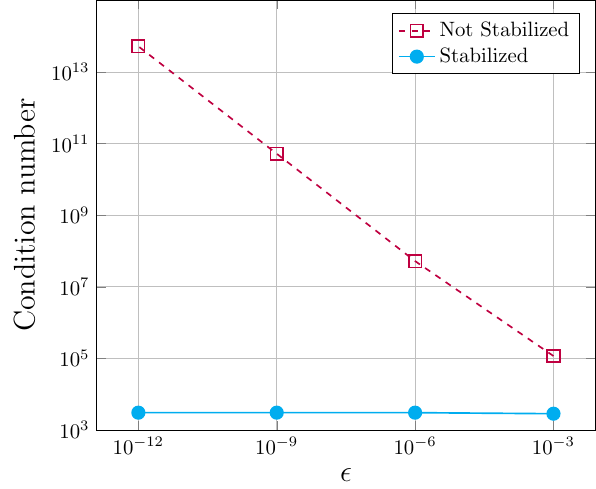}
\caption{Conditioning with respect to the trimming condition ($\epsilon$).}
\label{fig:two_cond}
\end{figure}

We next study the influence of stabilization on the conditioning of the resulting linear system. The test is carried out with a fixed mesh resolution but a varying trimming condition by changing $\epsilon$; see Figure~\ref{fig:two_in}. The symmetric flux is adopted. The quadratic discretization (for pressure) is used. Condition numbers are computed from preconditioned matrices, where we use the diagonal preconditioner. The results are summarized in Figure~\ref{fig:two_cond}. We observe that in the non-stabilized case, the condition numbers quickly blow up (in the order of $10^{5}\sim 10^{13}$) as $\epsilon$ decreases. In contrast, the condition numbers remain a constant ($\sim 10^3$) in the stabilized case regardless of the change of $\epsilon$, which numerically verifies the efficacy of the proposed stabilization method.

\rv{It is worth mentioning that regarding the conditioning issue, two types of cut elements are distinguished in a recent review article on stability and conditioning~\cite{dePrenter2023}: one that contains a single vertex (Case 1), and the other that contains no vertex at all (Case 2). According to the particular treatment of cut elements in~\cite{dePrenter2023}, the conditioning issue in Case 1 can be resolved by merely rescaling the stiffness matrix with a Jacobi preconditioner, whereas Case 2 needs both stabilization and rescaling. While our example here ``geometrically'' belongs to Case 1, it exhibits a fundamental difference due to the different treatment of cut elements. There is only one basis function in~\cite{dePrenter2023}, which is also the condition for rescaling only to resolve the conditioning issue, whereas there are always multiple basis functions in our example. Indeed, based on the results in Figure~\ref{fig:two_cond}, where the condition numbers are computed on rescaled matrices, we conclude that, unfortunately, rescaling only does not work in our case and stabilization is essential.}

\subsection{Multi-patch union}

\begin{figure}[!ht]
\centering
\begin{tabular}{cc}
\includegraphics[width=0.45\textwidth]{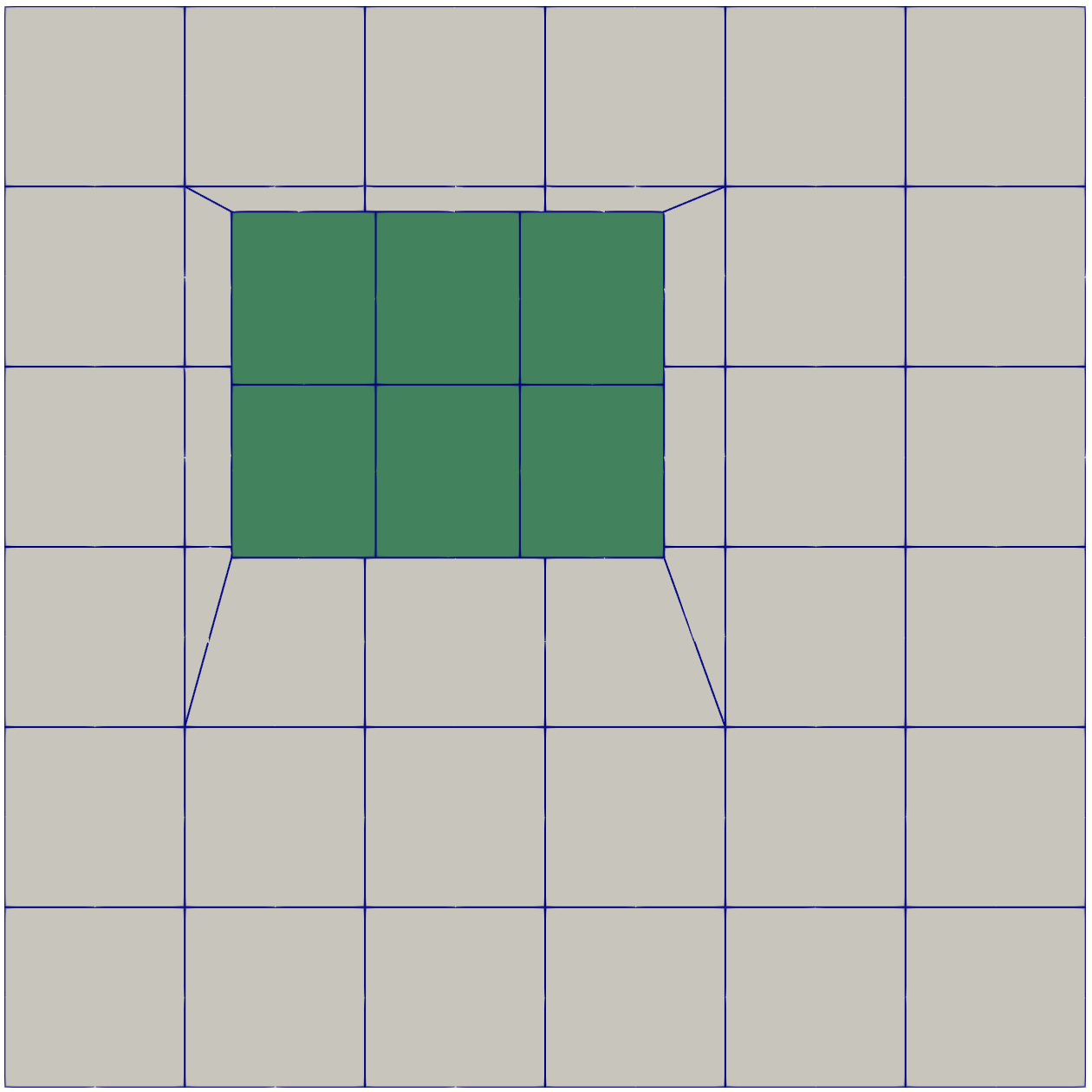}&
\includegraphics[width=0.45\textwidth]{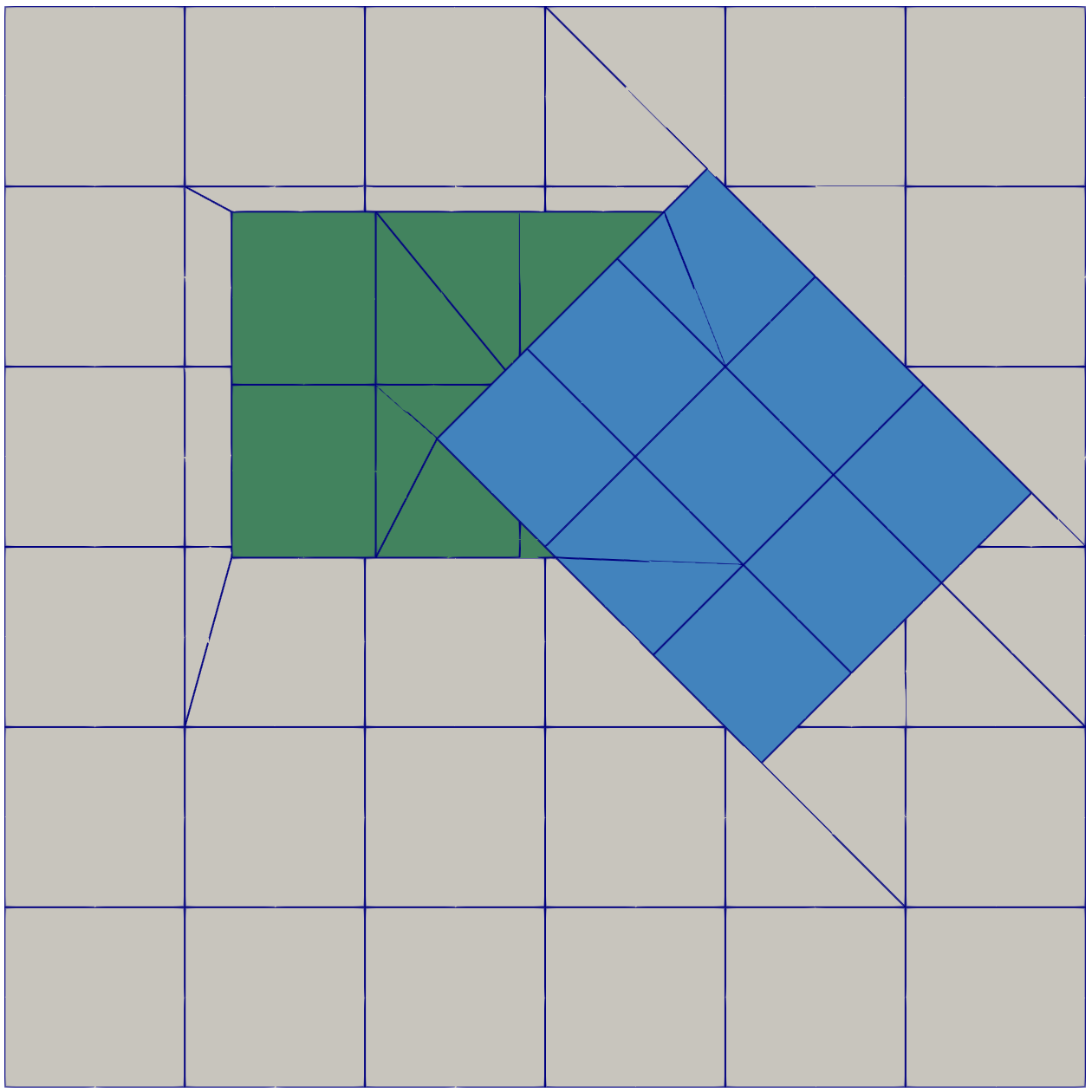}\\
(a) & (b)\\
\includegraphics[width=0.45\textwidth]{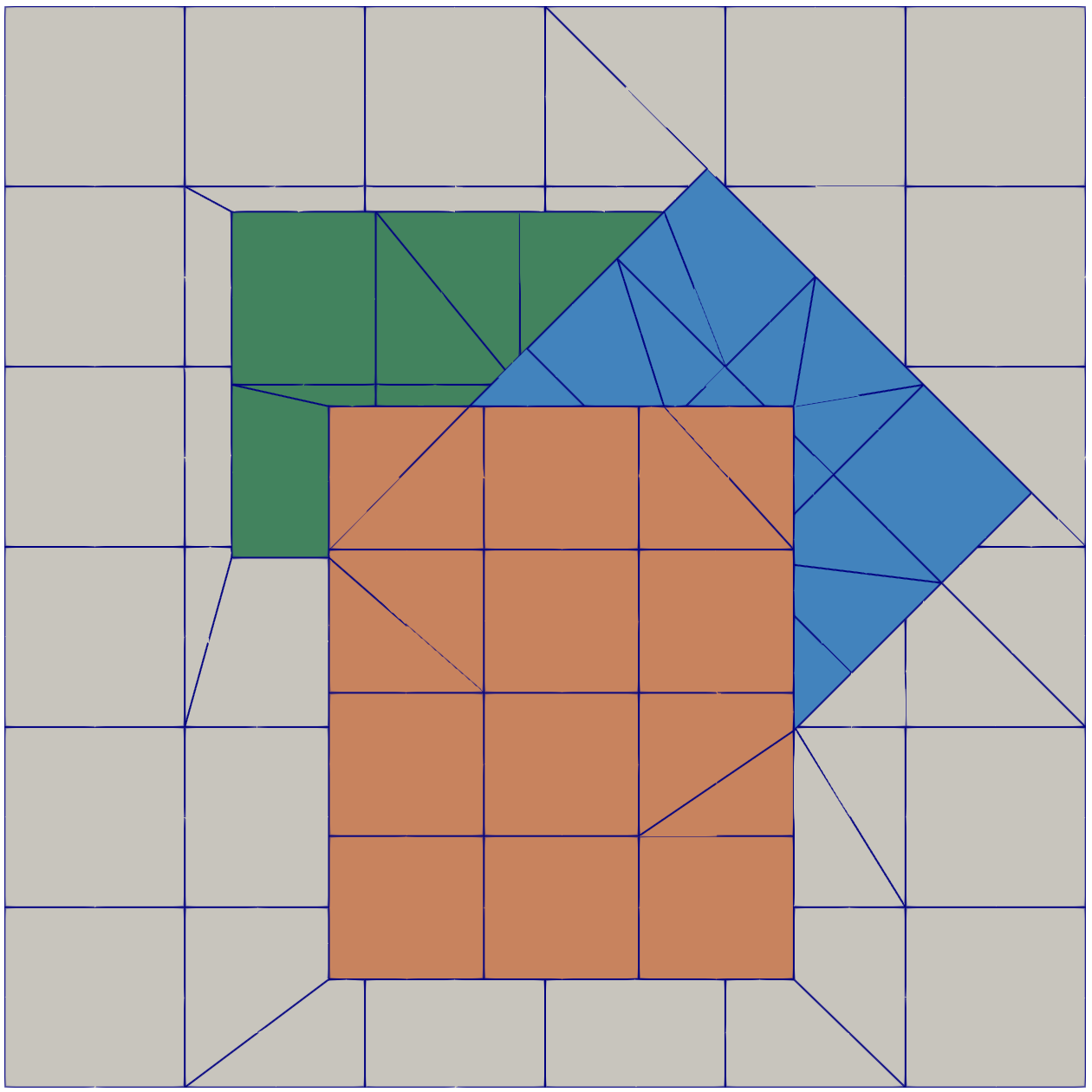}&
\includegraphics[width=0.45\textwidth]{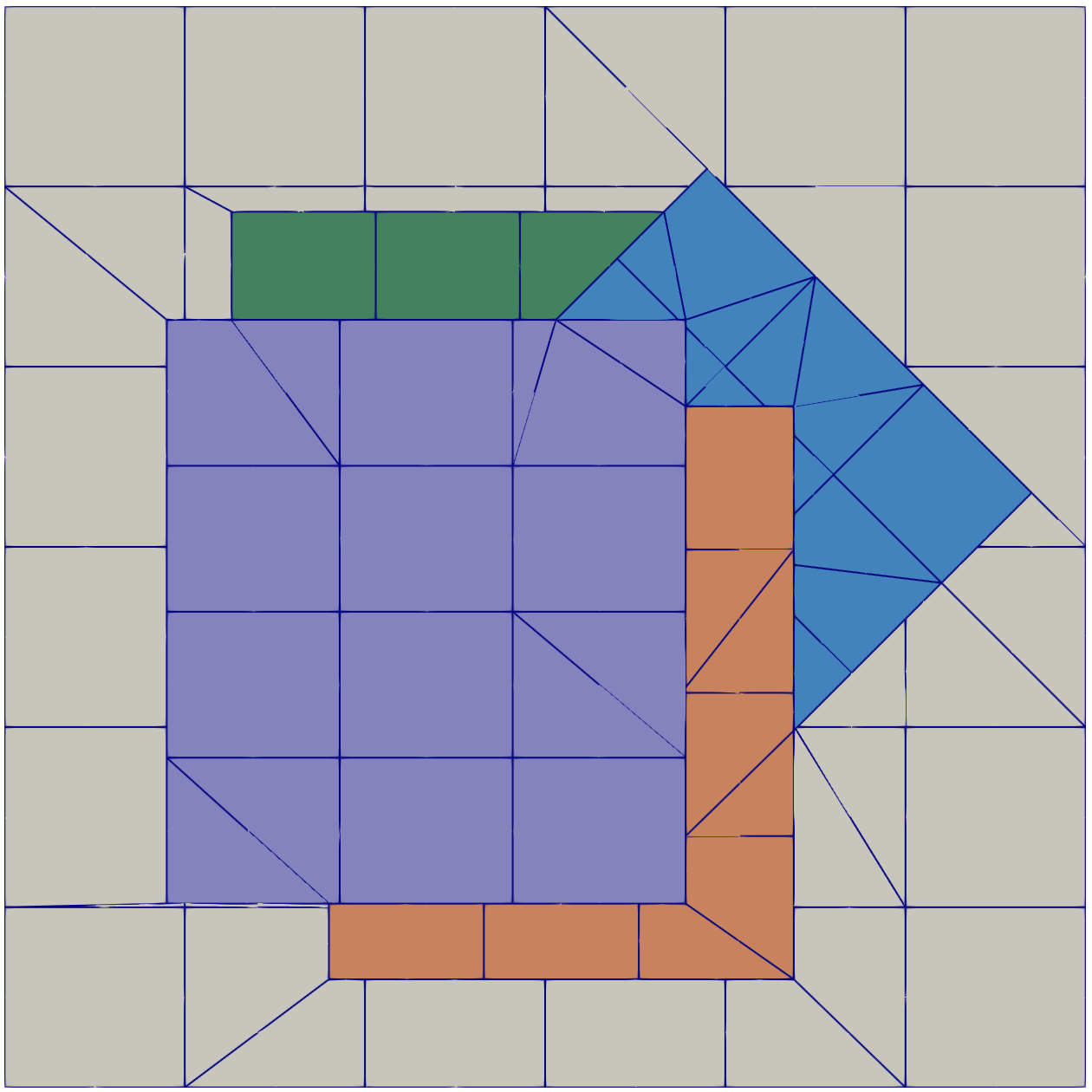}\\
(c) & (d)\\
\end{tabular}
\caption{Multiple overlapping patches: Initial meshes. Note that trimmed elements may be reparameterized into multiple integration cells.}
\label{fig:multi_in}
\end{figure}

Next we study whether the number of overlapping patches has a negative impact on accuracy. The overall geometry is simply a unit square, over which we gradually increase the number of overlapping patches from 2 to 5; see Figure~\ref{fig:multi_in}. To create a complex multi-patch overlapping scenario, each time we place a new patch in such a way that it overlaps with all the existing patches, which in practice becomes more and more difficult to achieve as the number of patches increases. All patches are intended to have a similar resolution. This way, the overall mesh, i.e., the collection of the visible meshes of all patches, can relatively retain a similar resolution even if new patches are added. 

\begin{figure}[!ht]
\centering
\begin{tabular}{cc}
\includegraphics[width=0.45\textwidth]{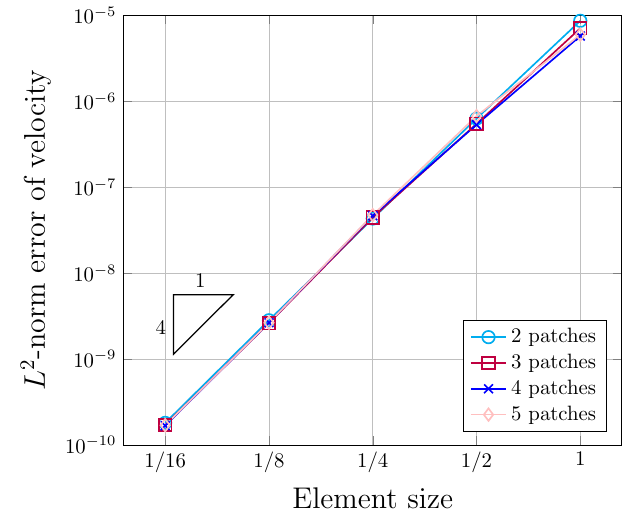}&
\includegraphics[width=0.45\textwidth]{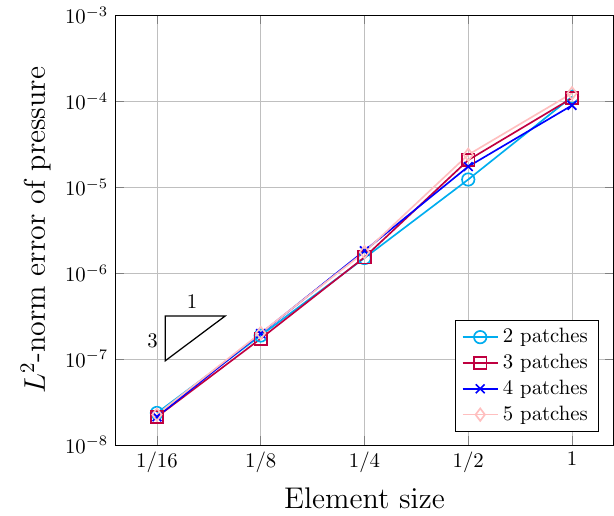}\\
(a) & (b)\\
\end{tabular}
\caption{Initial meshes of multiple overlapping patches, where the number of patches varies. Note that trimmed elements may be reparameterized into multiple integration cells.}
\label{fig:multi_conv}
\end{figure}

We compare the accuracy for meshes of different number of overlaps. Using the same manufactured solution in Section~\ref{sec:2patch}, we check the convergence as well as accuracy. The initial meshes are shown in Figure~\ref{fig:multi_in}. The convergence plots are summarized in Figure~\ref{fig:multi_conv}. We observe that expected optimal convergence is achieved in all cases. In fact, the accuracy of both velocity and pressure are almost indistinguishable, especially when the meshes are fine enough (e.g., $h\leq 1/4$ with $h$ being the element size). This means that the number of overlapping patches does not have an important effect on accuracy when the number is moderate. In practice, having a large number of overlaps (under the condition that each patch overlaps with all the patches below it) is quite rare to encounter.

\subsection{Airfoil}

\begin{figure}[!ht]
\centering
\begin{tabular}{c}
\includegraphics[width=0.6\textwidth]{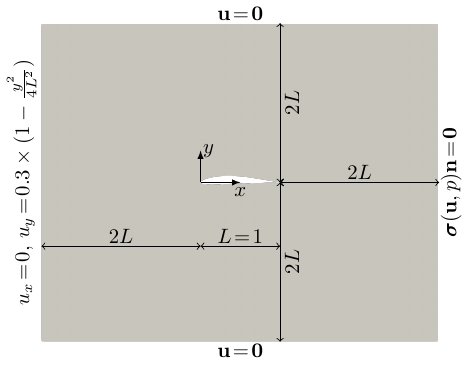}\\
(a) Setup
\end{tabular}
\begin{tabular}{cc}
\includegraphics[width=0.45\textwidth]{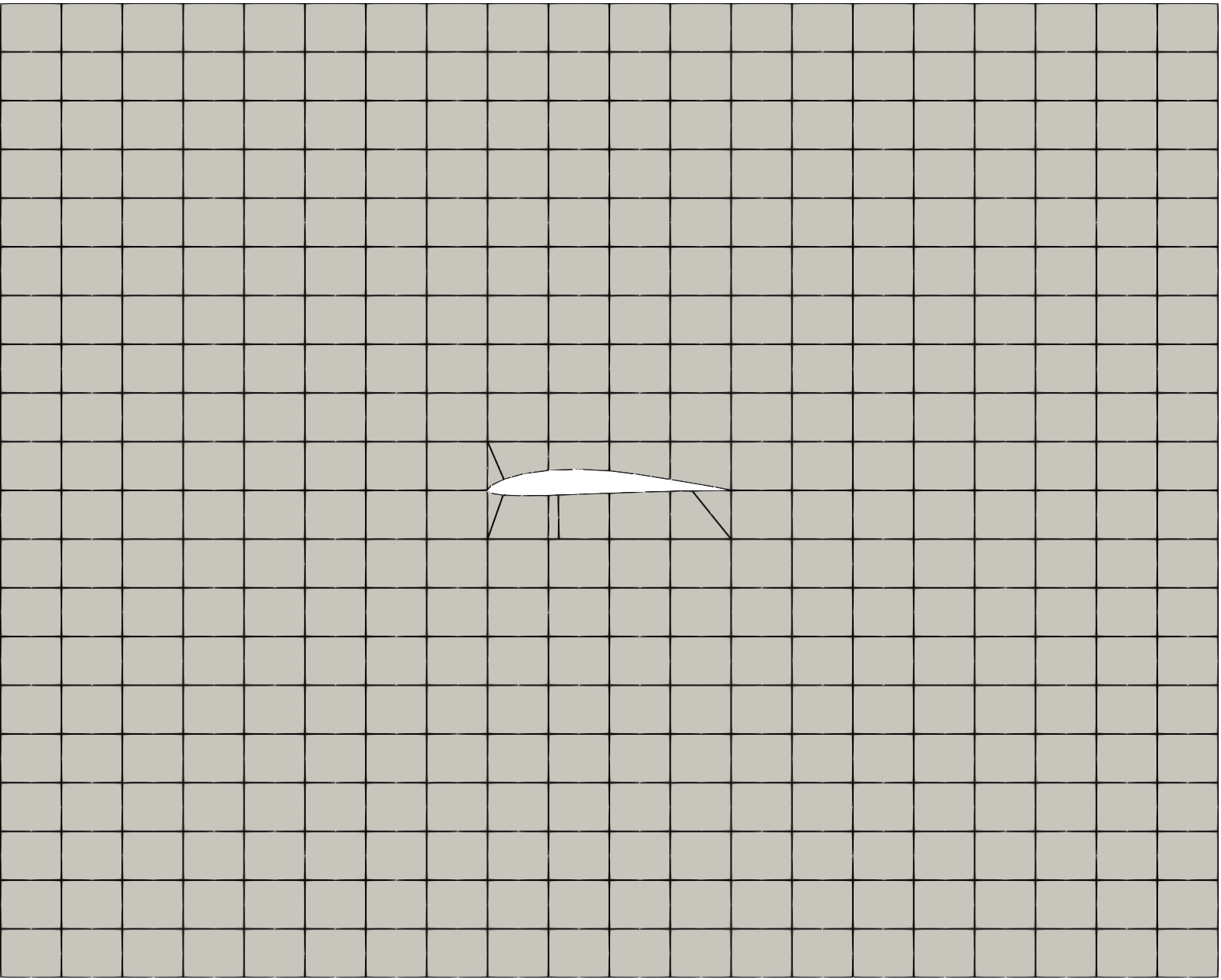} &
\includegraphics[width=0.45\textwidth]{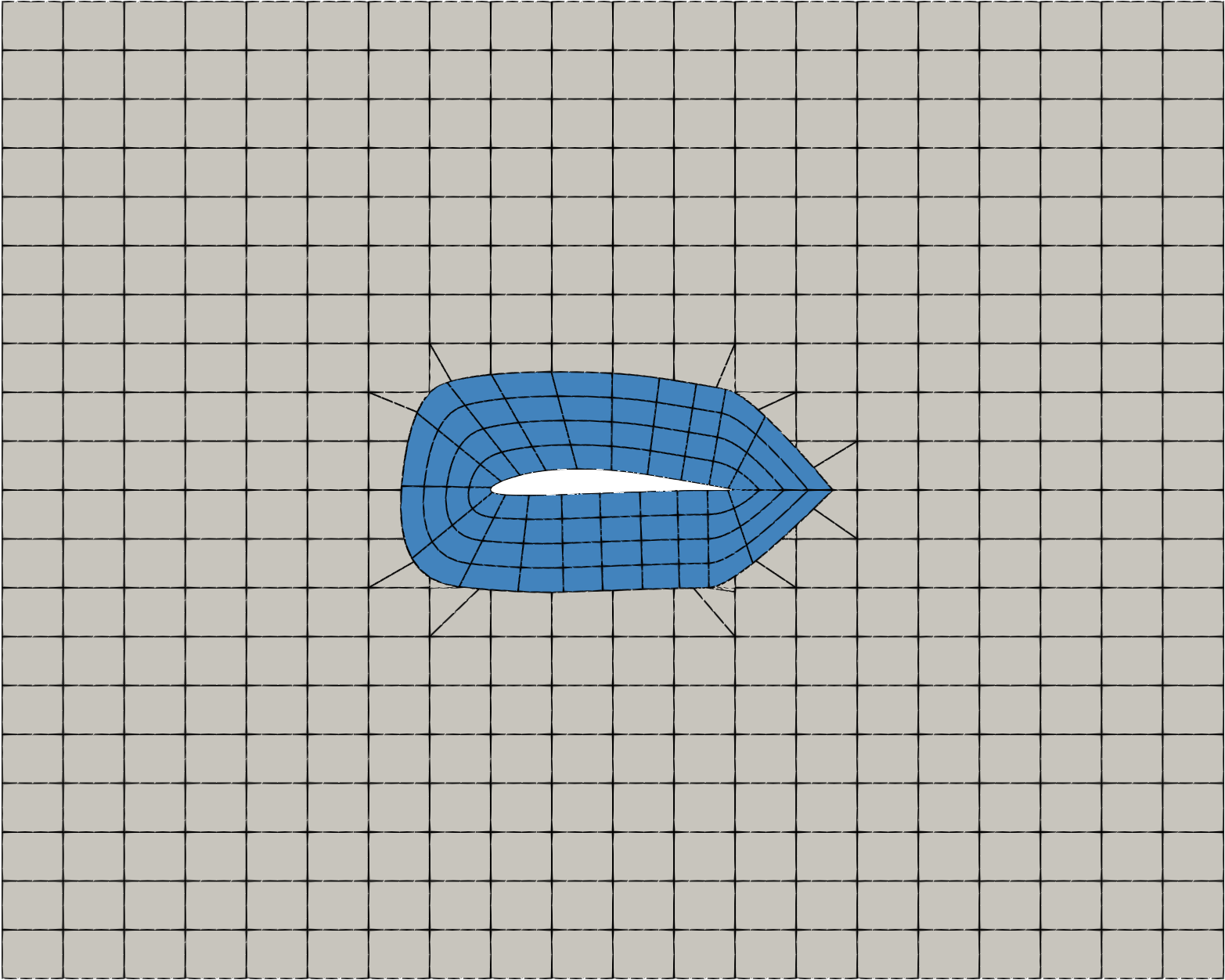} \\
(b) Trimming & (c) IBCM
\end{tabular}
\caption{The problem setup (a) and the two geometric constructions (b, c) of the airfoil. Compared to the geometric modeling via pure trimming (a), the immersed boundary-conformal method (IBCM) features an additional conformal mesh around the boundary (b). Note that trimmed elements are reparameterized into multiple integration cells.}
\label{fig:airfoil_input}
\end{figure}

In the last example, we solve the Stokes problem around an airfoil to show the flexibility and efficiency of the proposed method in capturing local features in the solution field. We compare two geometric constructions, one via straightforward trimming, and the other using the immersed boundary-conformal method (IBCM)~\cite{ref:wei21b}; see Figure~\ref{fig:airfoil_input}. IBCM is a special type of overlapping construction that features a layer of conformal mesh around geometric features (e.g., boundaries or interfaces), while keeping the majority of the computational domain described by the Cartesian grid. This way, IBCM represents a hybrid manner of geometric modeling that leverages the advantages of both the boundary-fitted method and the boundary-unfitted method.

The geometric modeling using IBCM for this example is introduced as follows. The initial data of the airfoil is a series of points, to which we fit a quadratic B-spline curve $\Gamma_{\mathrm{airfoil}}$. An offset B-spline curve $\Gamma_{\mathrm{offset}}$ is then obtained by offsetting the control polygon of $\Gamma_{\mathrm{airfoil}}$. This way we can readily have a layer of conformal mesh. We adopt such a simple way for offsetting because as already shown in~\cite{ref:wei21b}, the shape of the offset curve does not influence much the accuracy. Finally, $\Gamma_{\mathrm{offset}}$ is used to trim the background Cartesian grid. The computational domain consists of the conformal mesh and the trimmed Cartesian grid.

The boundary conditions are shown in Figure~\ref{fig:airfoil_input}(a). A parabolic Dirichlet condition is imposed on the left side, the homogeneous Dirichlet condition is imposed on the top and bottom boundaries as well as the airfoil, and the homogeneous Neumann condition is imposed on the right side. Quadratic spline discretization (for pressure) is used for this problem. We are particularly interested in the pressure field because it has a large gradient at the left tip of the airfoil and exhibits a singularity at the right tip.


\begin{figure}[!ht]
\centering
\begin{tabular}{cc}
\includegraphics[width=0.45\textwidth]{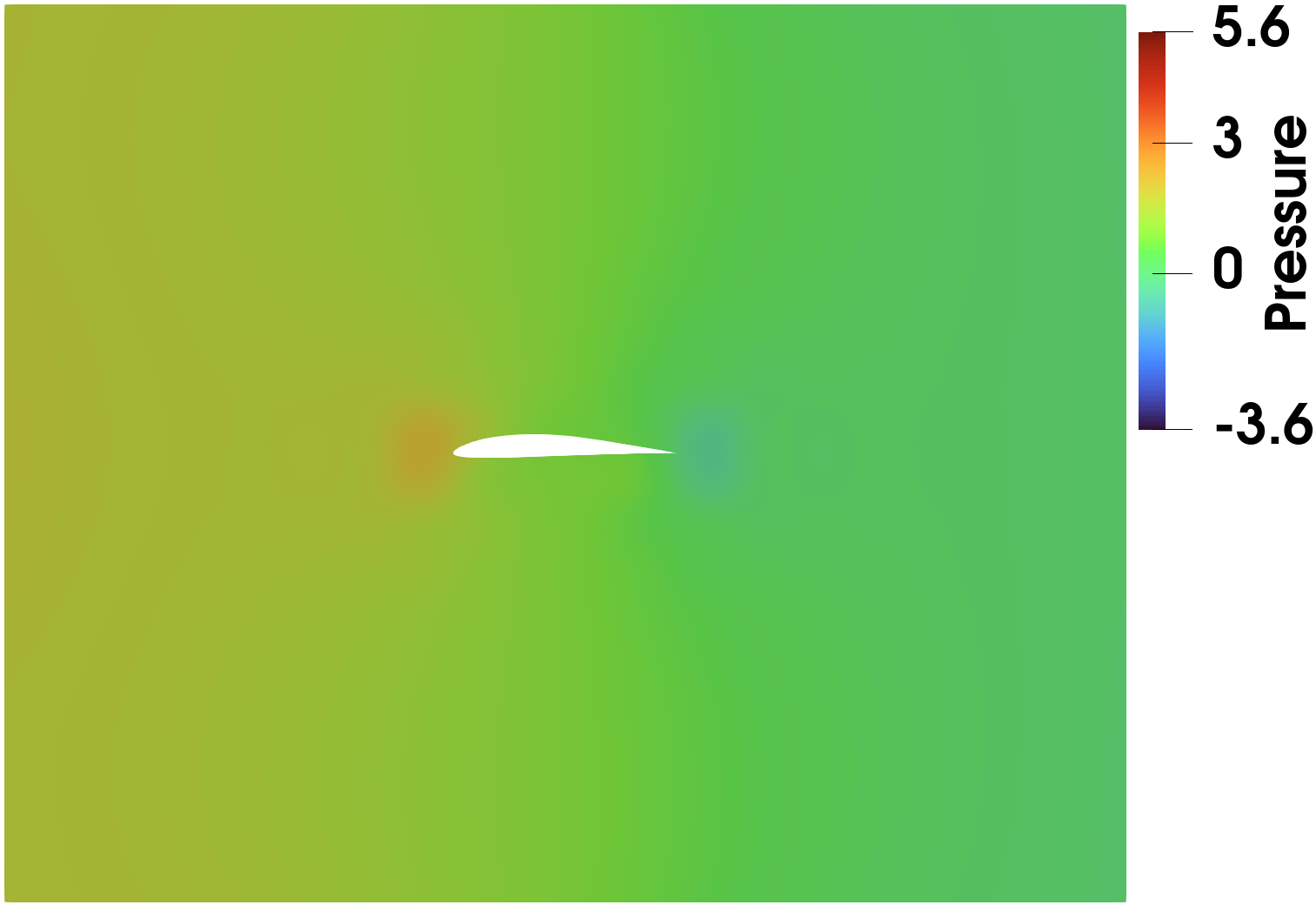} &
\includegraphics[width=0.45\textwidth]{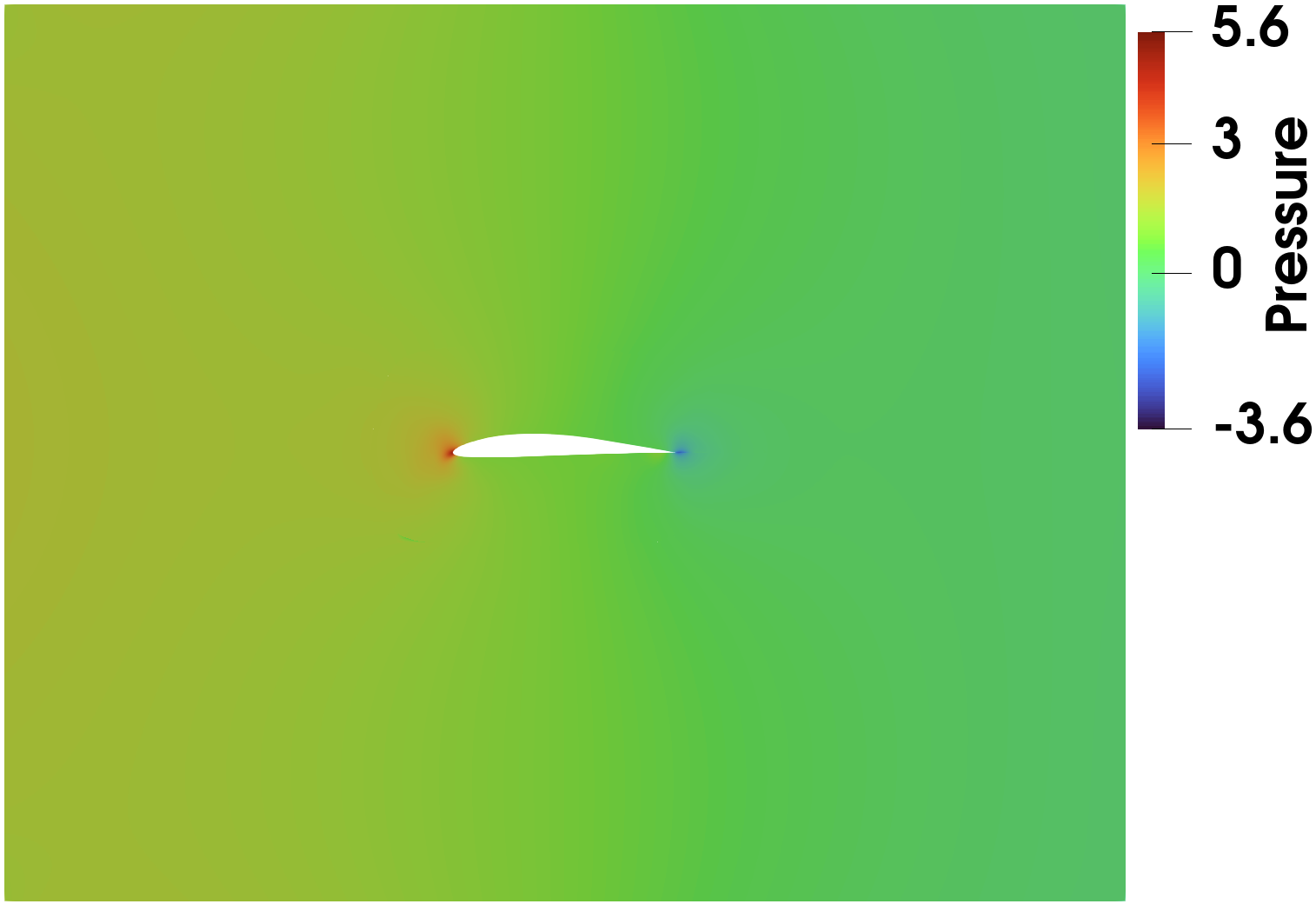} \\
(a) Trimming (initial mesh) & (b) IBCM (initial mesh) \\
\includegraphics[width=0.45\textwidth]{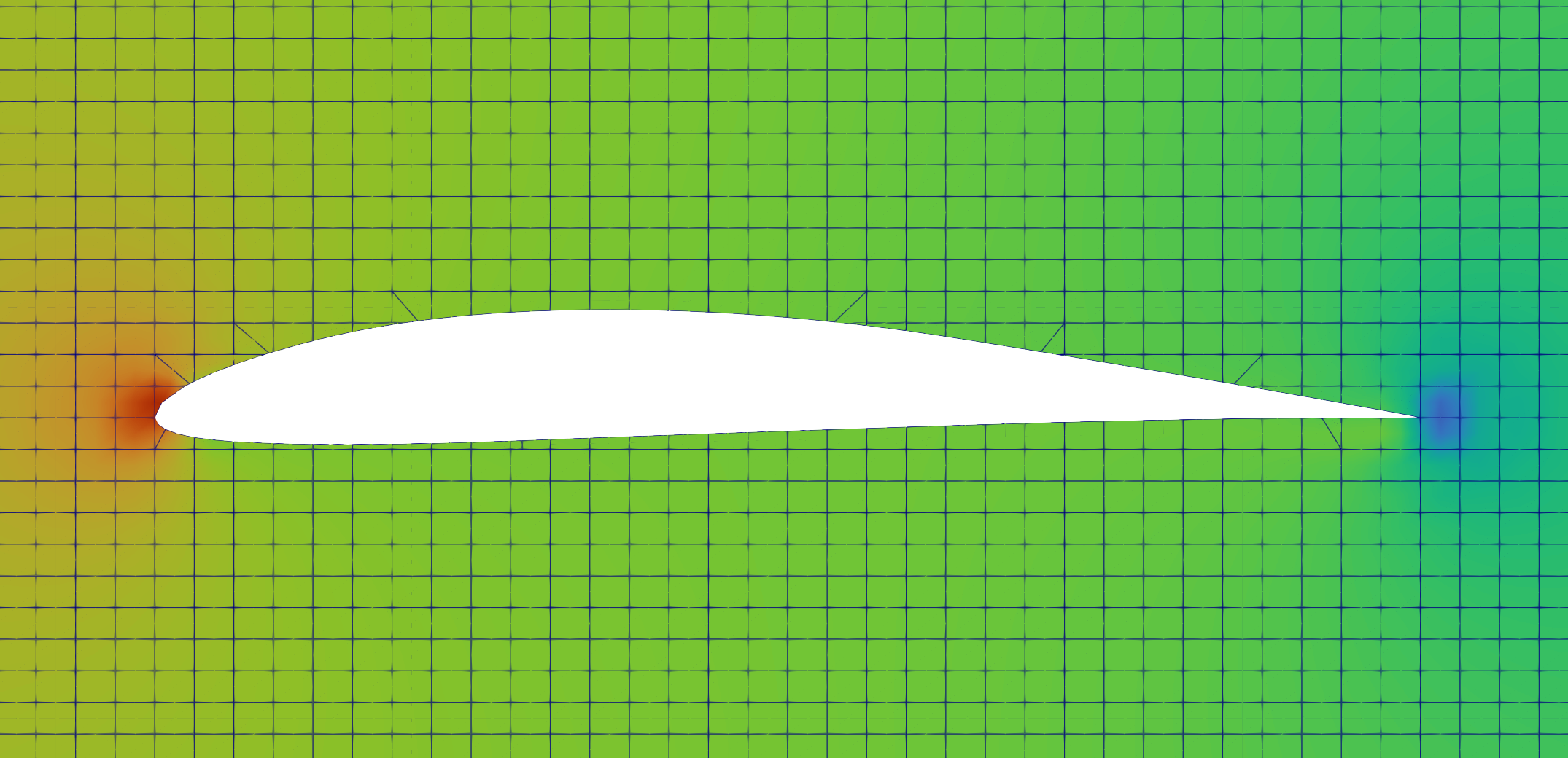} &
\includegraphics[width=0.45\textwidth]{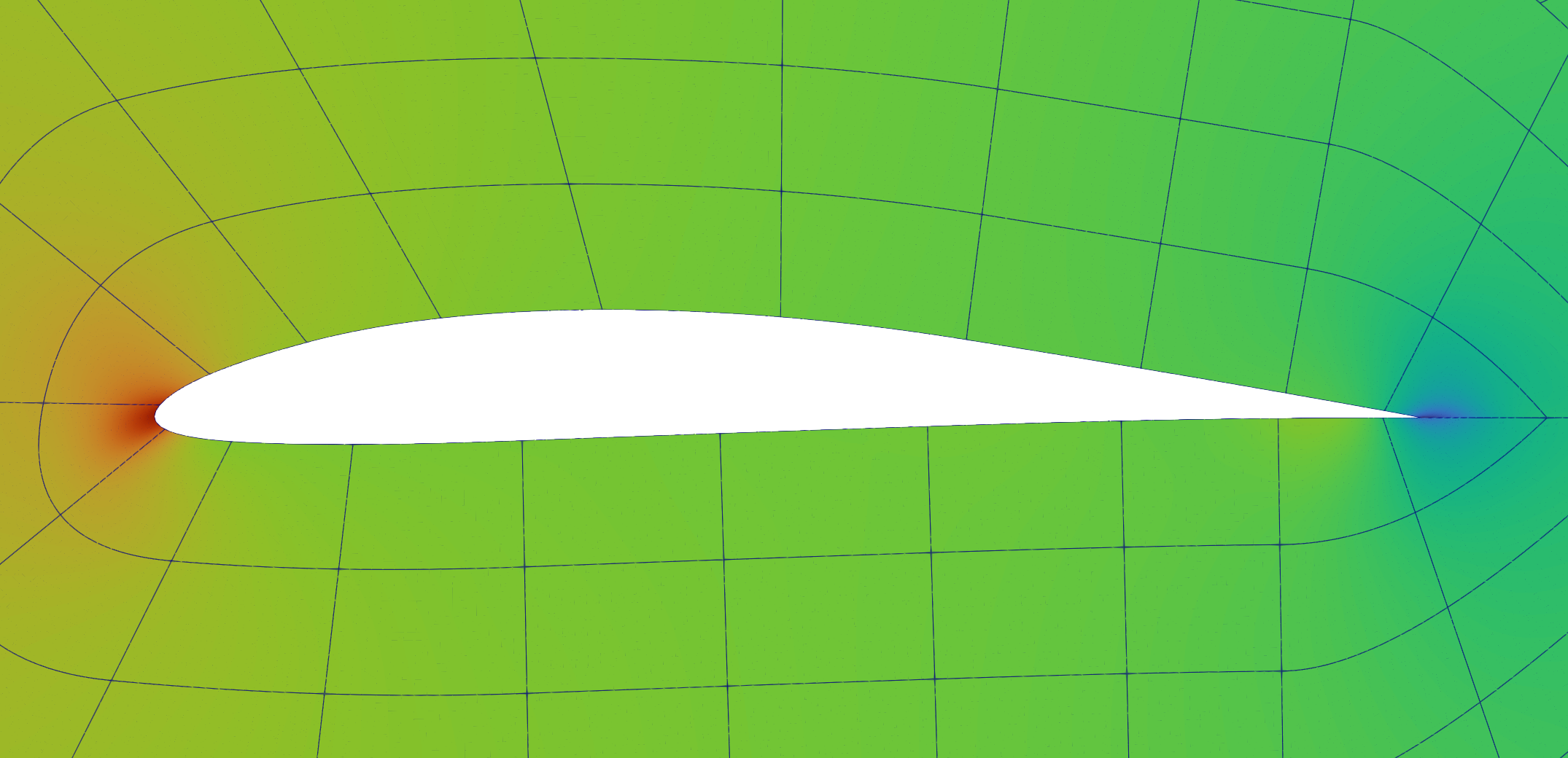} \\
(c) Trimming (refining 3 times) & (d) IBCM (initial mesh) \\
\end{tabular}
\caption{The pressure fields of the Stokes problem around the airfoil. (a, b) The results using the meshes in Figure~\ref{fig:airfoil_input}(b, c), respectively, (c) the zoom-in view of the pressure via trimming after refining the initial mesh three times, and (d) the zoom-in view of the pressure via IBCM using the initial mesh.}
\label{fig:airfoil_result}
\end{figure}

We first compare the pressure fields using the initial meshes shown in Figure~\ref{fig:airfoil_input}(b, c). The results are shown in Figure~\ref{fig:airfoil_result}(a, b). We observe that with trimming, as the mesh does not align with the airfoil geometry and the mesh resolution is not fine enough, the pressure at the two ends of the airfoil is poorly resolved. On the other hand, with IBCM, since its mesh aligns with the airfoil, it can already capture the large gradient even with a low mesh resolution. In fact, to achieve a comparable accuracy as IBCM (using the initial mesh), the trimming construction has to dyadically refine its initial mesh at least three times. The result is shown as a zoom-in view in Figure~\ref{fig:airfoil_result}(c). For easy comparison, the same zoom-in view of the IBCM result is shown in Figure~\ref{fig:airfoil_result}(d). 

\begin{figure}[!ht]
\centering
\includegraphics[width=0.6\textwidth]{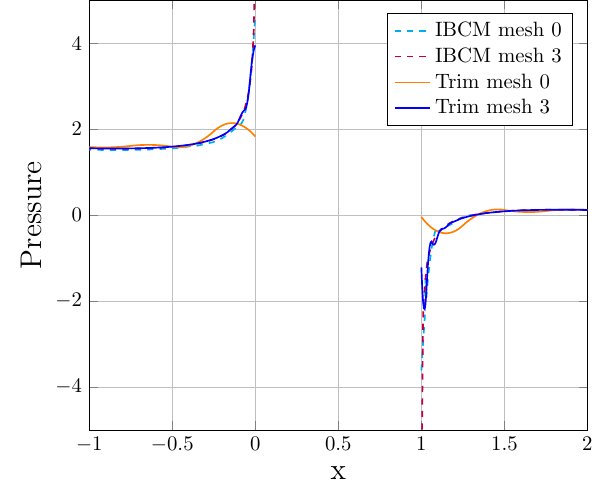}
\caption{The pressure plot along the line $y=0$, where the legends ``mesh 0'' and ``mesh 3'' mean the initial mesh and the mesh after three dyadic refinements, respectively. Note that the interval $(0,1)$ corresponds to the airfoil and only part of the entire domain is shown here.}
\label{fig:airfoil_plots}
\end{figure}

Moreover, in Figure~\ref{fig:airfoil_plots} we plot the pressure along the line $y=0$ using meshes of different resolutions. We observe that the results using IBCM converge fast (i.e., the results using ``mesh 0'' and ``mesh 3'' are very close), whereas the results using trimming vary a lot near the airfoil (i.e., around $x=0$ and $x=1$). In summary, with IBCM, we can leverage the geometric flexibility of the overlapping patches to efficiently capture the local features in the solution field without resorting to extensively refining the entire Cartesian grid.

\section{Conclusion}
\label{sec:con}

In this paper, we present a stabilized isogeometric formulation for the Stokes problem on overlapping multi-patches. The stabilization procedure is a generalization of \emph{minimal stabilization}, which, for the Stokes problem, treats the velocity and pressure separately. Stabilization of the velocity is only needed on interfaces, where flux terms from bad elements (i.e., cut elements that have small effective measures) are stabilized. On the other hand, the pressure space is modified on all bad elements and this provides a stable version of the Stokes problem on cut elements without any parameter to be set. Indeed, with one assumption (whose validity can be numerically assessed), we are able to prove the wellposedness and accuracy of our formulation in a completely general setting, allowing for overlaps and mesh cuts of any kind. 

In addition to the theoretical study, several numerical examples are presented to show the efficacy, accuracy, and efficiency of the proposed method. The expected convergence is achieved, which agrees with the theory. The accuracy of the pressure jump across interfaces and the conditioning of the linear system are significantly improved by the proposed stabilization procedure, which verifies the necessity of stabilization in practice. The efficiency and geometric flexibility of the proposed method are demonstrated with a special type of overlapping multi-patches that features a layer of conformal mesh around the geometry of interest (called the immersed boundary-conformal method, IBCM). It enables boundary-unfitted methods to naturally and efficiently capture the local solution features around domain boundaries.

In the future, it is promising to investigate minimal stabilization for Kirchhoff-Love shells, especially with the aid of IBCM for rapid geometric modeling and efficient solution procedure. Extension of minimal stabilization to the three-dimensional case via V-rep (the true volume representation as opposed to the dominant B-rep) is another promising and challenging direction. While the theory is mostly dimension-independent, it needs robust and efficient methods and algorithms to make it really happen.

\section*{Acknowledgment}

X. Wei is supported in part by the National Natural Science Foundation of China (No. 12202269). P. Antolin, and A. Buffa are partially supported by the SNSF of Switzerland through the project “Design- through-Analysis (of PDEs): the litmus test” n. 40B2-0 187094 (BRIDGE Discovery 2019), and the European Union’s Horizon 2020 research and innovation program under Grant Agreement n.862025 (ADAM2).



%% file: appendix.tex
\section{Technical proofs}\label{appendix:A1}
\begin{proof}[Proof of Theorem~\ref{thm:wellposedness}]
In order to prove~\eqref{wellposed:eq1}, let us show that there exist $c_1,c_2>0$ such that, for every $\left(\vv_h, q_h,\bmu_h \right)\in V_h\times\overline Q_h\times\Lambda_h$, there exists $\left(\w_h,r_h,\bm\eta_h \right)\in V_h\times\overline Q_h\times\Lambda_h$ such that $\overline{\mathcal A}_h\left( \left( \w_h,r_h,\bm\eta_h\right) ;\left(\vv_h,q_h,\bm\eta_h\right)\right) \ge c_1 \vertiii{\left(\vv_h,q_h,\bmu_h\right)}^2$ and $\vertiii{\left(\w_h,r_h,\bm\eta_h\right)}\le c_2 \vertiii{\left(\vv_h,q_h,\bmu_h\right)}$. Let us take $\left(\vv_h,q_h,\bmu_h\right)\in V_h\times\overline Q_h\times\Lambda_h$. It holds
\begin{align*}
&\overline{\mathcal A}_h \left( \left( \vv_h,-q_h,-\bmu_h \right) ;\left( \vv_h,q_h,\bmu_h\right) \right) = \sum_{i=0}^N\norm{D\vv_i}^2_{L^2(\Omega_i)}\\
&+\sum_{i=1}^N\sum_{j=0}^{i-1} \gamma^{-1} \left[ \norm{\h^{\frac{1}{2}}\left(-\bmu_h +\langle q_h\n \rangle_t\right)}^2_{L^2(\Gamma_{ij})} - \norm{\h^{\frac{1}{2}}\langle DR_{ij}^v(\vv_h)\n\rangle_t}^2_{L^2(\Gamma_{ij})}  \right] .
\end{align*}
Hence, by using Proposition~\ref{proposition_stability_velocity} and choosing $\gamma>0$ large enough,
\begin{equation}\label{appendix:eq:runliang1}
\begin{aligned}
\overline{\mathcal A}_h \left(\left( \vv_h,-q_h,-\bmu_h\right);\left(\vv_h,q_h,\bmu_h\right)\right)	
\ge&  \left(1-C\gamma^{-1}\right)\sum_{i=0}^N \norm{D\vv_i}^2_{L^2(\Omega_i)}
+ \gamma^{-1}\sum_{i=1}^N\sum_{j=0}^{i-1}\norm{\h^{\frac{1}{2}}\left(-\bmu_h+ \langle q_h\n\rangle_t\right)}^2_{L^2(\Gamma_{ij})} \\
\ge & \frac{1}{2}  \sum_{i=0}^N\norm{D\vv_i}^2_{L^2(\Omega_i)} +\gamma^{-1}\sum_{i=1}^N\sum_{j=0}^{i-1}\norm{\h^{\frac{1}{2}}\left(-\bmu_h+ \langle q_h\n\rangle_t\right)}^2_{L^2(\Gamma_{ij})}.
\end{aligned}	
\end{equation}
Now, let $\w_h^q\in V_h$ be the supremizer of Lemma~\ref{lemma:is}. Since $b_0(\w_h^q,q_h)\ge C\norm{q_h}^2_{0,h}$, we have
\begin{equation}\label{wellposed:eq3}
\begin{aligned}
\overline{\mathcal A}_h\left(\left(\w_h^q,0,\bm 0\right);\left(\vv_h,q_h,\bmu_h\right)\right)
\ge & \sum_{i=0}^N\int_{\Omega_i} D\vv_i:D\w_i^p  +C  \norm{q_h}^2_{0,h} + \sum_{i=1}^N \sum_{j=0}^{i-1}\int_{\Gamma_{ij}} \bmu_h [\w_h^q]\\
& - \sum_{i=1}^N \sum_{j=0}^{i-1}\gamma^{-1} \int_{\Gamma_{ij}}\h \langle DR_{ij}^v(\vv_h)\n\rangle_t \langle DR_{ij}^v(\w_h^q)\n\rangle_t\\
&+  \sum_{i=1}^N \sum_{j=0}^{i-1} \gamma^{-1} \int_{\Gamma_{ij}}\h \left( -\bmu_h + \langle q_h\n\rangle_t\right) \langle DR_{ij}^v(\w_h^q)\n\rangle_{t}.
\end{aligned}	
\end{equation}
By the Cauchy-Schwarz inequality
\begin{align*}
\sum_{i=1}^N \sum_{j=0}^{i-1} &	\int_{\Gamma_{ij}} \h \langle DR_{ij}^v(\vv_h)\n\rangle_t \langle DR^v_h(\w_h^q)\n\rangle_t\\
\le  & \sum_{i=1}^N \sum_{j=0}^{i-1} \norm{ \h^{\frac{1}{2}}\langle DR_{ij}^v(\vv_h)\n\rangle_t}_{L^2(\Gamma_{ij})}\norm{\h^{\frac{1}{2}}\langle DR_{ij}^v(\w_h^q)\n\rangle_t}_{L^2(\Gamma_{ij})},\\
\sum_{i=1}^N \sum_{j=0}^{i-1} &	\abs{ \int_{\Gamma_{ij}}\h \left(-\bmu_h+\langle q_h\n\rangle_t\right) \langle DR_{ij}^v(\w_h^q)\n  \rangle_{t}}\\
\le&\sum_{i=1}^N \sum_{j=0}^{i-1}  \norm{\h^{\frac{1}{2}} \left(-\bmu_h +\langle q_h\n\rangle_t \right)}_{L^2(\Gamma_{ij})}
\norm{\h^{\frac{1}{2}}\langle DR_{ij}^v(\w_h^q)\n\rangle_t}_{L^2(\Gamma_{ij})}.
\end{align*}	
By using Young's inequality, for $s,\delta>0$, Proposition~\ref{proposition_stability_velocity} and the construction of $\w_h^q$ in the proof of Lemma~\ref{lemma:is}, there exists $C_1>0$ such that
\begin{align*}
\sum_{i=1}^N \sum_{j=0}^{i-1}\gamma^{-1}\int_{\Gamma_{ij}} \h\langle DR_{ij}^v(\vv_h)\n\rangle_t \langle DR^v_{ij}(\w_h^q)\n\rangle_t  	\ge & -\frac{C_1}{2\gamma s} \sum_{i=0}^N\norm{D\vv_i}^2_{L^2(\Omega_i)} - \frac{sC_1}{2\gamma}\norm{q_h}^2_{0,h}, \\
\sum_{i=1}^N \sum_{j=0}^{i-1} \gamma^{-1}\int_{\Gamma_{ij}} \h \left(-\bmu_h+\langle q_h\n\rangle_t\right) \langle DR_{ij}^v(\w_h^q)\n  \rangle_{t} \ge & -\frac{C_1}{2\gamma\delta}\sum_{i=1}^N \sum_{j=0}^{i-1} \norm{\h^{\frac{1}{2}}\left(-\bmu_h+\langle q_h\n\rangle_t\right))}_{L^2(\Gamma_{ij})}^2\\
& -\frac{\delta C_1}{2\gamma} \norm{q_h}^2_{0,h}.
\end{align*}	
From the Cauchy-Schwarz inequality, the construction of $\w_h^q$, and Young's inequality, there exists $C_2>0$ such that, for $\eps>0$,
\begin{align*}
\sum_{i=0}^N \int_{\Omega_i}D\vv_i:D\w_i^p \ge -\frac{C_2}{2\eps} \sum_{i=0}^N \norm{D\vv_i}^2_{L^2(\Omega_i)} - \frac{\eps C_2}{2} \norm{q_h}^2_{0,h}.
\end{align*}	
In an analogous fashion, there exists $C_3>0$ such that, for $r>0$,
\begin{align*}
\sum_{i=1}^N \sum_{j=0}^{i-1}	\int_{\Gamma_{ij}} \bmu_h [\w_h^q] \ge -\frac{C_3}{2r} \sum_{i=1}^N \sum_{j=0}^{i-1} \norm{\mathsf h^{\frac{1}{2}}\bmu_h}^2_{L^2(\Gamma_{ij})} - \frac{rC_3 }{2} \norm{q_h}^2_{0,h}.
\end{align*}
Let us go back to~\eqref{wellposed:eq3}. We have
\begin{align*}
\overline{\mathcal A}_h& \left( \left(\w_h^q,0,\bm 0 \right)  ; \left(\vv_h,q_h,\bmu_h\right)\right)\ge -\frac{C_2}{2\eps} \sum_{i=0}^N\norm{D\vv_i}^2_{L^2(\Omega_i)}-\frac{\eps C_2}{2}\norm{q_h}^2_{0,h} +C\norm{q_h}^2_{0,h}\\
&-\frac{C_3}{2r}\sum_{i=1}^N \sum_{j=0}^{i-1}\norm{\h^{\frac{1}{2}}\bmu_h}^2_{L^2(\Gamma_{ij})}
- \frac{r C_3}{2}\norm{q_h}^2_{0,h}
-\frac{C_1 }{2\gamma s}\sum_{i=0}^N\norm{D\vv_i}^2_{L^2(\Omega_i)}\\
& - \frac{s C_1}{2\gamma}\norm{q_h}^2_{0,h}
- \frac{C_1}{2\gamma\delta}\sum_{i=1}^N\sum_{j=0}^{i-1} \norm{\h^{\frac{1}{2}} \left(-\bmu_h+\langle q_h\n\rangle_t \right)}_{L^2(\Gamma_{ij})} - \frac{\delta C_1 }{2\gamma}\norm{q_h}^2_{0,h}\\
= &\left( C-\frac{\eps C_2}{2}-\frac{r C_3}{2}-\frac{s C_1}{2\gamma} - \frac{\delta C_1 }{2\gamma}\right)\norm{q_h}^2_{0,h} + \left( -\frac{C_2}{2\eps}-\frac{C_1 }{2\gamma s}\right)\sum_{i=0}^N \norm{D\vv_i}^2_{L^2(\Omega_i)}\\
& - \frac{C_1}{2\gamma\delta} \sum_{i=1}^N\sum_{j=0}^{i-1}\norm{\h^{\frac{1}{2}}\left(-\bmu_h+\langle q_h\n\rangle \right)}^2_{L^2(\Gamma_{ij})} - \frac{C_3}{2r} \sum_{i=1}^N\sum_{j=0}^{i-1}\norm{\h^{\frac{1}{2}} \bmu_h}^2_{L^2(\Gamma_{ij})}.
\end{align*}
Let $\eps,r,s,\delta>0$ be small enough such that $\displaystyle C-\frac{\eps C_2}{2}-\frac{r C_3}{2}-\frac{s C_1}{2\gamma} - \frac{\delta C_1 }{2\gamma}\ge \frac{C}{2}$. Hence, there exist $C_4, C_5, C_6 >0$ such that
\begin{equation}\label{appendix:eq:runliang2}
\begin{aligned}
\overline{\mathcal A}_h&\left(\left(\w_h^q,0,\bm 0\right);\left(\vv_h,q_h,\bmu_h \right)\right)\ge \frac{C}{2} \norm{q_h}^2_{0,h}- C_4 \sum_{i=0}^N \norm{D\vv_i}^2_{L^2(\Omega_i)}\\
& -  C_5 \sum_{i=1}^N \sum_{j=0}^{i-1} \norm{\h^{\frac{1}{2}}\left( -\bmu_h+\langle q_h\n\rangle_t\right)}^2_{L^2(\Gamma_{ij})}
-C_6\sum_{i=1}^N \sum_{j=0}^{i-1}\norm{\h^{\frac{1}{2}}\bmu_h}^2_{L^2(\Gamma_{ij})}.
\end{aligned}
\end{equation}
Let $P_h:\bigoplus_{0\le j<i\le N} L^2(\Gamma_{ij})\to \Lambda_h$ be the $L^2$-orthogonal projection. From condition~\eqref{condition:multiplier_space}, we have $\h^{-1}\restr{[\vv_h]}{\Gamma_{ij}}=\restr{P_h \h^{-1}[\vv_h]}{\Gamma_{ij}}$, for every $0\le j<i\le N$. Therefore, from the Cauchy-Schwartz inequality, we have
\begin{align*}
\overline{\mathcal A}_h  \left( \left( \bm 0,0, P_h\h^{-1}[\vv_h] \right) ;\left(\vv_h,q_h,\bmu_h\right)\right)
\ge & \sum_{i=1}^N\sum_{j=0}^{i-1}  \norm{\h^{-\frac{1}{2}}[\vv_h]}^2_{L^2(\Gamma_{ij})}\\
&-\gamma^{-1}\bigg( \sum_{i=1}^N\sum_{j=0}^{i-1} \norm{\h^{\frac{1}{2}}\langle DR_{ij}^v(\vv_h)\n\rangle_t}_{L^2(\Gamma_{ij})} \norm{\h^{-\frac{1}{2}}[\vv_h]}_{L^2(\Gamma_{ij})}\\
& + \sum_{i=1}^N\sum_{j=0}^{i-1} \norm{\h^{\frac{1}{2}}\left( -\bmu_h+\langle q_h\n\rangle_t \right)}_{L^2(\Gamma_{ij})} \norm{\h^{-\frac{1}{2}}[\vv_h]}_{L^2(\Gamma_{ij})}\bigg).
\end{align*}
By Proposition~\ref{proposition_stability_velocity} and Young's inequality, there exists $C_7>0$ such that, for $a,b>0$,
\begin{align*}
- \sum_{i=1}^N\sum_{j=0}^{i-1} \norm{\h^{\frac{1}{2}}\langle  DR^v_h (\vv_h)\n\rangle_t}_{L^2(\Gamma_{ij})}\norm{\h^{-\frac{1}{2}}[\vv_h]}_{L^2(\Gamma_{ij})} \ge& -\frac{C_7 }{2a} \sum_{i=0}^N \norm{D \vv_i}^2_{L^2(\Omega_i)}\\
& - \frac{aC_7}{2} \sum_{i=1}^N\sum_{j=0}^{i-1} \norm{\h^{-\frac{1}{2}}[\vv_h]}^2_{L^2(\Gamma_{ij})},\\
- \sum_{i=1}^N\sum_{j=0}^{i-1} \norm{\h^{\frac{1}{2}}\left(-\bmu_h+\langle q_h \n \rangle_t\right)}_{L^2(\Gamma_{ij})}\norm{\h^{-\frac{1}{2}}[\vv_h]}_{L^2(\Gamma_{ij})}\ge& -\frac{1}{2b}  \sum_{i=1}^N\sum_{j=0}^{i-1} \norm{\h^{\frac{1}{2}}\left( -\bmu_h+\langle q_h\n\rangle_t\right)}^2_{L^2(\Gamma_{ij})}\\
&- \frac{b}{2}  \sum_{i=1}^N\sum_{j=0}^{i-1}\norm{\h^{-\frac{1}{2}}[\vv_h]}^2_{L^2(\Gamma_{ij})}.
\end{align*}
Thus,
\begin{equation*}
\begin{aligned}
\overline{\mathcal A}_h&\left( \left(\bm 0,0P_h\h^{-1}[\vv_h]\right);\left(\vv_h,q_h,\bmu_h\right)\right) \ge \left( 1-\frac{aC_7}{2\gamma}-\frac{b}{2\gamma}\right)  \sum_{i=1}^N\sum_{j=0}^{i-1} \norm{\h^{-\frac{1}{2}}[\vv_h]}^2_{L^2(\Gamma_{ij})}\\
& - \frac{C_7}{2a\gamma} \sum_{i=0}^N\norm{D\vv_i}^2_{L^2(\Omega_i)}-\frac{1}{2b\gamma}  \sum_{i=1}^N\sum_{j=0}^{i-1} \norm{\h^{\frac{1}{2}}\left(-\bmu_h+\langle q_h\n\rangle_t \right)}^2_{L^2(\Gamma_{ij})}.
\end{aligned}
\end{equation*}
Let us choose $a,b>0$ small enough such that $\displaystyle 1-\frac{aC_7}{2\gamma}-\frac{b}{2\gamma}\ge \frac{1}{2}$. Hence, there exist $C_8,C_9>0$ such that
\begin{equation}\label{appendix:eq:runliang3}
\begin{aligned}
\overline{\mathcal A}_h\left( \left(\bm 0,0,P_h\h^{-1}[\vv_h]\right);\left(\vv_h,q_h,\bmu_h\right)\right) \ge & \frac{1}{2}  \sum_{i=1}^N\sum_{j=0}^{i-1} \norm{\h^{-\frac{1}{2}}[\vv_h]}^2_{L^2(\Gamma_{ij})}
-C_8 \sum_{i=0}^N\norm{D\vv_i}^2_{L^2(\Omega_i)}\\
&-C_9  \sum_{i=1}^N\sum_{j=0}^{i-1} \norm{\h^{\frac{1}{2}}\left(-\bmu_h+\langle q_h\n\rangle_t \right)}^2_{L^2(\Gamma_{ij})}.
\end{aligned}
\end{equation}	
Let us put together~\eqref{appendix:eq:runliang1},~\eqref{appendix:eq:runliang2},~\eqref{appendix:eq:runliang3}. For $k,\eta>0$, we have
\begin{align*}
\overline{\mathcal A}_h&\left(\left( \vv_h + k \w^q_h,-q_h,-\bm\mu_h+\eta P_h \h^{-1} [\vv_h] \right) ;\left( \vv_h,q_h,\bmu_h\right)\right) 
\ge  \left( \frac{1}{2}-kC_4 -\eta C_8\right) \sum_{i=0}^N \norm{D\vv_i}^2_{L^2(\Omega_i)}\\
&+ \left(\frac{1}{\gamma}-kC_5-\eta C_9 \right) \sum_{i=1}^N \sum_{j=0}^{i-1} \norm{\h^{\frac{1}{2}}\left(-\bm\mu_h +\langle q_h\n\rangle_t\right)}^2_{L^2(\Gamma_{ij})} 
+ \frac{kC}{2} \norm{q_h}^2_{0,h}\\
& -kC_6 \sum_{i=1}^N \sum_{j=0}^{i-1}\norm{\h^{\frac{1}{2}}\bm\mu_h}^2_{L^2(\Gamma_{ij})} + \frac{\eta}{2} \sum_{i=1}^N \sum_{j=0}^{i-1} \norm{\h^{-\frac{1}{2}} [\vv_h]}^2_{L^2(\Gamma_{ij})} .
\end{align*}
From Proposition~\ref{prop_press_stability}, there exists $C_{10}>0$ such that
\begin{equation}\label{appendix:hongkong}
\begin{aligned}
\overline{\mathcal A}_h&\left(\left( \vv_h + k \w^q_h,-q_h,-\bm\mu_h+\eta P_h \h^{-1} [\vv_h] \right) ;\left( \vv_h,q_h,\bmu_h\right)\right) 
\ge  \left( \frac{1}{2}-kC_4 -\eta C_8\right) \sum_{i=0}^N \norm{D\vv_i}^2_{L^2(\Omega_i)}\\
&+ \left(\frac{1}{\gamma}-kC_5-\eta C_9 \right) \sum_{i=1}^N \sum_{j=0}^{i-1} \norm{\h^{\frac{1}{2}}\left(-\bm\mu_h +\langle q_h\n\rangle_t\right)}^2_{L^2(\Gamma_{ij})} 
+\frac{k}{4C_{10} }\sum_{i=1}^N\sum_{j=0}^{i-1} \norm{\h^{\frac{1}{2}}\langle q_h\rangle_t}^2_{L^2(\Gamma_{ij})}\\
&+ \frac{kC}{4} \norm{q_h}^2_{0,h} -kC_6 \sum_{i=1}^N \sum_{j=0}^{i-1}\norm{\h^{\frac{1}{2}}\bm\mu_h}^2_{L^2(\Gamma_{ij})} + \frac{\eta}{2} \sum_{i=1}^N \sum_{j=0}^{i-1} \norm{\h^{-\frac{1}{2}} [\vv_h]}^2_{L^2(\Gamma_{ij})} .
\end{aligned}
\end{equation}
Let $\displaystyle C_{11}:= \frac{4}{k C_{10}}$ and $\displaystyle C_{12}:=\frac{1}{\gamma}-kC_5-\eta C_9$, and observe that, by Young inequality, for any $\ell >0$,
\begin{equation*}
\norm{\h^{\frac{1}{2}}\left(-\bm\mu_h +\langle q_h\n\rangle_t\right)}^2_{L^2(\Gamma_{ij})} \ge 
(1-\ell) \norm{\h^{\frac{1}{2}}\langle q_h \n\rangle_t}^2_{L^2(\Gamma_{ij})} + (1-\frac{1}{\ell}) \norm{ \h^{\frac{1}{2}} \bmu_h}^2_{L^2(\Gamma_{ij})}.
\end{equation*}
By using $\norm{\h^{\frac{1}{2}}\langle q_h \n\rangle_t}^2_{L^2(\Gamma_{ij})} = \norm{\h^{\frac{1}{2}}\langle q_h \rangle_t}^2_{L^2(\Gamma_{ij})}$, we have
\begin{equation}\label{appendix:stress}
\begin{aligned}
&C_{11}\sum_{i=1}^N\sum_{j=0}^{i-1} \norm{\h^{\frac{1}{2}}\langle q_h\rangle_t}^2_{L^2(\Gamma_{ij})} + C_{12} \sum_{i=1}^N \sum_{j=0}^{i-1} \norm{\h^{\frac{1}{2}}\left(-\bm\mu_h +\langle q_h\n\rangle_t\right)}^2_{L^2(\Gamma_{ij})}\\
\ge & C_{12} \sum_{i=1}^N\sum_{j=0}^{i-1} \left(  \left(  \frac{C_{11}}{C_{12}} + 1 -\ell   \right) \norm{\h^{\frac{1}{2}} \langle q_h\rangle_t }^2_{L^2(\Gamma_{ij})}  + \left( 1-\frac{1}{\ell}\right) \norm{\h^{\frac{1}{2}} \bmu_h}^2_{L^2(\Gamma_{ij})}  \right).
\end{aligned}
\end{equation}
By plugging~\eqref{appendix:stress} back into~\eqref{appendix:hongkong}, we have
\begin{align*}
\overline{\mathcal A}_h&\left(\left( \vv_h + k \w^q_h,-q_h,-\bm\mu_h+\eta P_h \h^{-1} [\vv_h] \right) ;\left( \vv_h,q_h,\bmu_h\right)\right) \ge  \left( \frac{1}{2}-kC_4 -\eta C_8\right) \sum_{i=0}^N \norm{D\vv_i}^2_{L^2(\Omega_i)}\\
+& C_{12} \left( \frac{C_{11}}{C_{12}}+1 -\ell \right) \sum_{i=1}^N \sum_{j=0}^{i-1} \norm{\h^{\frac{1}{2}} \langle q_h  \rangle_t}^2_{L^2(\Gamma_{ij})} + \left( C_{12}\left(1-\frac{1}{\ell}\right) -k C_6\right) \sum_{i=1}^N\sum_{j=0}^{i-1} \norm{\h^{\frac{1}{2}} \bmu_h}^2_{L^2(\Gamma_{ij})} \\
+& \frac{kC}{4}\norm{q_h}^2_{0,h} + \frac{\eta}{2} \sum_{i=1}^N \sum_{j=0}^{i-1} \norm{\h^{-\frac{1}{2}} [\vv_h]}^2_{L^2(\Gamma_{ij})} .
\end{align*}
We require $\displaystyle 1<\ell < \frac{C_{11}}{C_{12}}+1$ and $k,\eta$ to be small enough so that $\displaystyle \frac{1}{2}-kC_4 -\eta C_8\ge \frac{1}{4}$, $\displaystyle C_{12} \left( \frac{C_{11}}{C_{12}}+1 -\ell \right) \ge C_{13}$, $\displaystyle C_{12}\left(1-\frac{1}{\ell}\right) -k C_6\ge C_{14}$, for some $C_{13}, C_{14}>0$. Note that the choice of $\ell$ depends on $k$ and $\eta$. On the other hand, there exist $C_{15},C_{16}>0$ such that $\displaystyle \frac{kC}{4}\ge C_{15}$ and $\displaystyle \frac{\eta}{2} \ge C_{16}$. Hence,
\begin{align*}
\overline{\mathcal A}_h&\left(\left( \vv_h + k \w^q_h,-q_h,-\bm\mu_h+\eta P_h \h^{-1} [\vv_h] \right) ;\left( \vv_h,q_h,\bmu_h\right)\right) \ge  \frac{1}{4} \sum_{i=0}^N \norm{D\vv_i}^2_{L^2(\Omega_i)}\\
+& C_{13} \sum_{i=1}^N \sum_{j=0}^{i-1} \norm{\h^{\frac{1}{2}} \langle q_h  \rangle_t}^2_{L^2(\Gamma_{ij})} + C_{14} \sum_{i=1}^N\sum_{j=0}^{i-1} \norm{\h^{\frac{1}{2}} \bmu_h}^2_{L^2(\Gamma_{ij})}\\
+& C_{15}\norm{q_h}^2_{0,h} +C_{16} \sum_{i=1}^N \sum_{j=0}^{i-1} \norm{\h^{-\frac{1}{2}} [\vv_h]}^2_{L^2(\Gamma_{ij})} .
\end{align*}
Finally, the stability property of $\w_h^q$, namely $\norm{\w_h^q}_{1,h}\le C \norm{q_h}_{0,h}$, entails
\begin{align*}
\vertiii{\left(\vv_h+k\w_h^q,-q_h,-\bm\mu_h+\eta P_h \h^{-1}[\vv_h]\right)}\le C  \vertiii{\left(\vv_h,q_h,\bm\mu_h\right)}.
\end{align*}
\end{proof}
\IfStandalone
{
\bibliographystyle{../elsarticle-num}
\bibliography{../bibliographyx}
}{}